\documentclass[a4paper,11pt,oneside,reqno]{amsart}
\usepackage{amsmath}
\usepackage{amsthm}
\usepackage{amssymb}
\usepackage{graphicx}
\theoremstyle{plain}
\newtheorem{teore}{Theorem}[section]
\newtheorem*{teore*}{Theorem}
\newtheorem*{coro*}{Corollary}
\newtheorem{defin}[teore]{Definition}
\newtheorem{lem}[teore]{Lemma}
\newtheorem{coro}[teore]{Corollary}
\newtheorem{propo}[teore]{Proposition}
\newtheorem*{claim}{Claim}
\newtheorem*{claim*}{Claim}
\theoremstyle{remark}
\newtheorem{ejemplo}[teore]{{\sc Example}}
\newtheorem{notas}[teore]{{\sc Remark}}
\newcommand{\nrm}[1]{\|#1\|}
\DeclareMathOperator{\conv}{conv} \DeclareMathOperator{\sgn}{sgn}

\newcommand{\prop}{\begin{propo}}
\newcommand{\fprop}{\end{propo}}
\newcommand{\cor}{\begin{coro}}
\newcommand{\fcor}{\end{coro}}

\newcommand{\defi}{\begin{defin}\rm}
\newcommand{\fdefi}{\end{defin}}
\newcommand{\eje}{\begin{ejemplo}}
\newcommand{\feje}{\end{ejemplo}}
\newcommand{\lema}{\begin{lem}}
\newcommand{\flema}{\end{lem}}
\newcommand{\teor}{\begin{teore}}
\newcommand{\fteor}{\end{teore}}
\newcommand{\nota}{\begin{notas}\rm}
\newcommand{\fnota}{ \end{notas}}
\newcommand{\clam}{\begin{claim}}
\newcommand{\fclam}{\end{claim}}
\newcommand{\clams}{\begin{claim*}}
\newcommand{\fclams}{\end{claim*}}

\newcommand{\lclam}{\begin{lclaim}}
\newcommand{\flclam}{\end{lclaim}}
\newcommand{\prucl}{\prue[Proof of Claim:]}
\newcommand{\fprucl}{\fprue}
\newcommand{\ben}{\begin{enumerate}}
\newcommand{\een}{\end{enumerate}}
\newcommand{\bit}{\begin{itemize}}
\newcommand{\eit}{\end{itemize}}
\newcommand{\cl}[2]{\overline{#1}^{#2}}
\newcommand{\mc}[1]{\mathcal{#1}}
\newcommand{\mm}[1]{\mathrm{#1}}

\newcommand{\casos}{\begin{itemize}}
\newcommand{\fcasos}{\end{itemize}\setcounter{cs}{1}}

\newcommand{\ptot}{\forall}

\newcommand{\pe}{\preceq}

\newcommand{\rou}{\al\ga\le \al\be\vee \be\ga}
\newcommand{\rod}{\al\be\le \al\ga\vee \be\ga}

\newcommand{\ro}{\varrho}

\newcommand{\conj}[2]{ \{ {#1}\,:\,{#2} \} }

\newcommand{\ou}{\omega_{1}}

\newcommand{\om}{\omega}
\newcommand{\Om}{\Omega}

\newcommand{\ka}{\kappa}

\newcommand{\buit}{\emptyset}

\newcommand{\ga}{\gamma}

\newcommand{\ko}{\ensuremath{K_{\omega_1}}}
\newcommand{\co}{\ensuremath{c_{00}(\omega_1)}}

\renewcommand{\labelenumi}{\theenumi.}

\newcommand{\Ga}{\Gamma}

\newcommand{\al}{\alpha}
\newcommand{\be}{\beta}
\newcommand{\de}{\delta}

\newcommand{\la}{\lambda}

\newcommand{\La}{\Lambda}
\renewcommand{\P}{\mathfrak{P}}
\newcommand{\sig}{\sigma}

\newcommand{\vep}{\varepsilon}

\newcommand{\eqs}{{\mathfrak X}}
\newcommand{\R}{{\mathbb R}}

\newcommand{\N}{{\mathbb N}}

\newcommand{\Q}{{\mathbb Q}}

\renewcommand{\P}{P}

\renewcommand{\ker}{\mathrm{Ker}}

\newcommand{\supp}{\mathrm{supp\, }}
\newcommand{\range}{\mathrm{ran\, }}
\newcommand{\ran}{\mathrm{ran\, }}

\newcommand{\oo}{\infty}

\newcommand{\con}{\subseteq}

\newcommand{\cones}{\varsubsetneq }

\newcommand{\prue}{\begin{proof}}
\newcommand{\fprue}{\end{proof}}

\makeindex

\begin{document}
\title{A class of
 Banach spaces with few non strictly singular operators}
\author{S. A. Argyros}
\address{Department of Mathematics, National Technical University of Athens,   15780 Athens, Greece}
\email{sargyros@math.ntua.gr}
\author{$\text{J. Lopez-Abad}^*$ }
\address{Equipe de Logique Math\'{e}matique, Universit\'{e} Paris VII, 2 place Jussieu, 75251 Paris Cedex 05, France}
\email{abad@logique.jussieu.fr}
\author{S. Todorcevic}
\address{CNRS-Universit\'{e} Paris VII, 2, place Jussieu, 75251 Paris Cedex 05, France}
\email{stevo@math.jussieu.fr} %

\thanks{ $^{*}$\text{Supported by a European Community Marie Curie Fellowship}}

\maketitle


The original motivation for this paper is based on  the natural question left open by the
Gowers-Maurey solution of the unconditional basic sequence problem for Banach spaces
(\cite{go-ma}). Recall that Gowers and Maurey have constructed a Banach space $X$ with a
Schauder basis $(e_n)_n$ but with no unconditional basic sequence. Thus, while every infinite
dimensional Banach space contains a sequence $(x_n)_n$ which forms a Schauder basis for its
closure $Y=\overline{\langle x_n \rangle_n}$,  meaning that every vector of $Y$ has a unique
representation $\sum_{n}a_n x_n$, one may not be able to get such $(x_n)_n$ such that the sums
$\sum_{n}a_n x_n$ converge unconditionally whenever they converge. The fundamental role of
Schauder basis and the fact that the notion is very much dependent on the order  lead to the
natural variation of the notion, the definition of transfinite Schauder basis
$(x_\al)_{\al<\ga}$, where  vectors of $X$ have a unique representations as  sums
$\sum_{\al<\ga}a_\al x_\al$. In fact, as it will be clear from some results in this paper,
considering a transfinite Schauder basis, even if one knows that $X$ has an ordinary Schauder
basis, can be an advantage. Thus, the natural question which originated the research of this
paper asks whether one can have Banach spaces with long (even of uncountable length) Schauder
bases but with no unconditional basic sequence. There is actually a more fundamental reason
for asking this question. As noticed originally by W. B. Johnson, the Gowers-Maurey space $X$
is hereditarily indecomposable which in particular yields that the space of operators on $X$
is very small in the sense that every bounded linear operator on $X$ can be written as $\la
\mathrm{Id}_X +S$, where $S$ is a strictly singular operator. On the other hand, if $X$ has a
transfinite Schauder basis $(e_\al)_{\al<\ga}$ of length, say, $\ga=\om^2$, it could no longer
have so small an operator space as projections on infinite intervals $\overline{\langle e_\al
\rangle_{\al\in I}}$ are all (uniformly) bounded. Thus one would like to find out the amount
of control on the space of non strictly singular operators that is possible in this case. In
fact, our solution of the transfinite variation of the unconditional basic sequence problem
has led us to many other new questions of this sort, has forced us to introduce several new
methods to this area, and has revealed several new phenomena that could have been perhaps
difficult to discover by working only in the context of ordinary Schauder bases.

To see the necessity for a new method we repeat that our first goal here is to construct a
Banach space $\eqs_{\ou}$ with a transfinite Schauder basis $(e_\al)_{\al<\ou}$ with no
unconditional basic sequence as well as to understand its separable initial segments
$\eqs_{\ga}=\overline{\langle e_\al\rangle_{\al<\ga}}.$ The original Gowers-Maurey method for
preventing unconditional basic sequences is to force the unconditional constants of initial
finite-dimensional  subspaces, according to the fixed Schauder basis, grow to infinity. Since
initial finite-dimensional subspaces according to our transfinite Schauder basis
$(e_\al)_{\al<\ou}$ are far from exhausting the whole space their method will not work here.
It turns out that in order to impose the conditional structure on our space(s) $\eqs_{\ga}$
($\ga\le \ou$) we needed to import a tool from another area of mathematics, a rather canonical
semi-distance function $\ro$ on the space $\ou$ of all countable ordinals (\cite{tod1}). What
$\ro$ does in our context here is to essentially identify the structure of finite-dimensional
subspaces of various $\eqs_{\ga}$'s which globally are of course very much different, since
for example $\eqs_\om$ is hereditarily indecomposable while, say, $\eqs_{\om^2}$ has a rich
space of non-strictly singular operators.

After solving this initial problem we went on and tried to show
that every bounded linear operator
$T$ on a given $\eqs_\ga$ is a sum of a
diagonal operator $D_T$ and a strictly singular one. There
are natural candidates for $D_T$ which would share the eigenvalues of $T$
 and have the property
that  $T-D_T$ is strictly singular. The problem is to show that $D_T$ is a
 bounded operator. This
forced us to a variation on the notion of $\ro$-function by adding to it certain universality
property. To see the need for this, suppose we are given a finitely supported vector $x$ such
that $\nrm{D_Tx}$ is very large in comparison with $\nrm{Tx}$. The vector $x$ has a natural
decomposition $x=x_1+\dots+x_n$ such that $D_Tx=\la_1x_1+\dots +\la_nx_n$ where $\la_i$'s are
eigenvalues of $T$. The universality of $\ro$ guarantees that $x_i$'s can all be
simultaneously moved (keeping the discrepancy between $\nrm{D_T x}$ and $\nrm{Tx}$) to be
almost equal to eigenvectors with eigenvalues $\la_i$'s giving us an impossibility. This also
gives us a new phenomenon, unprecedented in this area, that every finite dimensional subspace
of some $\eqs_{\ga}$ can be moved by an $(4+\vep)$-isomorphism to essentially any  region of
any other $\eqs_\de$.

Our attempt to extend the control of operators to arbitrary subspaces of $\eqs_{\ou}$ has
led us to a new phenomenon which a priori could have been discovered before since it
already has a solid basis in an old paper of Maurey-Rosenthal (\cite{mau-ros}). What we
discovered is that each $\eqs_\ga$ has an associated James-like space $J_{T_0}$ which is
minimally and canonically finitely block represented in $\eqs_{\ga}$ and which is
responsible for essentially all of its conditional and unconditional geometry, including
the complete structure of the corresponding space of bounded non-strictly singular
operators. In retrospect, what Maurey-Rosenthal \cite{mau-ros} have done in their attempt
to solve the unconditional basic sequence problem is to produce a space $X$ with a
Schauder basis $(e_n)_n$ such that every subsequence $(e_{n_k})_k$ finitely block
represents $J_{c_0}$, a fact which then they used to show that no subsequence of
$(e_n)_n$ is unconditional. The finite representability of $J_{T_0}$ and the global
control of block sequences provided by $\ro$ gives us a complete picture of the space of
bounded non-strictly singular operators defined not only on $\eqs_{\ga}$ ($\ga\le \ou$)
but also on their arbitrary subspaces. For example, we show that the space of all bounded
non-strictly singular operators on a given $\eqs_\ga$ is naturally isomorphic to the dual
of the corresponding James-like space $J_{T_0}$. We also discover subspaces $X$ of
$\eqs_\ga$ such that the non-strictly singular part of the operator space $\mc
L(X,\eqs_\ga)$ is quite rich but on the other hand every bounded operator $T:X\to X$ is a
strictly singular perturbation of a scalar multiple of the identity. Another new
phenomenon we found are hereditarily indecomposable subspaces of $\eqs_{\ou}$ that are
asymptotic versions of themselves.

We now pass to a more detailed presentation of the specific results of this paper. The first
section concerns extensions of some standard facts about Schauder basic sequences to the
transfinite case. For example we show that every subspace $Y$ of a space $X$ with a
transfinite basis contains a further subspace $Z$ isomorphic to a block subspace of $X$. We
should point out that this result is weaker than  the corresponding results  for Schauder
bases. This causes some problems when one tries to extend standard constructions into the
transfinite case. For example, if one considers the transfinite version of the Schlumprecht
space, or more generally spaces of the form $T_{\ga}[(1/m_j,n_j)_j],$ one runs into
difficulties when trying to prove arbitrarily distortion. We overcome this by adding a
property to the basis $(x_{\al})_{\al<\ga}$ which ensures that the block sequences approximate
in a strong sense the subspaces of the space $\eqs _{\omega_1}$. This condition permits us to
show that the spaces $\eqs_{\omega_1}$ and $\eqs_{\omega_1}^u$ are arbitrarily distortable. We
also give a characterization of  reflexivity analogous to the classical one due to James
\cite{jam}.

The second section is mainly devoted to the definition of the norming set $K_{\ou}$ of the
maximal space $\eqs_{\ou}$. This set is a subset of the norming set of the transfinite mixed
Tsirelson space $T_{\ou}[(1/m_j,n_j)]$. The norm can also be described by the following
implicit formula, for $x\in c_{00}(\ou),$
\begin{align*}
\nrm{x}_*= \max\{& \nrm{x}_\infty,\sup_j\{\sup \frac1{m_{2j}}\sum_{i=1}^{n_{2j}}
\nrm{E_ix}_*,\;E_1<\cdots<E_{n_{2j}}\} \vee\\
&\vee\sup\conj{\frac{1}{m_{2j+1}}\sum_{i=1}^{n_{2j+1}}\phi_i(x)}{\{\phi_i\}_{i=1}^{n_{2j+1}}
\text{ is a $n_{2j+1}$-special sequence}}\}\}.
\end{align*}
This definition shares the same components with the corresponding
definition of the separable hereditarily indecomposable spaces.
The crucial difference concerns the definition of $n_{2j+1}$-{\it
special sequences.} For this we introduce a new coding
$\sigma_{\ro}$ based on a $\ro$ function which while it cannot be
one-to-one anymore it does provide a tree-like interference
between pairs of special sequences sufficient to impose a strong
conditional structure on $\eqs_{\ou}$.

The aim of the third section is to explain how the new $\ro-$coding is used in proving some of
the basic properties of the space $\eqs_{\ou}$. Thus, postponing the proofs of some
estimations for the next section, we show that block subsequences of $(e_\al)_{\al<\ou}$
generate hereditarily indecomposable subspaces. Section four contains the basic estimations
which are analogous to the $\omega-$case. We also show that $\eqs_{\ou}$ is reflexive. The
fifth section contains the study of the bounded linear operators. As we have mentioned above
many of the results are based on the finite representability of the James-like space $J_{T_0}$
in the transfinite block subsequences of $\eqs_{\ou}$. There are two ways to define $J_{T_0}$.
The first is the Bellenot-Haydon-Odell definition (\cite{bel-hay-od}) of the Jamesfication of
the mixed Tsirelson space $T_0=T[(1/m_{2j},n_{2j})_j]$ and the second is the Tsirelson-like
space $T[G,(1/m_{2j},n_{2j})_j]$ with $G=\{\chi_I:\;I \text{ interval of }\N\}$. The space
$J_{T_0}$ is quasi reflexive and for every set of ordinals $A$ the space $J_{T_0}(A)$ is
defined similarly to \cite{E}. The study of $J_{T_0}(A)$ and the finite representability of
$J_{T_0}$ in $\eqs_{\ou}$ are contained in the first two subsections of section five. The
remaining  subsections are devoted to the study of the spaces of operators. The central notion
of \emph{step diagonal operator} is defined as follows. Let $X$ be a subspace of $\eqs_{\ou}$
generated by a transfinite block sequence $(x_\alpha)_{\alpha<\gamma}.$ A bounded linear
diagonal operator $D:X\to X$ is a  {step diagonal operator} if $\lambda_\alpha=\lambda_\beta$
for all $\alpha,\beta$ with $\al\le\beta< \alpha+\omega\le \gamma.$ For example, if
$\gamma=\omega$ then $D=\lambda \mathrm{Id}_X$ and if $\gamma=\omega^2$ then
$D=\sum_{n}\lambda_n \mathrm{Id}_n$ where $\mathrm{Id}_n$ is the identity of
$X_n=\overline{\langle (x_\alpha)_{\alpha\in [\omega(n-1), \omega n)} \rangle}$. We prove the
following result.

\begin{teore*}
There exists a universal constant $C>0$ such that for every countable limit
ordinal $\gamma$ there is a set of ordinals
$A_\gamma$ such that: For every transfinite
block sequence $(x_\alpha)_{\alpha<\gamma}$ in $\eqs_{\ou},$ the
algebra $\mathcal{D}(\overline{\langle
(x_\alpha)_{\alpha<\gamma}\rangle})$ of the step diagonal
operators is $C$-isomorphic to $J^*_{T_0}(A_\kappa)$.
\end{teore*}

There are several consequences of this theorem: It follows readily that the structure of
$\mathcal{D}(X)$ for $X$ generated by a transfinite block sequence
$(x_\alpha)_{\alpha<\gamma}$ depends only on the ordinal $\gamma$. The dimension of
$\mathcal{D}(X)$ is equal to the cardinality of the set $ A_\gamma$, and for every $D\in
\mathcal{D}(X)$ and $\varepsilon>0$ there is an operator of the form $\sum_{i=1}^n
\lambda_i P_{I_i}$, with $\{I_i\}_{i=1}^n$ intervals of $\gamma$, which
$\varepsilon$-approaches $D$. Furthermore the following holds.

\begin{teore*}
 There exists a universal constant $C>0$ such that for every
 subspace $X$ of $\eqs_{\omega_1}$ generated by
 transfinite block sequence $(x_\alpha)_{\alpha<\gamma}$ of  $\eqs_{\ou}$ and every $T \in
\mathcal{L}(X, \eqs_{\ou})$ we have:
\begin{enumerate}
\item[(i)] $T=D_T+S_T$ where $D_T\in \mathcal{D}(X),$ $\|D_T\|
\leq C\|T\|$ and $S:X\to \eqs_{\ou}$ is strictly singular.
\item[(ii)] Every $D\in \mathcal{D}(X)$ is extendable to a
$\tilde{D}\in \mathcal{D}(\eqs_{\ou})$ with $\| \tilde{D} \| \leq
C \| D \|$. \item[(iii)] $\mathcal{L}(X, \eqs_{\ou})\cong
J^*_{T_0}(A_\kappa) \oplus \mathcal{S}(X, \eqs_{\ou})$.
\end{enumerate}
\end{teore*}
We also introduce the notion of  asymptotically equivalent subspaces of $\eqs_{\omega_1}$
which permits us to extend  part (iii) of the above theorem to arbitrary subspaces of
$\eqs_{\ou}$. Namely for every subspace $X$ of $\eqs_{\ou}$ there exists a set of ordinals
$A_X$ such that  $\mathcal{L}(X, \eqs_{\ou})\cong J^*_{T_0}(A_X) \oplus \mathcal{S}(X,
\eqs_{\ou})$. We are not able however to provide a sufficient description of $\mathcal{L}(X)$
for an arbitrary subspace $X$ of $\eqs_{\ou}$. What we have noticed is that in general $\mc
L(X,\eqs_{\ou})/\mc S(X,\eqs_{\ou}) \not\cong \mathcal{L}(X)/\mathcal{S}(X)$. For strictly
singular operators on $\eqs_{\ou}$ we give the following characterization.

\begin{teore*}
An operator $S:\eqs_{\ou}\to \eqs_{\ou}$ is strictly singular iff the sequence $
(\|S(e_\alpha)\|)_{\alpha<\omega_1}  \in c_0(\omega_1)$.
\end{teore*}

\begin{coro*}
Every $T\in \mathcal{L}(\eqs_{\ou})$ has the form $T=\lambda
\mathrm{Id}_{\eqs_{\ou}}+D+S$ where $D\in \mathcal{D}(X_\gamma)$ for some
$\gamma<\omega_1$, and $S$ is strictly singular. In particular $T=\lambda
\mathrm{Id}_{\eqs_{\ou}}+Q$ where $Q$ has separable range.
\end{coro*}
We mention that  nonseparable  spaces $X$ such that all $T\in\mathcal{L}(X)$ are of the form
$\lambda \mathrm{Id}_X +Q$ with  the range of  $Q$ separable have been constructed before in
\cite{sh}, \cite{sh-st} and \cite{wark}. However, those constructions  are quite different
from ours as they, in particular, offer no information about operators on separable subspaces
of the resulting space $X$.

Furthermore, we
show that for $I, J$ disjoint intervals of $\omega_1$ the spaces
$\eqs_I$ and $\eqs_J$ are totally incomparable and the space
$\eqs_{\omega_1}$ is arbitrarily distortable. Moreover, modulo strictly singular perturbations,
the space $\eqs_{\ou}$ admits a
unique resolution of the identity. Out of the rich sources of
examples of subspaces of $\eqs_{\omega_1}$ with interesting spaces
of operators we mention the following

\begin{teore*}
There exists a separable reflexive Banach space $\eqs$ admitting an infinite dimensional
Schauder decomposition $\eqs=\bigoplus_n \eqs_n$ such that, denoting by
$\mathcal{D}(\eqs)$ the class of bounded operators $D:\eqs\to \eqs$ with the property
$D|_{\eqs_n}=\lambda_n \mathrm{Id}_{\eqs_n}$ for all $n$, the following hold:
\begin{enumerate}
\item[(i)] $\mathcal{L}(\eqs)\cong \mathcal{D}(\eqs)\oplus
\mathcal{S}(\eqs)\cong J^*_{T_0}\oplus \mathcal{S}(\eqs)$.
\item[(ii)] For every subspace $X$ of $\eqs$ there exists
$A\subseteq \mathbb{N}$ which is either an initial finite interval
or is equal to $\mathbb{N}$ such that
$\mathcal{L}(X,\eqs)\cong J^*_{T_0}(A)\oplus
\mathcal{S}(X,\eqs). $
\item[(iii)] There is a subspace $X$
of $\eqs$ such that $\mathcal{L}(X)\cong \langle I_X\rangle \oplus \mathcal{S}(X)$ while
$\mathcal{L}(X,\eqs)\cong J^*_{T_0}\oplus \mathcal{S} (X,\eqs)$.
\end{enumerate}
\end{teore*}
\noindent For example, the space $\eqs =\eqs _{\omega^2}$ has all
these properties. It is worth pointing out that $D(\eqs)$ is a
natural class of operators which behaves similarly to the class of
operators of the form $\lambda I+K$ with $K$ compact diagonal. For
example if $x_n\in \eqs_n$ with $\|x_n\|=1$ and
$X=\overline{\langle (x_n)_n \rangle }$ then for every $D\in
\mathcal{D}(\eqs)$ we have that $D|_X=\lambda I+K$. The
isomorphism between $\mathcal{D}(\eqs)$ and $J^*_{T_0}$ endows
$J^*_{T_0}$ with an equivalent norm under which $J^*_{T_0}$ with
the pointwise multiplication becomes a commutative Banach algebra.
This should be compared with results from \cite{andr-gre}.

Sections six and seven concern two new properties that can be simultaneously imposed on a
$\ro$-function and the resulting properties of $\eqs_{\ou}$. First, we present a construction
of a universal $\ro$-function where universality roughly speaking means that for every
infinite interval $I$ of $\omega_1$ the finite $\ro-$closed subsets of $I$ realize all
isomorphism types of finite submodels of all possible $\ro$-functions. As we have mentioned
before, our initial motivation for introducing the universal $\ro$-function was to understand
the structure of $\mathcal{L}(\eqs_{\ou})$. However it turns out that using the $\ro -$coding
with a universal $\ro$ we obtain some new properties on $\eqs_{\ou}$ which have their own
interest, even for the space $\eqs_\omega=\overline{\langle e_n: n<\omega \rangle}$. Indeed
$\eqs_{\ou}$ admits a nearly subsymmetric transfinite basis and moreover $\eqs_\omega$, which
is an hereditarily indecomposable space, is an asymptotic version of itself \cite{Ma-Mi-To}.
The results concerning subsymmetric transfinite sequences and asymptotic versions are
presented in section seven. Section six also contains the construction of smooth
$\ro$-functions and the following result. If the coding $\sigma_\ro$ is based on a smooth
$\ro$-function then  every countable ordinal $\gamma<\omega_1$ can be re-ordered as
$(\al_n)_{n<\omega}$ such that $(e_{\al_n})_{n<\omega}$ defines a Schauder basis of the space
$X_\gamma$. Section eight contains a unified approach of the proof of two important results,
the basic inequality and the nontrivial direction of the finite representability of $J_{T_0}$.
Their proofs share some common features, and so we attempt to develop a general theory that
includes both results and that could be useful elsewhere. The last section is devoted to the
unconditional counterpart of the space $\eqs_{\ou}$ denoted by $\eqs_{\ou}^u$. The relation of
$\eqs_{\ou}^u$, which is a space with an unconditional basis $(e_\alpha)_{\alpha<\omega_1}$,
with the space $\eqs_{\ou}$ is same as that of Gowers-Maurey space with Gowers unconditional
space \cite{G2}. We study the structure of $\mathcal{L}(\eqs_{\ou}^u)$ and the structure the
subspaces of $\eqs_{\ou}^u$.

We extend our thanks to A. Tolias for his help during the
preparation of this paper. The results of this paper have been
announced in  \cite{arg-lop-tod}.


\section{Transfinite basic sequences}

The first section concerns the presentation of some preliminary
results related to transfinite (Schauder) bases. We recall one of
the equivalent formulations of their definition. For a detailed
presentation we refer the reader to \cite{sing}.

\defi
Let $X$ be a Banach space, and $\ga$ be  an ordinal number.

\noindent 1.  A total family $(x_\al)_{\al<\ga}$ of elements of
$X$ (i.e., a family such that $X=\overline{\langle
x_\al\rangle_{\al<\ga}}$) is said to be a \emph{transfinite basis}
if there exists a constant $C\ge 1$ such that for every interval
$I$ of $\ga$ the naturally defined map on the linear span of
$(x_\al)_{\al<\ga}$
\begin{align*}
\sum_{\al<\ga}\la_\al x_\al \mapsto\sum_{\al\in I}\la_\al x_\al
\end{align*}
extends to a bounded projection $P_I:X\to X_I=\overline{\langle x_\al\rangle_{\al\in I}}$ of
norm at most $C$.

\noindent 2. A transfinite basis $(x_\al)_{\al<\ga}$ of $X$ is said to be \emph{bimonotone} if for
each interval $I$ of $\ga$, the corresponding projection $P_I$ has norm 1.

\noindent 3. A transfinite basis $(x_\al)_{\al<\ga}$ of $X$ is said to be \emph{unconditional}
if there exists a constant $C\ge 1$ such that for all  subsets $A$ of $\ga$, the corresponding
$P_A$ has norm at most $C$.

\noindent 4. A transfinite basis $(x_\al)_{\al<\ga}$ of $X$ is said to be \emph{1-subsymmetric} if
for every $n\in \N$, every $\al_1<\al_2<\dots<\al_n<\ga$ and every $(\la_i)_{i=1}^n\in \R^n$,
$\nrm{\sum_{i=1}^n \la_i x_i}=\nrm{\sum_{i=1}^n \la_i x_{\al_i}}$.
\fdefi

\nota

\noindent 1. As in the case of the usual Schauder basis (i.e., $\ga=\om$) the above
definition is equivalent to the fact that each $x\in X$ admits a unique representation as
$\sum_{\al<\ga}\la_\al x_\al$, where the convergence of these series is recursively defined.

\noindent 2. The definition of  $\sum_{\al<\ga} \la_\al x_\al$ easily yields that for
each convergent series $\sum_{\al<\ga}\la_\al x_\al$  with $(x_\al)_{\al<\ga}$  a bounded family, the
sequence of coefficients $(\la_\al)_{\al<\ga}$ belongs to $c_0(\ga)$. Furthermore, for every $\vep>0$ there
exists a finite subset $F$ of $\ga$ such that $\nrm{\sum_{\al\notin F} \la_\al x_\al}<\vep$.

\noindent 3. For every transfinite basis $(x_\al)_{\al<\ga}$ the dual basis
$(x_\al^*)_{\al<\ga}$ is also well defined. Just like the usual Schauder bases,
$(x_\al^*)_{\al<\ga}$ is a $w^*$-total subset of $X^*$ and each $x^*$ in $X^*$ has a unique
representation of the form $\sum_{\al<\ga}x^*(x_\al) x_\al^*$ where the series is
$w^*$-convergent.

\noindent 4. If $(x_\al)_{\al<\ga}$ is a transfinite basis for the space $(X,\nrm{\cdot})$,
then there exists an equivalent norm $|\nrm{\cdot}|$ on $X$ such that $(x_\al)_{\al<\ga}$ is a bimonotone basis for
the space $(X,|\nrm{\cdot}|)$. This norm is defined by $|\nrm{x}|=\sup\conj{\nrm{P_I(x)}}{I\text{ interval of
}\ga}$.
\fnota

In the sequel, for every ordinal $\ga$ we shall denote by
$c_{00}(\ga)$ the vector space of all sequences $(\la_\al)_{\al\in
\ga}$ of real numbers such that the set
$\conj{\al<\ga}{\la_\al\neq 0}$ is finite. We also denote by
$(e_\al)_{\al<\ga}$ the natural Hamel basis of $c_{00}(\ga)$. It
is an easy observation that every space $X$ with a transfinite
basis $(x_\al)_{\al<\ga}$ is isometric to the completion of
$c_{00}(\ga)$ endowed with an appropriate norm. Moreover if $K$ is
a subset  of $c_{00}(\ga)$ with the properties (a)
$\{e_\al^*\}_{\al<\ga}\con K$ and (b)  for every $\phi\in K$,
$\nrm{\phi}_{\infty}\le 1$ and for every interval $I$ of $\ga$,
the restriction $\phi_I=\phi \cdot \chi_I$ of $\phi$ to $I$ is
also a member of $K$, then the norm defined on $c_{00}(\ga)$ by
\begin{align*}
\nrm{x}_K=\sup\conj{|\phi(x)|=\langle \phi,x\rangle}{\phi\in K}
\end{align*}
has $(e_\al)_{\al<\ga}$ as a transfinite bimonotone basis for the completion of
$(c_{00}(\ga),\nrm{\cdot}_K)$.

%
%

Fix $X$ with a transfinite basis $(x_\al)_{\al<\ga}$. The support $\supp x$ of $x\in X$ is the set
$\conj{\al<\ga}{x_\al^*(x)\neq 0}$. For a given interval $I\con \ga$, let $X_I=P_I X$, and for $\al<\ga$, let
$X_\al=X_{[0,\al)}$. For $x,y\in X$ finitely supported, we write $x<y$ to denote that $\max \supp x <\min
\supp y.$

A sequence $(y_{\al})_{\al<\xi}$ is called a transfinite block
subsequence of $(x_\al)_{\al<\ga}$ if and only if for all
$\al<\xi$, $y_\al$ is finitely supported and for all
$\al<\be<\xi$, $y_\al<y_\be$. Notice that a transfinite block
subsequence of a transfinite basis is always a transfinite basis
of its closed linear span.

%

Fix two Banach spaces $X$ and $Y$.   A bounded operator $T:X\to Y$ is an \emph{isomorphism}
iff $TX$ is closed and $T$ is 1-1.  $T$  is called \emph{strictly singular} if it is not an
isomorphism when restricted to any infinite dimensional closed subspace of $X$ (i.e., for all
infinite dimensional closed subspace $X'$ of  $X$, either $TX'$ is not closed or $T|X'$ is not
$1-1$). This is equivalent to say that for all  infinite dimensional closed subspace $Y$ of $
X$ and $\vep>0$, there is an infinite dimensional closed subspace $Y'$ of $Y$ such that
$\nrm{T|Y'}\le \vep$.

It is well known that most of the structure of the infinite dimensional closed subspaces of a
separable Banach space $X$ with a basis $(x_n)_n$ is described by its block sequences. Namely
that for every infinite dimensional closed subspace $Y$ of $X$ and every $\vep>0$ there exists
a normalized sequence in $Y$ and a block sequence $(w_n)_n$ of $(x_n)_n$ which are
$1+\vep$-equivalent. The method used for the proof of this result is called the \emph{gliding
hump} argument (\cite{li-tza}). This result is not extendable in the case of the transfinite
block sequences. For example, consider a biorthogonal basis $(x_\al)_{\al<\om\cdot 2}$ of a
Hilbert space and let $Y$ be the subspace generated by the sequence $(x_n+x_{\om+n})_n$.

We now describe  how block sequences are connected to  subspaces in the transfinite case.
\prop\label{gliding}
Let $(x_\al)_{\al<\ga}$ be a transfinite basis of $X$ and $Y$ an infinite  dimensional closed subspace $X$.
Then there exists a $\la\le \ga$ and a closed subspace $Z$ of $Y$ such that

\noindent 1. $P_{\la}:Z\to X_\la$ is an isomorphism.

\noindent 2. For every $\vep>0$ there exists a semi-normalized  block sequence $(w_n)_n$ in $X_\la$
and a normalized sequence $(z_n)_n$ in $Z$ such that $\sum_n \nrm{P_\la z_n-w_n}<\vep$.

\noindent 3. There exists a subspace $Z'$ of $Z$ isomorphic to a block subspace of $X$.

\noindent 4. If we additionally assume that $Y$ has a Schauder basis $(y_n)_n$, then the sequence
$(z_n)_n$ in 2. can be chosen to be a block sequence of $(y_n)_n$.
\fprop
\prue
We assume that $(x_\al)_{\al<\ga}$ is a bimonotone basis. Let
\begin{equation}
\be_0=\min\conj{\be}{P_\be:Y\to X_{\be} \text{ is not strictly singular}}.
\end{equation}
Let us show that $\la=\be_0$ is the required ordinal. Notice that $\be_0$ has to be
necessarily a limit ordinal. Since $P_{\be_0}$ is not strictly singular on $Y$,  there exists
a subspace $Z$ of $Y$ such that $P_{\be_0}:Z\to X_{\be_0}$ is an isomorphism. On the other
hand for every $\ga<\be_0$, $P_\ga:Y\to X_\ga$ is strictly singular hence for every $\vep>0$
and every subspace $Z'$ of $Z$ there exists $W\hookrightarrow Z'$\footnote{We will write $X
\hookrightarrow Y$ to denote that $X$ is an infinite dimensional closed subspace of the Banach
space $Y$.} such that $\nrm{P_\ga |W}<\vep$. Now we are ready to apply a modified gliding hump
argument to obtain $(z_n)_n$, $(w_n)_n$ as they are required in 2. Indeed for a given $\vep$
we choose $(\vep_n)_n$ such that $\vep_n>0$, $\sum \vep_n<\vep/4$. We choose a normalized
$z_1\in Z$ . Since $\be_0$ is a limit  ordinal, there must exist $\ga_1<\be_0$ such that
$\nrm{P_{[\ga_1,\be_0)}z_1}<\vep_1$. Hence setting $w_1=P_{\ga_1}z_1$ we have that
$\nrm{w_1-P_{\be_0}z_1}<\vep_1$. Since $P_{\ga_1} :Z\to X_{\ga_1}$ is strictly singular there
exists a normalized $z_2\in Z$ with $\nrm{P_{\ga_1}z_2}<\vep_2$. Choose $\ga_2>\ga_1$ such
that $\nrm{P_{[\ga_2,\be_0)}z_2}<\vep_2$ and set $w_2=P_{[\ga_1,\ga_2)}z_2$. Observe that
$\nrm{P_{\be_0}z_2-w_2}<2\vep_2$ and $w_1<w_2$. Continuing in this manner we obtain $(z_n)_n $
and $(w_n)_n$ such that for all $n$, $\nrm{P_{\be_0}z_n-w_n}\le 2\vep_n$, hence
\begin{equation}\label{eqthree}
\sum_n \nrm{P_{\be_0}z_n-w_n}\le \vep/2.
\end{equation}
Since we assume that the transfinite basis $(x_\al)_{\al<\ga}$ is bimonotone, (\ref{eqthree})
implies that $(P_{\be_0}z_n)_n$ and $(w_n)_n$ are equivalent.  Property 3. follows from 2.,
while 4. results from a careful choice of $(z_n)_n$ in 2. \fprue

As we have mentioned in the introduction the manner that block subspaces saturate the
subspaces of $X$ is weaker than the corresponding result for spaces $X$ with a basis
$(x_n)_n$. In the next proposition we provide a sufficient condition which ensures the
complete extension of the result from  Schauder bases to transfinite Schauder bases fulfilling
the additional condition.

\prop\label{stronggliding} Let $(x_{\al})_{\al<\gamma}$ be a transfinite basis of $X$.
Assume that for all  disjoint intervals $I,J$ of $\gamma$ the spaces $X_I$ and $X_J$ are
totally incomparable. Then for every  closed infinite dimensional subspace $Y$ of $X$ and
every $\varepsilon>0$ there exist  normalized  sequences $(y_n)_n$, $(z_n)_n$ such that
$(y_n)_n\subseteq Y$, $(z_n)_n$ is a block sequence of $(z_\al)_{\al<\gamma}$ and
$\sum_n\|y_n-z_n\|<\varepsilon$.
 \fprop
\prue From Proposition \ref{gliding} there exists a subspace $Z$
of $Y$ and $\lambda\le\gamma$ such that $P_\lambda:Z\to
X_{\lambda}$ is an isomorphism. Assume that $\lambda<\gamma$ and
set $I=[1,\lambda)$ and $J=[\lambda,\gamma)$. Then $P_{J}:Z\to
X_{J}$ is a strictly singular operator. Hence we may find $(w_n)$,
$(z_n)$ as in Proposition \ref{gliding} (2) such that
$\sum_n\|P_{J}(z_n)\|<\varepsilon$ which yields that
$\sum_n\|z_n-w_n\|<2\varepsilon$. \fprue

\defi
A transfinite basis $(x_\al)_{\al<\ga}$ is called \emph{shrinking
} iff for all $(\al_n)_n \uparrow$, $(x_{\al_n})_n$ is shrinking
in the usual sense (i.e., $(x_{\al_n}^*)$ generates in norm the
dual of the closed span of $(x_{\al_n})_n$).

 It is called \emph{boundedly
complete} iff for all $(\al_n)_n \uparrow$, $(x_{\al_n})_n$ is
boundedly complete in the usual sense (i.e., for all sequence of
scalars $(\la_n)_n$, if there is some $C>0$ such that for all $n$,
$\nrm{\sum_{i=1}^n \la_i x_{\al_i}}\le C$, then $\sum_i \la_i
x_{\al_i}$ converges in norm).

\fdefi

The above definitions are simpler and easier checked than the corresponding ones cited in
\cite{sing}. The following result  is the extension of the well-known James' characterization
(\cite{li-tza}) of reflexivity in the general setting of a Banach space with a transfinite
basis.

\prop\label{jamesgencrit} Let $(x_\al)_{\al<\ga}$ be a transfinite basis of $X$. Then $X$
is reflexive iff $(x_\al)_{\al<\ga}$ is shrinking and boundedly complete. \fprop
\prue
The direct implication is consequence of the  James' characterization (\cite{jam}). The
opposite requires the following two Claims: \clam If $(x_\al)_{\al<\ga}$ is shrinking,
then the biorthogonal basis $(x_\al^*)_{\al<\ga}$ generates in the norm topology the dual
space $X^*$. \fclam \prucl Assume the contrary. Then there exists $x^*\in X^*$ not in the
closed  linear span $Y$ of $(x_\al^*)_{\al<\ga}$. Set $\be_0=\min \conj{\be\le
\ga}{P_\be^* x^*\notin Y}$. Then $P_{\be_0}^*x^*\notin Y$ but for all $\ga<\be_0$,
$P_\ga^*x^*\in Y$. Therefore there exists an increasing sequence of successive disjoint
intervals $I_1<I_2<\dots <I_n<\dots <\be_0$ and $\vep>0$ such that for each $n\in \N$,
$P_{I_n}^*x^*\in Y$ and $\nrm{P_{I_n}^* x^*}\ge \vep$. Observe that if $x^*\in X^*$,
$x^*= w^*-\sum_{\al<\ga}\mu_\al x_\al^*$, where for each $\al<\ga$, $\mu_\al=x^*(x_\al)$.
Moreover if $I$ is an interval of $\ga$ such that $P_I^* x^*\in Y$ and $\vep'>0$, then
there is a finite subset $F_{\vep'}$ of $I$ such that $\nrm{y_{\vep'}^*-x^*}<\vep$ where
\begin{align*}
y_{\vep'}^*=w^*-\sum_{\al\in \ga \setminus I}\mu_\al x_\al^*+\sum_{\al\in F_{\vep'}} \mu_\al
x_\al^*.
\end{align*}
Using this observation we inductively select finite sets $F_1\con I_1$,\dots, $F_n\con I_n$ such that setting
\begin{equation}\label{noone}
y_n^*=\sum_{i=1}^n \sum_{\al\in F_i} \mu_\al x_\al^*+P_{\be_0\setminus \bigcup_{i=1}^n I_n} x^*,
\end{equation}
we have that
\begin{equation}\label{notwo}
\nrm{P_{\be_0}^* x^*-y_n^*}<\vep_n<\frac\vep 4.
\end{equation}
Set $y^*=w^*-\lim_n y_n^*$ and (\ref{noone}) and (\ref{notwo}) yield that $\supp y_n^*\con
\bigcup_n F_n$ and also $\nrm{P_{F_n}^*y^*}>\vep/2$. Since each $F_n$ is a finite  set we can
enumerate $\bigcup_n F_n$ as $(\al_n)_n\uparrow$ and clearly $y^*$ yields that the sequence
$(x_{\al_n})_n$ is not a shrinking Schauder basis, yielding a contradiction. \fprucl
 \clam If $(x_\al)_{\al<\ga}$ is boundedly complete, then for every
$x^{**}\in X^{**}$, the series $\sum_{\al<\ga}x^{**}(x_\al^*)x_\al$ converges in norm.
\fclam \prucl Suppose the contrary and fix $x^{**}\in X^{**}$ but not in $X$.  The proof
is similar to the previous one.
 For each $\al<\ga$, let $\la_\al=x^{**} x_\al^*$ and let
 \begin{align*}
 \be_0=\min\conj{\be<\ga}{P_\be^{**}x^{**}\notin X}.
 \end{align*}
 Using a similar argument we can choose an increasing sequence $(F_n)_n$ of finite subsets of $\ga$ such that
 $w^*-\sum_{\al\in \bigcup_n F_n} \la_\al x^*_\al$ exists and for every $n$, $\nrm{\sum_{\al\in F_n} \la_\al x_\al
 }>\vep>0$. This yields that the sequence $(x_\al)_{\al\in \bigcup_n F_n}$ is not boundedly complete, a
 contradiction.
\fprucl

\fprue
\section{Definition of the space $\eqs_{\ou}$}
This section is devoted to the definition of the norm of the space $\eqs_{\ou}$. This
norm will be induced by a set of functionals, denoted by $K_{\ou}$, on the space
$c_{00}(\ou)$. Then $\eqs_{\ou} $ will be the completion of it. We start with a short
presentation of the unconditional frame, which is a   mixed Tsirelson space with
1-subsymmetric transfinite basis of a given length $\ga$. The aforementioned set
$K_{\ou}$ will be selected as a subset of $B_{Y^*}$ where $Y$ is the corresponding mixed
Tsirelson space.
\subsection{The space $T_\ga[(1/m_j,{n_j})_j]$}
Throughout the paper we fix two infinite sequences $(m_j)_j$, $(n_j)_j$ defined recursively as
follows:
\begin{enumerate}
\item $  m_1=2, \text{ and } m_{j+1}=m_j^4$
\item $ n_1=4, \text{ and }  n_{j+1}=(4n_j)^{s_j}\text{
 where } s_j=\log_2 m_{j+1}^3.$
\end{enumerate}
Let $\ga$ be an infinite ordinal. Consider the norm
$\nrm{\cdot}_*$ on $c_{00}(\ga)$ described by the implicit formula
\begin{align*}
\nrm{x}_*=\max\{\nrm{x}_\infty,\sup_j\sup \frac1{m_j}\sum_{i=1}^{n_j} \nrm{E_i}_*\},
\end{align*}
where for $E\con \ga$, $x\in c_{00}(\ga)$  $Ex$ denotes the restriction  of $x$ to the set $E$
(i.e., $Ex=P_E x=\langle \chi_E, x \rangle$) and the inside supremum is taken over all
sequences $E_1<\dots <E_{n_j}$ of subsets of $\ga$.

The existence of a norm satisfying the above formula is provided, as the case of Tsirelson space, by an
inductive argument (e.g. \cite{li-tza}). It is also easy to see that the usual basis $(e_\al)_{\al<\ga}$ of
$c_{00}(\ga)$ defines a 1-subsymmetric and 1-unconditional basis for the space
\begin{align*}
T_\ga[(m_j^{-1},n_j)_j]=\overline{(c_{00}(\ga),\nrm{\cdot}_*)}.
\end{align*}
The first variation of the original Tsirelson construction is due
to Th. Schlumprecht \cite{schlu} who introduced the space
$S=T_\om[(1/\log_2(j+1),j)]$ providing the first known example of
an arbitrarily distortable Banach space. The space $S$ is one of
the key ingredients in the Gowers-Maurey construction \cite{go-ma}
of a Banach space with no unconditional basic sequence. The
general definition of a mixed Tsirelson space
$T_\ga[(m_j^{-1},n_j)_j]$ for $\ga=\om$ was introduced in
\cite{arg-del} using the slightly different notation $T[(\mc
A_{n_j},1/m_j)_j]$ which stresses the use of the family $\mc
A_{n_j}$ of all subsets of the index-set (in their case $\omega$)
of cardinality at most $n_j$ and indicates the possibility to use
some other compact family instead of $\mc A_{n_j}$ (see e.g.
\cite{arg-del-man} and \cite{arg-tol}). Since in this paper we are
not going to vary the definition in this direction we suppress the
$\mc A$  as this will give us some notational advantages at some
latter points of the paper.
\nota \noindent 1. It follows readily from the definition of the norm that for $A\con \ga$ with
order type of $A$ equal to the ordinal $\la$ the space $X_A=\overline{\langle e_\al
\rangle_{\al\in A}}$ is isometric to $T_\la[( m_j^{-1},{n_j})_j]$. Therefore granting that
$T_\om[( m_j^{-1}, {n_j})_j]$ is reflexive (e.g. \cite{schlu}, \cite{arg-mano})  Proposition
\ref{jamesgencrit} yields that for each $\ga$ the space $T_\ga[( m_j^{-1},{n_j})_j]$ is also
reflexive.

\noindent 2. A possible variation of the norm of $T_\ga[( m_{j}^{-1}, {n_j})_j]$ is to
allow sequences $(E_{1},\dots,E_{n_j})$ consisting of disjoint sets (i.e., not
necessarily successive). Such spaces are called modified mixed Tsirelson spaces and they
are denoted by $T^{\mc M}[(m_{j}^{-1}, {n_j})_j]$. Schlumprecht has shown that $S^\mc M$
contains $\ell^1$ while such spaces have been studied in \cite{arg-del-kut-man},
\cite{man}, \cite{arg-del-man}. The situation for the spaces $T^\mc M_\om[(
m_j^{-1},{n_j})_j]$ remains unclear. Namely, we do not know if there exists a sequence $(
q_j^{-1},{l_j})_j$ such that the space  $T^\mc M_\om[( q_j^{-1},{l_j})_j]$ is reflexive
and not containing any $\ell_p$, $1< p<\infty$.\fnota

There exists an alternative definition of the norm of $T_\ga[(m_{j}^{-1}, {n_j})_j]$ which is
close to the definition of the norm of $\eqs_{\ou}$. This goes as follows.

Let $L_\ga \con c_{00}(\ga)$ be the minimal subset $L$ of $c_{00}(\ga)$ satisfying the
following four properties:
\begin{enumerate}
\item For every $\phi\in L$ and
every $E\con \ga$, $E\phi\in L $.
\item For every $\al<\ga$, $\pm e_{\al}^*\in L$.
\item For every $j\in \N$ and every $\phi_1<\dots <\phi_{n_j}$ in $L$,
$(1/m_j)\sum_{i=1}^{n_j} \phi_i $ also belongs to $L$.
\item $L$ is closed under rational convex combinations.
\end{enumerate}
The third property is also described by saying that $L$ is closed
in all $( m_j^{-1},{n_j})$-operations. It is not difficult to see
that the norm induced on $c_{00}(\ga)$ by the set $L_{\ga}$ (i.e.,
for $x\in c_{00}(\ga)$, $\nrm{x}=\sup_{\phi\in L_\ga}\{\phi
x=\langle \phi,x \rangle\}$) is exactly the norm $\nrm{\cdot}_*$.
\nota Let $L'_\ga$ be the minimal subset of $c_{00}(\ga)$ satisfying  1., 2. and 3. It is
not difficult to prove that $L_\ga={\conv}_{\Q}(L_\ga')$. This means that $L'_\ga$ norms the
space $T_\ga[( m_j^{-1},{n_j})_j]$.
 \fnota

\nota
\noindent 1. It follows from the minimality of $L_\ga$ that each $\phi\in L_\ga$ is
either equal to $\pm e_\al$ for some $\al<\ga$ or is of the form $\phi=(1/m_j)\sum_{i=1}^d \phi_i$, $d\le
n_j$ and $\phi_1<\dots <\phi_d$ all in $L_\ga$. Furthermore the set
\begin{align*}
L_{\ga,j}=\conj{\phi\in L_\ga}{\phi=\frac1{m_j}\sum_{i=1}^d \phi_i}
\end{align*}
defines an equivalent norm, denoted by  $\nrm{\cdot}_{*,j}$, on the space  $T_\ga[(m_j^{-1}, {n_j})_j]$.
The important property of the mixed Tsirelson spaces results from a fine balance of the sequences
of norms $(\nrm{\cdot}_{*,j})_j$. Namely for every block sequence $(x_n)_n$ and for every $j$ there
exists a normalized vector $y_j$ in the linear span of $(x_n)_n$ such that $\nrm{y_j}_{*,j}>1/4$
and for every $j'\neq j$ $\nrm{y_j}_{*,j'}<6/m_{j'}$ if $j'<j$ and $\nrm{y_j}_{*,j'}<4/m_{j}^2$
otherwise.
 \fnota

\subsection{The norming set $K_{\omega_1}$}\label{definitionnorming} The maximal space in our class
$\eqs_{\omega_1}$ will be defined
as the completion of $(c_{00},\|\cdot\|_{\oo})$ under the norm $\|\cdot\|_{\oo}$ induced
by a set of functionals $K_{\ou} \subseteq c_{00}(\ou)$.

The set $K_{\omega_1}$ is the minimal subset of $c_{00}(\omega_1)$ satisfying that:
\begin{enumerate}
\item[\emph{(1)}] It contains $(e_{\gamma}^*)_{\gamma<\omega_1}$,
is symmetric (i.e., $\phi\in K$ implies $-\phi\in K$) and is
closed under the restriction on intervals of $\omega_1$.
\item[\emph{(2)}] For every
$\{\phi_i:\;i=1,\ldots,n_{2j}\}\subseteq \ko$ with $\supp\phi_1<\dots<\supp\phi_{n_{2j}}$, the
functional $\phi=({1}/{m_{2j}})\sum_{i=1}^{n_{2j}}\phi_i\in\ko$. We say that $\phi$ is a
result of a $(m_{2j}^{-1},n_{2j})$-operation.
\item[\emph{(3)}] For every special sequence
$(\phi_1,\ldots,\phi_{n_{2j+1}})$ (for a definition, see subsection \ref{sseq}), the
functional $\phi=({1}/{m_{2j+1}})\sum_{i=1}^{n_{2j+1}}\phi_i$ is in $K_{\ou}$. We call  $\phi$
a \emph{special functional} and say that $\phi$ is a result of a
$(m_{2j+1}^{-1},{n_{2j+1}})$-operation. \item[\emph{(4)}] It is rationally convex.
 \end{enumerate}
The norm on $c_{00}(\omega_1)$ is defined as $\nrm{x}=\sup\{\phi(x)=\sum_\al\phi(\alpha)\cdot
x (\alpha):\;\phi\in\ko\}$ and $\eqs_{\omega_1}$ is the completion of $(\co,\nrm{\cdot})$. Each  of
the above four properties provides certain features in the space $\eqs_{\omega_1}$. The first
makes the family $(e_{\al})_{\al<\omega_1}$ a transfinite bimonotone basis of
$\eqs_{\omega_1}$. The second saturates $\eqs_{\omega_1}$ with local unconditional structure.
This property will be responsible for the existence of semi-normalized averages in every block
sequence of $\eqs_{\omega_1}$. The third property saturates $\eqs_{\omega_1}$ with conditional
structure and will make it impossible for $\eqs_{\omega_1}$ to contain any unconditional basic
sequence.
Finally, the fourth property is a tool for proving properties of the space of operators from an
arbitrary subspace $X$ of $\eqs_{\omega_1}$ into $\eqs_{\omega_1}$. The above definition, with the
exception of the missing definition of special sequences, is based on the corresponding
definitions from \cite{arg-mano} and \cite{arg-tol} which in turn are variants of the construction
from \cite{go-ma}. By the minimality of $K_{\omega_1}$ each $\phi\in\ko$ has one of the following
forms:
\begin{enumerate}
\renewcommand{\labelenumi}{(\roman{enumi})}
\item[(i)] $\phi$ is of  \emph{type 0} if $\phi =\pm e_\al^*$.
\item[(ii)] $\phi$ is of  \emph{type I} if $\phi =\pm Ef$ for  $f$
a result of one $(m_{j}^{-1},{n_{j}})$-operation and $E$ an interval. In this case we say
that the weight $w(\phi)$ of $\phi$ is $m_j$.
 \item[(iii)] $\phi$ is of \emph{type II} if $\phi$ is a rational
convex combination of type 0 and type I functionals.
\end{enumerate}
An alternative description of the norm is the following: For a given $x\in \eqs_{\ou}$,
\begin{align*}
\nrm{x}_*=\max\{\nrm{x}_\infty,\sup_j\sup \frac1{m_{2j}}\sum_{i=1}^{n_{2j}}
\nrm{E_ix}_*,\;E_1<\cdots<E_{n_{2j}}\} \vee\\
\sup\conj{\frac{1}{m_{2j+1}}\sum_{i=1}^{n_{2j+1}}\phi_i(x)}{\{\phi_i\}_{i=1}^{n_{2j+1}}
\text{ is a $n_{2j+1}$-special sequence}}.
\end{align*}

\nota\label{basbound} From the definition of the norming set
$K_{\ou}$ it follows easily that $(e_{\al})_{\al<\ou}$ is a
bimonotone basis of $\eqs_{\ou}$. Also, it is not difficult to see
using  \emph{(2)} from the definition of $K_{\ou}$ that the basis
$(e_\al)_{\al<\ou}$ is boundedly complete. Indeed, for $x\in
c_{00}(\ou)$ and $E_1<\dots <E_{n_{2j}}$ intervals of $\ou$
\emph{(2)} of the norming set yields that $\nrm{x}\ge
(1/{m_{2j}})\sum_{i=1}^{n_{2j}}\nrm{E_i x}$. Also, from the choice
of the sequence $(m_j)_j$, $(n_j)_j$, it follows that
$n_{2j}/m_{2j}$ increases to infinity. From these observations it
follows that the basis $(e_\al)_{\al<\ou}$ is boundedly complete.
To prove that the space $\eqs_{\ou}$ is reflexive we need to show
that the basis is shrinking.
 \fnota

\defi
For $\phi\in\ko$, we say that  $m_j\in \N$ is a  \emph{weight} of
$\phi$, or $w(\phi)=m_j$, if $\phi$ can  be obtained as a result
of the $(m_{j}^{-1},{n_{j}})$-operation. Notice that $\phi\in\ko$
may have many weights. \fdefi The definition of the special
sequences will, as in the case \cite {go-ma}, depend crucially on
certain coding $\sigma_\ro$. The essential difference is that now
$\sigma_{\ro}$ is not an injection, a crucial property on which
the proofs in \cite {go-ma} rely. Our proofs on the other hand
will rely on a ``tree-like property" of our coding which we now
describe. First we notice that each $2j+1$-special sequence
$\Phi=(\phi_1,\phi_2,\ldots,\phi_{n_{2j+1}})$ is of the form
$\supp\phi_1<\dots<\supp\phi_{n_{2j+1}}$ with each $\phi_i$ of
type I. The  \emph{tree-like property} is the following: For any
pair of $2j+1$-special sequences
$\Phi=(\phi_1,\phi_2,\ldots,\phi_{n_{2j+1}})$,
$\Psi=(\psi_1,\psi_2,\ldots,\psi_{n_{2j+1}})$ there exist $1\le
\ka_{\Phi,\Psi}\le\lambda_{\Phi,\Psi}\le n_{2j+1}$ such that

  \begin{enumerate}
\renewcommand{\labelenumi}{(\roman{enumi})}
\item[(i)] If $1\le k< \ka_{\Phi,\Psi}$ then $\phi_k=\psi_k$ and
if $\ka_{\Phi,\Psi}<k<\lambda_{\Phi,\Psi}$, then
$w(\phi_k)=w(\psi_k).$ \item[(ii)]
$(\cup_{\ka_{\Phi,\Psi}<k<\lambda_{\Phi,\Psi}}\supp \phi_k)\cap
(\cup_{\ka_{\Phi,\Psi}<k<\lambda_{\Phi,\Psi}}\supp\psi_k)
=\emptyset.$ \item[(iii)]
$\{w(\phi_k):\;\lambda_{\Phi,\Psi}<k<n_{2j+1}\}\cap
\{w(\psi_k):\;\lambda_{\Phi,\Psi}<k<n_{2j+1}\}=\emptyset$.
\end{enumerate}

Comparing the above tree-like property with the corresponding property from \cite {go-ma}, we
notice that the new ingredient is the number $\ka_{\Phi,\Psi}$. Its occurrence is a byproduct
of the fact that the coding $\sigma_{\ro}$ is not one-to-one. The property (ii) will however
give  sufficient control of our special functionals. The coding $\sig_\ro$ is based on the
following mapping introduced in \cite{tod1} (see also \cite{tod2}).

\subsection{$\ro$-functions}
A function $\varrho:[\ou]^2\to \omega$ such that:
\begin{enumerate}
\item $\varrho(\al,\ga)\le \max\{\varrho(\al,\be),\varrho(\be,\ga)\}$ for all
$\al<\be<\ga<\ou$.
\item $\varrho(\al,\be)\le \max\{\varrho(\al,\ga),\varrho(\be,\ga)\}$ for all
$\al<\be<\ga<\ou$.
\item $\conj{\al<\be}{\varrho(\al,\be)\le n}$ is finite for all $\be<\ou$ and $n\in \N$.
\end{enumerate}
is called a $\ro$-function. The reader is referred to \cite{tod1} and \cite{tod2} for
full discussion of this notion and constructions of various $\ro$-functions. In Section
\ref{difros} we shall give yet another construction of a $\ro$-function with certain
universality property.

Let $\varrho:[\ou]^2\to \om$ be a $\varrho$-function fixed from now on, and all definitions
and facts that follow should be relative to this choice of $\ro$.

\defi
Given a finite set $F\con \ou$, let $p_F=p_\varrho(F)=\max_{\al,\be\in F}\varrho(\al,\be)$.
For a finite  set $F\con \ou$ and $p\in \N$,  let
$$ \cl{F}{p}=\conj{\al\le \max F}{\text{there
is }\be \in F \text{ s.t. }\al \le \be \text{ and  } \varrho(\al,\be)\le p}.$$ Notice that by
condition 3., $\cl{F}{p}$ is a finite set of countable ordinals. We say that $F$ is
\emph{$p$-closed} iff $\cl{F}{p}=F$, and that $F$ is $\ro$-closed iff it is $p_F$-closed.

\fdefi

\nota\label{roclosrem} \noindent 1. Note that $\cl{\cdot}{p}$ is a
monotone and idempotent operator and so, in particular, every
$\cl{F}{p}$ is a $p$-closed set: It  is clear that if $F\con G$,
then $\cl{F}{p}\con \cl{G}{p}$. Let us show now that
$\cl{\cl{F}{p}}{p}=\cl{F}{p}$. Let $\al \in \cl{\cl{F}{p}}{p} $.
This implies that $\varrho(\al,\al_0)\le p$, for some $\al_0\in
\cl{F}{p}$, $\al\le \al_0$. Choose $\al_1\ge \al_0$, $\al_1\in F$
such that $\ro(\al_0,\al_1)\le p$. Then, $\ro(\al,\al_1)\le
\max\{\ro(\al,\al_0),\ro(\al_0,\al_1)\}\le p$.

\noindent 2. Suppose that $F\con \ou$ is finite and suppose that $p\ge p_F$. Then $p_{\cl
F p}\le p$. Indeed, let  $\al<\be$ such that both   belong to $\cl F p$. Let $\al'\ge
\al$, $\be'\ge \be$ such that $\al,\be'\in F$ and $\ro(\al,\al'),\ro(\be,\be')\le p$.
Then we distinguish the following cases:

\noindent(a) If $\al\le \al'\le \be \le \be'$, then $\ro(\al,\be)\le
\max\{\ro(\al,\al'),\ro(\al',\be)\}\le \max\{\ro(\al,\al'),\ro(\al',\be'),\ro(\be,
\be')\}\linebreak \le p$.

\noindent(b) If $\al\le  \be \le \al' \le \be'$, then $\ro(\al,\be)\le
\max\{\ro(\al,\al'),\ro(\be,\al')\} \le\max\{\ro(\al,\al'),\ro(\be,\be'), \ro(\al',
\be')\}\linebreak \le p$.

\noindent(c) If $\al   \le \be\le \be'\le \al'$, use a similar proof to case (a). \fnota

\prop \label{easyro} Let $F,G\con \ou$ be two finite sets and $p\ge p_F,p_G$. Then:
\begin{enumerate}
\item For every
ordinal $\al\le \ou$, $\cl{F\cap \al}{p}=\cl{F}{p}\cap \al$ and $\cl{F\cap \al}{p}$ is an
initial part of $\cl{F}{p}$. Therefore, if $F$ is $p$-closed, so is $F\cap \al $.
\item For every $\al\in F\cap G$, we have that  $\cl{F\cap (\al+1)}{p}=\cl{G\cap (\al+1)}{p}$. Hence,
if $F$ and $G$ are in addition $p$-closed, then $F\cap (\al+1)=G\cap (\al+1)$.
\item  $\cl {F\cap G} p=\cl F p \cap \cl G p$. Therefore, if $F$ and $G$ are $p$-closed
then $F\cap G$ is also $p$-closed and it is an initial part of both $F$ and $G$.
\end{enumerate}
\fprop \prue 1: Since $F\cap \al\con F,\al$, it follows that $\cl{F\cap \al}{p}\con
\cl{F}{p}\cap \al$. Now let $\be\in \cl{F}{p}\cap \al$. Then there is some $\ga\in F$,
$\ga\ge be$ such that $\ro(\be,\ga)\le p$. If $\ga<\al$, then we are done. If not, let
$\de=\max F\cap \al\in F$ and since $\be\le \de<\ga$ we have that
\begin{equation}\label{kjuie1}
\ro(\be,\de)\le \max\{\ro(\be,\ga),\ro(\de,\ga)\}\le \max\{p,p_F\}=p,
\end{equation}
the last equality using our assumption that $p\ge p_F$. (\ref{kjuie1}) shows that $\be\in
\cl{F\cap \al}{p}$. Suppose now  that $F$ is $p$-closed. Then we have  just shown that
$\cl{F\cap \al }{p}=\cl{F}{p}\cap \al=F\cap \al$, and we are done.

2: Fix $\al\in F\cap G$.  Let $\be\in \cl{F\cap (\al+1)}{p}=\cl{F}{p}\cap (\al+1)$. Let
$\ga\in F\cap (\al+1)$, $\ga\ge \be$ be such that $\ro(\be,\ga)\le p$.  Then
$\ro(\be,\al)\le \max\{\ro(\be,\ga),\ro(\ga,\al)\}\le \max\{p,p_F\}=p$. Since $G$ is
$p$-closed, and $\al\in G$, we can conclude that $\be\in \cl{G\cap (\al+1)}{p}$. This
shows that $\cl{F\cap (\al+1)}{p}\con \cl{G\cap (\al+1)}{p}$. The other inclusion follows
by symmetry. The last part of 2. follows easily.

3: Let $\al=\max F\cap G$. Then by 2., $\cl{F\cap G}p=\cl {F\cap G \cap (\al+1)}p= \cl {F
\cap (\al+1)}p=\cl {F}p \cap (\al+1)$ and   $\cl{F\cap G}p=\cl G p\cap (\al+1)$.
Combining the above equalities we get $\cl{F\cap G} p=\cl F p\cap \cl G p\cap (\al+1)=
\cl F p\cap \cl G p$, the last equality because $\cl F p\cap \cl G p\con F\cap G\con \max
(F\cap G)+1=\al+1$. \fprue

\subsection{The $\sig_\ro$-coding and the special sequences}\label{sseq}
 We denote by $\Q_s(\ou)$ the set of finite sequences
$(\phi_1,w_1,p_1,\phi_2,w_2,p_2,\ldots,\phi_d,w_d,p_d)$ such that
\begin{enumerate}
\item for all $i\le d $, $\phi_i\in c_{00}(\ou)$ and $\phi_1<\phi_2<\dots<\phi_d$,
\item $(w_i)_{i=1}^d, \,(p_i)_{i=1}^d \in \N^d$ are strictly increasing, and
\item $p_i \ge  p_{(\cup_{k=1}^i\supp\phi_k)}$ for every $i\le d$.
\end{enumerate}
Let  $\Q_s$  be the set of finite sequences $(\phi_1,w_1,p_1,\phi_2,w_2,p_2,\ldots,\phi_d,w_d,p_d)$
satisfying 1., and 2. above and in addition for every $i\le d$, $\phi_i\in c_{00}(\N)$. Notice that
$\Q_s$ is a countable set.
 Fix a
one-to-one function $\sigma:\Q_s\to\{2j:\;j \text{ odd}\}$ such that
\[  \sigma(\phi_1,w_1,p_1,\phi_2,w_2,p_2,\ldots,\phi_d,w_d,p_d)>
\max\{p_d^2,\frac{1}{\varepsilon^2},\max\supp \phi_d\}, \] where
$\varepsilon=\min\{|\phi_k(e_\alpha)|:\; \alpha\in\supp \phi_k,\;k=1,\ldots ,d\}$.  Given a
finite subset $F$ of $\omega_1$, we denote by $\pi_F:\{1,2,\ldots,\# F\}\to F$ the natural
order preserving map. Given $\Phi=(\phi_1,w_1,p_1,\phi_2,w_2,p_2,\ldots, \linebreak
\phi_d,w_d,p_d)\in\Q_s(\ou)$ we set $G_{\Phi}=\cl{\cup_{i=1}^d\supp \phi_i}{p_d}$ and then we
consider the family
\[
\pi_{G_{\Phi}}(\Phi)=(\pi_G(\phi_1),w_1,p_1,\pi_G(\phi_2),w_2,p_2,\ldots,\pi_G(\phi_d),w_d,p_d)
\in\Q_s,\] where
\begin{align*}
\pi_G(\phi_k)(n)=\left\{\begin{array}{ll} \phi_k(\pi_{G_\Phi}(n)) & \text{if }n\in
G_\Phi \\
0 & \text{otherwise.}
\end{array}
\right.
\end{align*}
Finally, $\sig_\ro:\Q_s(\ou)\to \conj{2j}{j \text{ odd}} $ is defined as
$\sigma_\ro(\Phi)=\sigma(\pi_G(\Phi))$.

A sequence $\Phi=(\phi_1,\dots,\phi_{n_{2j+1}})$ of functionals of $K_{\omega_1}$ is said to be a
$2j+1$-\emph {special sequence} if:

 \noindent (1) $\supp \phi_1<\supp \phi_2<\dots<\supp
\phi_{n_{2j+1}}$, each $\phi_k$ is of type I, $w(\phi_k)=m_{2j_k}$ and $w(\phi_1)=m_{2j_1}$ with
$j_1$ even and satisfying $m_{2j_1}> n_{2j+1}^2$.

\noindent (2) There exists a strictly increasing sequence $(p_1^{\Phi},\dots,
p^{\Phi}_{n_{2j+1}-1})$ of natural numbers such that for all $1\le i\le n_{2j+1}-1$ we
have that $w(\phi_{i+1})=m_{\sigma_\ro(\Phi_i)}$ where  $\Phi_i=(\phi_1,
w(\phi_1),p_1^{\Phi}, \phi_2,w(\phi_2),p_2^{\Phi}, \ldots ,\linebreak \phi_i ,
w(\phi_i),p_i^{\Phi})$.

As we have mentioned before, the weight  of a type I element of $K_{\ou}$ is not uniquely
determined. However in the case of the elements $\phi_i$ of a $2j+1$-special sequence
$\Phi$, $w(\phi_i)$ will denote the unique weight involved in the definition of the
special sequence $\Phi$.

%

\lema[\emph{Tree-like interference of a pair of special sequences}]\label{0011}
Let $\Phi=(\phi_1,\dots,\phi_{n_{2j+1}})$ and $\Psi=(\psi_1,\dots,\psi_{n_{2j+1}})$ be two $2j+1$-special
sequences. Then there are two numbers $0 \le \ka_{\Phi,\Psi}\le  \la_{\Phi,\Psi}\le n_{2j+1}$ such that the
following conditions hold:
\begin{enumerate}
\item[TP.1] For all $i\le \la_{\Phi,\Psi}$, $w(\phi_i)=w(\psi_i)$ and
$p_i^{\Phi}=p_i^{\Psi}$.
\item[TP.2] For all $i< \ka_{\Phi,\Psi}$, $\phi_i=\psi_i$.
\item[TP.3]  For all   $\ka_{\Phi,\Psi}<i< \la_{\Phi,\Psi}$
\begin{align*}
&\supp \phi_i \cap \cl{\supp \psi_1\cup \dots \cup\supp\psi_{\la_{\Phi,\Psi}-1}}
{p_{\la_{\Phi,\Psi}-1}}=
\buit \\
\text{and }&\supp \psi_i
 \cap \cl{\supp \phi_1\cup \dots \cup \supp\phi_{\la_{\Phi,\Psi}-1}}{p_{\la_{\Phi,\Psi}-1}}=\buit.
\end{align*}
\item[TP.4] $\conj{w(\phi_i)}{\la_{\Phi,\Psi}<i\le n_{2j+1}} \cap \conj{w(\psi_i)}{i\le n_{2j+1}}=\buit$
  and $\conj{w(\psi_i)}{\la_{\Phi,\Psi} <i\le n_{2j+1}} \cap \conj{w(\phi_i)}{i\le n_{2j+1}}=\buit.$
\end{enumerate}

\begin{center}
\begin{figure}[h]
\includegraphics[scale=0.8]{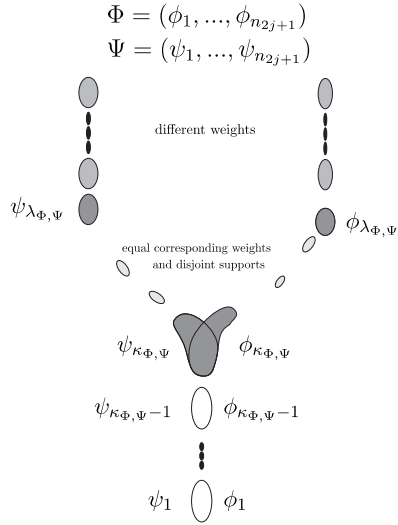}
\caption{The tree-like  interference between a pair of special sequences.}\label{figu1}
\end{figure}
\end{center}

\flema We refer to the reader to Figure \ref{figu1} for a description of the conclusion
of this Lemma. \prue First we observe that for $i\neq l$, $w(\phi_i)\neq w(\phi_l)$.
Indeed if $i=1$ and $l>1$ then $w(\phi_1)=m_{2j_1}$ with $j_1$ even while
$w(\psi_l)=m_{2j_l'}$ with $j_l'$ odd. If $i$ and $l$ are greater than $1$ then
$w(\phi_i)=m_{2j_i}\neq m_{2j_l'}=w(\psi_l)$  as consequence of the fact that they code
sequences of different lengths $i-1$ and $l-1$ respectively.

 Let $\la_{\Phi,\Psi}$ be the maximum of all $i\le n_{2j+1}$ such
that $w(\phi_i)=w(\psi_i)$ if defined. If not, we set
$\la_{\Phi,\Psi}=\ka_{\Phi,\Psi}=0$.
%
%
Suppose now that $\la_{\Phi,\Psi}>0$. Define $\ka_{\Phi,\Psi}$ by
\begin{align*}
\ka_{\Phi,\Psi}=\min\conj{i<\la_{\Phi,\Psi}}{\phi_i\neq \psi_i},
\end{align*}
if defined and $\ka_{\Phi,\Psi}=0$ if not. In this last case it is trivial to check our
requirements. So assume that $\ka_{\Phi,\Psi}>0$. (TP.2) and (TP.4) follows easily from
the properties of the coding $\sig_\ro$. We show  (TP.3). Let
\begin{align*}
G= \cl{\bigcup_{i=1}^{\la_{\Phi,\Psi}-1}\supp \phi_i}{p_{\la_{\Phi,\Psi}-1}} \text{ and } G'=
\cl{\bigcup_{i=1}^{\la_{\Phi,\Psi}-1}\supp \psi_i}{p_{\la_{\Phi,\Psi}-1}}.
\end{align*}
And let   $\pi_G:G\to \{1,\dots,\#G\}$ and $\pi_G':G'\to \{1,\dots,\#G\}$ be the unique
order-preserving bijections. \clam \noindent (a)  $\#G=\#G'$.

\noindent (b) $\pi_G |(G\cap G')=\pi_{G'}|(G\cap G')$ and $(G\cap
G')\phi_{\ka_{\Phi,\Psi}}=(G\cap G')\psi_{\ka_{\Phi,\Psi}}$.

\noindent (c) $\max (G\cap G')<\min \{\max \supp \phi_{\ka_{\Phi,\Psi}},\max \supp
\psi_{\ka_{\Phi,\Psi}}\}$. \fclam \prucl (a): Notice that
\begin{equation}
\text{$\# G=\max \supp \pi_G(\phi_{\la_{\Phi,\Psi}-1})$ and  $\# G'=\max \supp
\pi_{G'}(\psi_{\la_{\Phi,\Psi}-1})$.}
\end{equation}  Since
$\sig_\ro((\phi_i,w(\phi_i),p_i)_{i=1}^{\la_{\Phi,\Psi}-1} )=
\sig_\ro((\psi_i,w(\psi_i),p_i)_{i=1}^{\la_{\Phi,\Psi}-1} )$, then
$\pi_G(\phi_{\la_{\Phi,\Psi}-1})=\pi_{G'}(\psi_{\la_{\Phi,\Psi}-1})$ and hence $\#G=\#G'$, as
desired.  (b): It follows from the properties of $\ro$ that $\pi_G |(G\cap G')=\pi_{G'}|(G\cap
G')$ (see Remark \ref{roclosrem} and Proposition \ref{easyro}). Fix now $\al\in G\cap G'$.
Since $\pi_G(\al)=\pi_{G'}(\al)$ we have that
\begin{equation}
\phi_{\ka_{\Phi,\Psi}}(e_\al)=\psi_{\ka_{\Phi,\Psi}}(e_{\pi_G(\pi_{G'}^{-1}\al)})=
\psi_{\ka_{\Phi,\Psi}}(e_{\al}),
\end{equation}
as desired. (c): Suppose not.   W.l.o.g.  assume that $\max G\cap G'\ge \max \supp
\phi_{\ka_{\Phi,\Psi}}$. (b) yields that
\begin{equation}
 \phi_{\ka_{\Phi,\Psi}}=(G\cap G' ) \phi_{\ka_{\Phi,\Psi}}= (G\cap G')\psi_{\ka_{\Phi,\Psi}},
\end{equation}
and since $\#\supp \phi_{\ka_{\Phi,\Psi}}= \#\supp \psi_{\ka_{\Phi,\Psi}}$ we obtain that
$\phi_{\ka_{\Phi,\Psi}}=\psi_{\ka_{\Phi,\Psi}}$, a contradiction. \fprucl To complete the
proof choose $\ka_{\Phi,\Psi}<i< \la_{\Phi,\Psi}$. Then the previous Claim yields that $\supp
\phi_i\con G\setminus (G\cap G')$ and hence $\supp\phi_i\cap G'=\buit$. \fprue

\subsection{Tree-analysis of functionals}
Computing the norm of a vector $x$ from $\eqs_{\ou}$ is typically not an easy task. From the
definition of the norming set $K_{\ou}$ one observes that each $\phi\in K_{\ou}$ is
constructed from the basic functionals $e_\al^*$ after finitely many steps where at each step
one applies some $( m_j^{-1},{n_j})$-operation, or one takes some convex combination. The
tree-analysis defined below describes this procedure and it will be a very useful tool in
estimations of the norm of certain vectors of $\eqs_{\ou}$.

\defi
Let $\phi$ be a functional of the norming set $K_{\ou}$. A
\emph{tree-analysis} of $\phi$ is a mapping $ \mc F :\mc T  \to
K_{\ou}$, $ t \mapsto \mc F(t)=\phi_t$ such that the following
conditions are satisfied:
\begin{enumerate}
\item[1.]  $\mc T=(\mc T, \prec)$ is a finite tree with a unique
root $\buit\in\mc{T}$, and   $\phi_\buit=\phi$.
 \item[2.] If $t\in\mc{T}$ is a maximal node of $\mc{T}$, then $\phi_t=\pm e_\al^*$,
 for some  $\al<\ou$. We say in this case that $\phi_t=\pm e_\al^*$ has \emph{type 0}.
 \item[3.] If $t\in \mc{T}$ is not a  maximal node, and denoting by $S_t$ the set of
 immediate
 successors of  $t$, $S_t$  satisfies exactly  one of the following  two:
 \begin{enumerate}
 \item[(3.a)] There is a unique ordering of $S_t=\{s_1<_t \dots <_t s_d\}$  defined by
 $\phi_{s_1}<\dots<\phi_{s_d}$, there exists an integer $j\in\mathbb{N}$ such that $d\le n_j$ and
 $\phi_t=({1}/{m_j})\sum_{i=1}^d \phi_{s_i}$.
 \item[(3.b)] For every $s\in S_t$, 
 $\phi_s$ is either of type 0 or I, and   there is a sub-convex family   $\{r_u\}_{u\in S_t}$ of positive
  rational numbers such that
$\phi_t=\sum_{s\in S_t}r_s \phi_s$.
 \end{enumerate}
 \item[4.] For every $s\pe t$, $\ran \phi_t\con \ran \phi_s$.
 \end{enumerate}
 For a given $\phi\in K_{\ou}$, whenever we write $w(\phi)$ we implicitly assume that $\phi$
 is
 of type I. In many cases we will use the explicit notation $(\phi_t)_{t\in \mc T}$ to
 denote a tree-analysis.
 \fdefi

\nota \label{CnormsD} \noindent 1. The minimality of $K_{\ou}$ easily yields that for any
functional $f\in K_{\ou}$ there is a tree $(f_t)_{t\in \mc T}$ satisfying conditions 1-3.
Such a tree $(f_t)_{t\in \mc T}$ for $f$ can be refined to a tree-analysis of $f$. The
proof goes as follows: Given a tree-analysis $(f_t)_{t\in \mc T}$ of  $f$ we show by
downwards induction over $\mc T$ that every $f_t$ has a tree-analysis as desired. The
only non trivial case is when $f_t$ is of type II, $f_t=\sum_{s\in S_t}r_s f_s$. Let
$E=\ran f_t$, and let $f_s'=f_s|E$ for $s\in S_t$. Then, $f_t= f_t|E=\left(\sum_{s\in
S_t}r_s f_s\right)|E=\sum_{s\in S_t}r_s f_s'$. Since $\ran f_s'\con E=\ran f_t$, the
inductive hypothesis finishes the proof.

\noindent 2. Observe that the subset of $K_{\ou}$ consisting on functionals of type $0$
and $I$ is also a
 norming set for the space: Given a finitely supported vector   $x$, and $\phi=\sum_i r_i \phi_i$ of type II with $\phi_i$
 of type 0 or I,  $|\langle \phi,x \rangle|=|\sum_i r_i\langle \phi_i,x \rangle|\le \max_i
  |\langle \phi_i,x \rangle|$.

\noindent 3. Observe that a given $\phi\in K_{\ou}$ may have many trees as well as
weights. \fnota

\section{ $\eqs_{\ou}$ has no unconditional basic sequences}

\defi
A pair $(x,\phi)$ with $x\in \eqs_{\ou}$ and $\phi\in K_{\ou}$ is said to be a
${(C,j)}$\emph{-exact pair} if \emph{(a)} $\nrm{x}\le C$, $w(\phi)= m_j$  and $\phi(x)=1$, and
\emph{(b)} for every $\psi\in K_{\ou}$ of type I and $w(\psi)=m_{i}$, $i\neq j$ we have:
\begin{equation}
|\psi(x)|\le \left\{\begin{array}{ll} \frac {2C}{m_i} & \text{if }i <j \\
 \frac C{m_j^2}& \text{if } i>j.
\end{array}\right.
\end{equation}
\fdefi $(C,j)$-exact pairs are one of the basic ingredients for the study of mixed
Tsirelson spaces as well as of hereditarily indecomposable  spaces built on a frame of a
mixed Tsirelson space. The next proposition ensures their existence everywhere. \prop Let
$(x_n)_n$ be a block sequence in $\eqs_{\ou}$. Then for each $j\in \N$ there exists
$(x,\phi)$ such that $x\in \langle  x_n \rangle_n$, $\phi\in K_{\ou}$ and $(x,\phi)$ is a
$(6,j)$-exact pair. \fprop

The existence of $(6,j)$-exact pairs it is proved by a similar argument to that for the
Gowers-Maurey space \cite{go-ma}. It is primarily based on the unconditional part of the
definition of $K_{\ou}$ (i.e., property 2.). A simple example of a $(6,2j)$-exact pair is
the pair $(x,\phi)$ where $x=(m_{2j}/n_{2j})\sum_{\al\in  F}e_\al$,
$\phi=(1/m_{2j})\sum_{\al\in F}e_{\al}^*$ and $\#F=n_{2j}$.

\defi\label{1dep}
Let $j\in \N$. A sequence $(x_1,\phi_1,\dots,x_{n_{2j+1}},\phi_{n_{2j+1}})$ is said to be a
\emph{$(1,j)$-dependent sequence} if:
\begin{enumerate}
\item[\emph{DS.1}] $\supp x_1\cup \supp \phi_{1}<\dots <\supp x_{n_{2j+1}}\cup \supp
\phi_{n_{2j+1}}$. \item[\emph{DS.2}] The sequence $\Phi=(\phi_1,\dots,\phi_{n_{2j+1}})$
is a $2j+1$-special sequence. \item[\emph{DS.3}] $(x_i,\phi_i)$ is a $(6,2j_i)$-exact
pair  with $\#\supp x_i \le m_{2j_{i+1}}/n_{2j+1}^2 $ for every $1\le i \le n_{2j+1}$.
\item[\emph{DS.4}]  For every $(2j+1)$-special
sequence  $\Psi=(\psi_1,\dots,\psi_{n_{2j+1}})$ we have that
\begin{equation}
\bigcup_{\ka_{\Phi,\Psi}<i<\la_{\Phi,\Psi}}\supp x_i \cap
\bigcup_{\ka_{\Phi,\Psi}<i<\la_{\Phi,\Psi}} \supp \psi_i =\buit.
\end{equation}
\end{enumerate}

\fdefi \prop\label{depseqex} For every $(y_n)_n$, a block sequence of
$\eqs_{\ou}$, and every $j\in \N$ there exists $(1,j)$-dependent sequence
$(x_1,\phi_1,\dots,x_{n_{2j+1}},\phi_{n_{2j+1}})$  such that $x_i\in \langle y_n \rangle_n$
for every $i=1,\dots,n_{2j+1}$. \fprop \prue Let $(y_n)_n$ and $j$ be given. We inductively
produce $\{(x_i,\phi_i)\}_{i=1}^{n_{2j+1}}$ as follows. For $i=1$ we choose a $(6,2j_1)$-exact
pair $(x_1,\phi_1)$ such that  $m_{2j_1}>m_{2j+1}^2$, $j_1$ even (see the definition of
special sequences) and $x_1\in \langle y_n \rangle_n$. Assume that
$\{(x_l,\phi_l)\}_{l=1}^{i-1}$ has been chosen such that there exists $(p_l)_{l=1}^{i-2}$
satisfying

\noindent \emph{(a)} $\supp x_1\cup \supp \phi_1<\dots <\supp x_{i-1}\cup \supp
\phi_{i-1}$, each $x_l\in \langle y_n \rangle_n$ and $(x_l,\phi_l)$ is a $(6,2j_l)$-exact
pair.

\noindent \emph{(b)} For $1<l\le i-1 $,
$w(\phi_l)=\sig_\ro(\phi_1,w(\phi_1),p_1,\dots.,\phi_{l-1},w(\phi_{l-1}), p_{l-1} ) $.

\noindent \emph{(c)} For $1\le l < i-1 $, $p_{l}\ge \max\{p_{l-1},p_{F_{l}}\}$, where
$F_l=\bigcup_{k=1}^l \supp \phi_k\cup \supp x_k $.

To define $(x_i,\phi_i)$ we choose $p_{i-1}\ge \max\{p_{i-2},p_{F_{i-1}},n_{2j+1}^2\cdot  \#\supp x_i \}$ and
we set
\begin{align*}
2j_i=\sig_\ro(\phi_1,w(\phi_1),p_1,\dots.,\phi_{i-1},w(\phi_{i-1}),p_{i-1}).
\end{align*}
Choose a $(6,2j_i)$-exact pair $(x_i,\phi_i)$ such that $x_i\in \langle  y_n \rangle_n$
and $\supp x_{i-1}\cup \supp \phi_{i-1}<\supp x_i\cup \supp \phi_i$. This completes the
inductive construction. (\emph{DS.1})-(\emph{DS.3})  easily holds, while (\emph{DS.4})
follows from \emph{(c)} and (3) of Lemma \ref{0011}.
\fprue

\nota Suppose that $(y_n)_n$ and $(z_n)_n$ are  block sequences such that $\sup_n  \max
\supp  y_n = \sup_n \max    \supp z_n$. Then for every $j\in \N$ there is a
$(1,j)$-dependent sequence $(x_1,\phi_1,\dots,x_{n_{2j+1}},\linebreak \phi_{n_{2j+1}})$
with the property that $x_{2i-1}\in \langle y_n \rangle_n $ and $x_{2i}\in \langle z_n
\rangle_n$ for every $i=1,\dots,n_{2j+1}/2$. \fnota

\lema\label{depseqest} Fix a $(1,j)$-dependent sequence
$(x_1,\phi_1,\dots,x_{n_{2j+1}},\phi_{n_{2j+1}})$, and a sequence of scalars
$(\la_i)_{i=1}^{n_{2j+1}}$
 such that $\max_i|\la_i|\le 1$. Suppose that for every   $\psi\in K_{\ou}$ such that
$w(\psi)=m_{2j+1}$, and every interval of integers  $E\con [1,n_{2j+1}]$  it holds that
\begin{equation}\label{jkioe1}
|\psi(\sum_{i\in E}\la_i x_i)|\le 12(1+\frac{\#E}{n_{2j+1}^2}).
\end{equation}
Then,
\begin{equation}\label{jkioe2}
\nrm{\frac1{n_{2j+1}}\sum_{i=1}^{n_{2j+1}}\la_ix_i}\le  \frac{1}{m_{2j+1}^2}.
\end{equation}
\flema
We postpone the proof of this lemma to the end of Subsection \ref{cbi}, since involves non
trivial estimates.

\prop
If  $(x_1,\phi_1,\dots,x_{n_{2j+1}},\phi_{n_{2j+1}})$ is a $(1,j)$-dependent sequence, then
\begin{equation}
\nrm{\frac1{n_{2j+1}}\sum_{i=1}^{n_{2j+1}}x_i}\ge \frac1{m_{2j+1}} \text{ and }
\nrm{\frac1{n_{2j+1}}\sum_{i=1}^{n_{2j+1}}(-1)^{i+1}x_i}\le \frac{1}{m_{2j+1}^2}.
\end{equation}
\fprop \prue The first estimate is clear since the functional
$\phi=(1/m_{2j+1})\sum_{i=1}^{n_{2j+1}}\phi_i\in K_{\ou}$ and $\phi((1/n_{2j+1})
\sum_{i=1}^{n_{2j+1}}x_i)=1/m_{2j+1}$. For the second, we use Lemma \ref{depseqest} applied to
the sequence of scalars $((-1)^{i+1})_i$, and the desired estimate will follow from
(\ref{jkioe2}). Fix $\psi \in K_{\ou}$ with $w(\psi)=m_{2j+1}$, and  an interval $E\con
[1,n_{2j+1}]$. Set $\Psi=(\psi_1,\dots,\psi_{n_{2j+1}})$ and $x=\sum_{i\in E}(-1)^{i+1}x_i$,
where $\psi=(1/m_{2j+1})\sum_{i\in E}\psi_i$. Notice that
\begin{equation}\label{eeirue1}
|\psi(x)|= |\frac1{m_{2j+1}}\sum_{i=1}^{\ka_{\Phi,\Psi}-1}\phi_i(x) +
\frac1{m_{2j+1}}\sum_{i=\ka_{\Phi,\Psi}}^{n_{2j+1}}\psi_i(x)|\le \frac1{m_{2j+1}}+
|\frac1{m_{2j+1}}\sum_{i=\ka_{\Phi,\Psi}}^{n_{2j+1}}\psi_i(x)|.
\end{equation}
We shall show that the following hold
\begin{enumerate}
\item[\emph{(a)}] $|\psi_{\ka_{\Phi,\Psi}}(\sum_{i\in E}(-1)^{i+1}x_i)|\le 1+
12(\#E-1)/n_{2j+1}^2 $, \item[\emph{(b)}] $|\psi_{\la_{\Phi,\Psi}}(\sum_{i\in
E}(-1)^{i+1}x_i)|\le  1+ 12(\#E-1)/n_{2j+1}^2 $, and \item[\emph{(c)}]
$|(\sum_{l>\ka_{\Phi,\Psi},l\neq \la_{\Phi,\Psi}}\psi_l)(x_i)|\le 12/n_{2j+1} $ for every
$1\le i\le n_{2j+1}$.
\end{enumerate}
Let us show first \emph{(a)}. Let $2j_i$ be such that $w(\phi_i)=m_{2j_i}$. Notice that for $i\neq
\ka_{\Phi,\Psi}$ we have that
\begin{equation}\label{jyew1}
|\psi_{\ka_{\Phi,\Psi}}(x_i)|\le \left\{ \begin{array}{ll}\frac {12}{w(\psi_{\ka_{\Phi,\Psi}})}
 & \text{if } i>\ka_{\Phi,\Psi} \\
\frac {6}{ m_{2j_i}^2} & \text{if } i<\ka_{\Phi,\Psi}.
\end{array}
\right.
\end{equation}
By the properties of the sequences $(m_l)_l$, $(n_l)_l$ and the fact that
$n_{2j+1}^2<w(\psi_{\ka_{\Phi,\Psi}}), \,m_{2j_i}$, (\ref{jyew1}) yields that
$|\psi_{\ka_{\Phi,\Psi}}(x_i)|\le {12}/{n_{2j+1}^2}$ for  $i\neq \ka_{\Phi,\Psi}$. Hence
\begin{equation}
|\psi_{\ka_{\Phi,\Psi}}(\sum_{i\in E}x_i)|\le |\psi_{\ka_{\Phi,\Psi}}(x_{\ka_{\Phi,\Psi}})|+
|\psi_{\ka_{\Phi,\Psi}}(\sum_{i\in E,\, i\neq \ka_{\Phi,\Psi}}x_i)|\le 1+ \frac{12
(\#E-1)}{n_{2j+1}^2}.
\end{equation}
\noindent\emph{(b)} has a similar proof to that of \emph{(a)}. We check now \emph{(c)}. Fix
$l>\ka_{\Phi,\Psi}$, $l\neq \la_{\Phi,\Psi}$.  Suppose that $l>\la_{\Phi,\Psi}$. Since
$w(\psi_{l})\neq w(\phi_i)$ for all $i\le n_{2j+1}$, we obtain that $|\psi_{l}(x_i)|\le
{12}/{n_{2j+1}^2}$. Now suppose that $\ka_{\Phi,\Psi}<l<\la_{\Phi,\Psi}$. By (\emph{DS.4}) we
have that $\psi_{l}(x_i)=0$ for every $\ka_{\Phi,\Psi}<i<\la_{\Phi,\Psi}$. And for  $i\notin
(\ka_{\Phi,\Psi},\la_{\Phi,\Psi})$, using the fact that $w(\psi_l)\neq w(\phi_i)$, we can
conclude that $|\psi_{l}(x_i)|\le {12}/{n_{2j+1}^2}$. Hence, $(\sum_{l>\ka_{\Phi,\Psi},l\neq
\la_{\Phi,\Psi}}\psi_l)(x_i)\le {12}/{n_{2j+1}}$ for every $1\le i\le n_{2j+1}$, as desired.

Combining \emph{(a)}, \emph{(b)} and\emph{ (c)} we obtain that
\begin{equation}\label{eeirue2}
|\frac1{m_{2j+1}}\sum_{i=\ka_{\Phi,\Psi}}^{n_{2j+1}}\psi_i(x)|\le  1 +\frac {\#E}{n_{2j+1}^2}.
\end{equation}
From (\ref{eeirue1}) and (\ref{eeirue2}) we conclude that $|\psi(x)|\le
12(1+{\#E}/{n_{2j+1}^2})$, as desired. \fprue

\prop\label{blosunc} The closed linear span of a block sequence of $\eqs_{\ou}$ is
hereditarily indecomposable. \fprop \prue Fix a block sequence $(y_n)_n$ of $\eqs_{\ou}$,
two block subsequences $(z_n)_n$ and $(w_n)_n$ of $(y_n)_n$ and $\vep>0$. Let $j$ be
large enough such that $m_{2j+1}\vep
>1$. By Proposition \ref{depseqex} we can choose a $(1,j)$-dependent sequence $(x_1,\phi_1,\dots,x_{n_{2j+1}},\phi_{n_{2j+1}})$
such that  $x_{2i-1}\in \langle z_n \rangle_n$, and $x_{2i}\in \langle w_n \rangle_n$. Set
$z=(1/n_{2j+1})\sum_{i=1,i \text{ odd}}^{n_{2j+1}}x_i$ and $w=(1/n_{2j+1})\sum_{i=1,i \text{
even}}^{n_{2j+1}}x_i$. Notice that $z\in \langle z_n \rangle_n$ and $w\in \langle w_n \rangle_n$.
By Proposition \ref{depseqest}, we know that $\nrm{z+w}\ge {1}/{m_{2j+1}}$   and $ \nrm{z-w}\le
{1}/{m_{2j+1}^2}$. Hence $\nrm{z-w}\le \vep\nrm{z+w}$.
\fprue
\cor \label{dis0} \emph{(a)} The distance between the unit spheres of every two normalized block
sequences $(x_n)$ and $(y_n)$ in $\eqs_{\ou}$ such that $\sup_n \max \supp x_n= \sup_n \max \supp
y_n$ is 0.

\noindent \emph{(b)} There is no unconditional basic sequence in $\eqs_{\ou}$.

\noindent \emph{(c)} Every infinite dimensional closed subspace of $\eqs_{\ou}$ contains an
hereditarily indecomposable subspace.

\noindent \emph{(d)} The distance between the unit spheres of  two nonseparable subspaces of
$\eqs_{\ou}$ is equal to 0. \fcor
\prue
(b): follows from Proposition \ref{blosunc} and 4. of Proposition \ref{gliding}.
 (c): This result follows from (b) and Gowers' dichotomy. Moreover, every subspace
of $\eqs_{\ou}$ isomorphic to the closed linear span of a block sequence with respect to the
basis $(e_\al)_{\al<\ou}$ is hereditarily indecomposable. (d): Fix two nonseparable closed
subspaces $X$ and $Y$ of $\eqs_{\ou}$.  Now we can find a sequence $(z_n)_n$ of normalized
vectors such that for every $n$ (a) $z_{2n-1}\in X $, $z_{2n}\in Y$ and (b) $\supp z_{n}<\supp
z_{n+1}$.  Notice that the supports $\supp z_n$ are not necessarily finite. Now approximating
$(z_n)_n$ by a normalized block sequence $(w_n)_n$ as close as needed we obtain the desired
result. \fprue

\section{Basic estimations and further properties of ${\eqs_{\ou}}$}
In this section we introduce some of the standard tools of this
area (see \cite{schlu}, \cite{go-ma},\cite{arg-tol},
\cite{arg-del-man}) which will be quite useful in our analysis of
the space $\eqs_{\ou}$. We also obtain that the space $\eqs_{\ou}$
is reflexive.

\subsection{Rapidly Increasing Sequences. The basic inequality}

\defi \emph{(Rapidly Increasing Sequences (RIS))}
Let $C,\vep>0$. A block sequence  $(x_k)_k$ of $X$ is called  a $(C,\vep)$- \emph{rapidly
increasing sequence} ($(C,\vep)$-RIS in short) iff there is an increasing sequence $(j_k)_k$
of integers such that for all $k$,
\begin{enumerate}
\item $\nrm{x_k}\le C$
\item $|\supp x_k| \le m_{j_{k+1}}\vep$ and
\item For all type I functionals $\phi$ of $K$ with $w(\phi)<m_{j_k}$, $|\phi(x_k)|\le C/w(\phi)$.
\end{enumerate}
\fdefi

\nota\label{RRIS} \noindent \emph{1.} Notice that given
$\vep'<\vep$,  every $(C,\vep)$-RIS  has a subsequence which is $(C,\vep')$-RIS. Notice also
that for every strictly increasing sequence $\{\al_n\}_n $, and every $\vep>0$,
$(e_{\al_n})_n$ is a $(1,\vep)$-RIS.  \emph{2.} For every $(1,j)$-dependent sequence
$(x_1,\phi_1,\dots,x_{n_{2j+1}},\phi_{n_{2j+1}})$   the corresponding  block sequence
$(x_1,\dots,x_{n_{2j+1}})$ is a $(12,1/n_{2j+1}^2)$-RIS. \fnota A primary reason for the
usefulness of the notion of RIS is that one has good estimates of upper bounds on $|\langle
f,x \rangle|$, for $f\in K_{\ou}$ and $x$ averages of an RIS.

\noindent{\bf Notation.} In the sequel we shall denote  by $W$  the minimal subset of
$c_{00}(\N)$ which contains $\{e_n^*\}_{n\in \N}$, is symmetric, and is closed in rational
convex combinations, closed in restriction to intervals, and closed for the $(
m_{j}^{-1},{4n_{j}})$-operations.
 \nota\label{oipwer1} By minimality of $W$, every
element $f$ of $W$ has  a
 tree-analysis $(f_t)_{t\in \mc T}$. Using induction over the tree-analysis, it is not difficult to show that
 every $f\in K$ is the convex combination $f=\sum_i r_i f_i$, with every $f_i$ in the norming set of $T[( m_j^{-1},
 {4n_j})_j]$
 and in the case that $f$ is of type I, then each $f_i$ can be chosen such that $w(f_i)=w(f)$.
 Hence, $W$ norms the mixed Tsirelson  space $T[( m_j^{-1},{4n_j})_j]$.
\fnota


The following Lemma gives a very useful tool for reducing for a given $f\in K_{\ou}$ and a RIS
$(x_k)_k$, upper bound estimates of $|\langle f,\sum_k b_k x_k \rangle|$ to upper bounds of
$|\langle g,\sum_k |b_k| e_k \rangle|$ where $g$ is a functional of the auxiliary space
$T[(m_j^{-1}, {4n_j})_j]$ and $(e_k)_k$ is its basis.

\lema[\emph{Basic Inequality for  RIS}]\label{bin} Let $(x_n)_n$ be a   $(C,\vep)$-RIS
sequence and fix $(b_k)_k\in c_{00}(\N)$. Suppose that $j_0\in \N$ is such that for all $f\in
K_{\ou}$  with weight $w(f)=m_{j_0}$ and all intervals $E$,
\begin{equation}  \label{negligible}
\left|f(\sum_{k\in E} b_k x_k)\right|\le C\left(\max_{k\in E}|b_k|+\vep\sum_{k\in E}
 |b_k|\right).
\end{equation}
(We say in this case that  $(x_n)_n$ makes ${j_0}$ \emph{negligible} for $(b_k)_k$.) Then for
every $f\in K_{\ou}$ of type I there exists $g_1, g_2\in c_{00}(\N)$ such that

$$   |f(\sum b_k x_k)|\le C(g_1+g_2)(\sum |b_k|e_k),$$
where $g_1=h_1$ or $g_1=e_t^*+h_1$, $t\notin \supp h_1$, and $ h_1\in W$ is such that $h_1\in
\conv_\Q\conj{h\in W}{w(h)=w(f)}$ and  with $m_{j_0}$ not appearing  as a weight of  a node of a
tree-analysis of $h_1$, and
 $\nrm{g_2}_\infty\le \vep$.

\flema
We postpone the proof of this result until Subsection \ref{dwaknvnbvb}.

\nota Notice that any  finite $(C,\vep)$-RIS sequence $(x_k)_k$ is going to be
$j_0$-negligible for large $j_0$. \fnota

\subsection{Estimates on the basis}

%
\prop\label{estbasis}
Fix a functional $f$ of type I,  either in $W$ or in $K_{\ou}$,  $j\in \N$ and $l\in
[n_j/m_j,n_j]$. Then for every set $\#F=l$
\begin{equation}
|f(\frac1{l}\sum_{\al\in F}e_\al)|\le \left\{ \begin{array}{ll} \frac 2{w(f)m_j} &
\text{if }w(f)<m_j \\
\frac 1{w(f)}& \text{if } w(f)\ge m_j.
 \end{array} \right.
\end{equation}
If the tree-analysis of $f$ does not  contain nodes with weight $m_j$,  then
\begin{equation}
|f(\frac1{l}\sum_{\al\in F} e_\al)|\le \frac2{m_j^3},
\end{equation}
where in each case we interpret $(e_\al)_{\al\in F}$ in the obvious way. \fprop \prue Fix
$f\in W$ of type I. By Remark  \ref{oipwer1} we can assume that $f$ belongs to the norming set
of $T[( m_j^{-1},{4n_j})_j]$, i.e., $f$ admits a tree-analysis with no convex combinations.
The result is proved in the same manner as Lemma 4.2 of \cite{arg-mano}.

The result for $f\in K_{\ou}$ follows easily from the following. Let us denote by $\nrm{\cdot}_l$ the
norm of the natural extension of $T[( m_j^{-1},{4n_j})_j]$ to $\ou$. It is clear that for this
norm the natural Hamel basis $(e_\al)_{\al<\ou}$ of $c_{00}(\ou)$ is 1-subsymmetric, and also that
$\nrm{\cdot}_l$ dominates the norm $\nrm{\cdot}_{\eqs_{\ou}}$.
\fprue

%

\subsection{Consequences of the basic inequality }\label{cbi}

We start this subsection with the following estimates on averages of RIS.
\prop\label{niceasris}
Let $(x_k)_{k}$ be a $( C,\vep)$-RIS for  $\vep\le 1/n_j$, $l\in [n_{j}/m_j,n_j]$ and let
$f\in K$ of type I. Then,
\begin{equation}
|f(\frac1{l}\sum_{k=1}^{l}x_k)|\le \left\{ \begin{array}{ll}
\frac {3C}{w(f)m_j} & \text{if }w(f)<m_j \\
\frac C{w(f)}+\frac{2C}{n_j} & \text{if } w(f)\ge m_j.
 \end{array} \right.
\end{equation}
Consequently, if $(x_k)_{k=1}^{l}$ is a normalized $(C,\vep)$-RIS with $\vep\le 1/n_{2j}$ and
$l\in [n_{2j}/m_{2j},n_{2j}]$, then
\begin{equation}
\frac 1{m_{2j}}\le \nrm{\frac1{l} \sum_{k=1}^{l} x_k }\le \frac {2C}{m_{2j}}.
\end{equation}
\fprop

\prue This follows from the basic inequality and the estimates on the basis of
$T[(m_j^{-1},{4n_j})_j]$ given in Proposition \ref{estbasis}. For the last consequence, notice
that if for every $k\le l$ we consider $x_k^*$ in $K$ such that $x_k^* x_k=1$ and $\ran
x_k^*\con \range x_k$, then $x^*=(1/m_{2j})\sum_{k=1}^{l}x_k^*$ belongs to $K$, and
$x^*((1/n_{2j})\sum_{k=1}^{l} x_k )=1/m_{2j}$. \fprue

\defi

Let $C>0$ and $ k\in \N$. A normalized vector $y$ is called a $C-\ell_1^k$-average iff there is a finite
block sequence $(x_1,\dots,x_n)$ such that
  $y=(x_1+\dots + x_k)/k$ and  $\nrm{x_i}\le C$.

Observe that since  $K_{\ou}$ is closed under the $(m_{2j}^{-1},{n_{2j}})$-operation, for
every normalized block sequences $(y_n)_n$ and every $k$, there are $z_1<\dots < z_k$ in
$\langle y_n \rangle_n$ such that $(z_1+\dots + z_k)/k$ is a
 $2-\ell_1^k$-average (for a detailed proof see for example \cite{arg-mano}).
\fdefi

\prop\label{risav}
Suppose that $y$ is a  $C-\ell_1^k$-average  and suppose that $E_1<\dots < E_n$ are
intervals with $n< k$. Then, $ \sum_{i=1}^n \nrm{E_i y}\le C(1+{2n}/k)$. As a
consequence, if $y$ is a  $C-\ell_1^{n_j}$-average  and    $\phi\in K$ is  with
$w(\phi)<m_j$ , then $|\phi(y)|\le {3C}/{2w(\phi)}$.

In particular, for $2-\ell_1^{n_j}$-averages we get that $|\phi(y)|\le 3/w(\phi)$ if
$w(\phi)<m_j$.

\fprop

\prue
See  \cite{schlu} or  \cite{go-ma}.
\fprue

\nota Suppose that $(x_k)_k$ is such that there is  a strictly
increasing sequence $(j_k)_k$ and $\vep>0$ such that for all $k$, (a) $x_k$ is a
$2-\ell_1^{n_{j_k}}$-average and (b)  $\#\supp x_k < \vep m_{j_{k+1}}$. Then Proposition
\ref{risav} shows that $(x_k)_k$ is a $(3,\vep)$-RIS. In this case we will say that $(x_k)_k$
is a $(3,\vep)$-RIS of $\ell_1$ averages. These remarks yield the following. \fnota

\prop Any block sequence in $\eqs_{\ou}$ has a further normalized block subsequence which
is a $(3,\vep)$-RIS. \qed \fprop

\prop
Let $(x_n)_n$ be a block sequence in $\eqs_{\ou}$. Then for each $j\in \N$ there exists a
$(6,2j)$-exact pair $(x,\phi)$ such that $x\in \langle  x_n \rangle_n$.
\fprop
\prue Fix a block sequence $(x_n)_n$ of $\eqs_{\ou}$ and an integer $j$. By the previous
proposition we can find a normalized  $(3,1/n_{2j})$-RIS $(y_n)_n$ in $\langle x_n \rangle_n$.
For each $1\le i\le n_{2j}$ choose $\phi_i\in K_{\ou}$ such that $\phi_i(y_i)=1$, and
$\phi_i<\phi_{i+1}$. Set $\phi=(1/m_{2j})\sum_{i=1}^{n_{2j}}\phi_i\in K_{\ou}$, and
$x=(m_{2j}/n_{2j})\sum_{i=1}^{n_{2j}}y_i$. Then $\phi(x)=1$ and  estimates in Proposition
\ref{niceasris} yield
\begin{equation}
|f(x)|\le \left\{ \begin{array}{ll}
\frac {9}{w(f)} & \text{if }w(f)<m_{2j} \\
\frac {3 m_{2j}}{w(f)}+\frac{6m_{2j}}{n_{2j}} & \text{if } w(f)\ge m_{2j},
 \end{array} \right.
\end{equation}
and $\nrm{x}\le 6$. Hence $(x,\phi)$ is a $(6,2j)$-exact pair.\fprue

To finish this subsection we  show Lemma \ref{depseqest}: \prue
(Of Lemma \ref{depseqest}.) Fix a $(1,j)$ dependent sequence
$(x_1,\phi_1,\dots,x_{n_{2j+1}},\phi_{n_{2j+1}})$ and a sequence
$(\la_i)_{i=1}^{n_{2j+1}}$ with $\max_i |\la_i|\le 1$ such that
for every $\psi$ with weight $m_{2j+1}$, and every interval $E\con
[1,n_{2j+1}]$,
\begin{equation}  \label{jkiode1}
|\psi(\sum_{i\in E}\la_i x_i)|\le 12(1+\frac{\#E}{n_{2j+1}^2}).
\end{equation}
Since $(x_i)_i$ is a $(12,1/n_{2j+1}^2)$-RIS (see Remark \ref{RRIS}), (\ref{jkiode1}) tells
that $(x_i)_i$ makes $2j+1$ negligible for $(\la_i)_i$. From the conclusion of the basic
inequality and the estimates on the basis of $T( {4n_j},1/m_j)$, it follows  that for every
$f\in K_{\ou}$
\begin{equation}
|f(\frac1{n_{2j+1}}\sum_{i=1}^{n_{2j+1}} \la_i x_i)|\le 12
(\frac2{n_{2j+1}}+ \frac2{m_{2j+1}^3})\le \frac 1{m_{2j+1}^2},
\end{equation}
as required.
\fprue


\prop The basis $(e_\al)_{\al<\ou}$ is shrinking and boundedly
complete. Therefore $\eqs_{\ou}$ is reflexive. \fprop \prue Since the basis
$(e_\al)_{\al<\ou}$ is boundedly complete (see Remark \ref{basbound}), we only need to prove
that it is also shrinking. Suppose not. Then there exists a strictly increasing sequence
$A=\{\al_n\}_n$ of ordinals, scalars $(\la_n)_n$ and $x^*=w^*-\lim_n \sum_{n=1}^\infty \la_n
e_{\al_n}^*$ with $x^*\notin \overline{\langle e_{\al_n}^* \rangle_n}$. Thus there exist
$\vep>0$ and successive intervals $(E_n)_n$ such that for all $n$, $\nrm{E_n x^*}>\vep$.
Choose $(x_n)_n$ in $\eqs_A$ with $\supp x_n\con E_n$, $\nrm{x_n}=1$ and $x^*(x_n)>\vep$ for
all $n$. It follows that every convex combination $\sum_n \mu_n x_n$ satisfies $\nrm{\sum_n
\mu_n x_n}>\vep$. Now for sufficient large $j\in \N$ we may construct a $(2\vep,
1/n_{2j+1})$-RIS $(y_n)_n$ of $\vep$-normalized averages such that every $y_n$ is an average
of $(x_k)_k$. Proposition \ref{niceasris} yields that $\nrm{1/({n_{2j}})\sum_{i=1}^{n_{2j}}
y_i }\le (4\vep)/{m_{2j}}<\vep$, a contradiction. \fprue

\section{The operator spaces}
In this section we state and prove the main results about
operators on $\eqs_{\ou}$ and its subspaces. The new basic  tool
is the finite interval representability of a James-like space into
$\eqs_{\ou}$. The section is divided into six subsections. The
first  concern James like spaces. In the second  the finite
interval block representability of $J_{T_0}$ is defined and the
structure of the space of the step diagonal operators is studied.
In the third subsection the spaces $\mathcal{L}(\eqs_{\gamma})$
are studied and some consequences concerning the structure of the
subspaces of $\eqs_{\omega_1}$ are obtained. In the fourth
subsection the concept of asymptotically equivalent subspaces of
$\eqs_{\omega_1}$ is introduced and the structure of the spaces
$\mathcal{L}(X,\eqs_{\omega_1})$ with $X$ subspace of
$\eqs_{\omega_1}$ is studied. In the fifth subsection a
construction of  subspaces $X$ is presented such that $\dim
\mathcal{L}(X)/\mathcal{S}(X)=1$ while
$\mathcal{L}(X,\eqs_{\omega_1})/\mathcal{S}(X,\eqs_{\omega_1})$ is
of infinite dimension. The last subsection concerns some further
results about the operators related to the Fredholm theory of
strictly singular operators.

\subsection{James-like spaces}\label{jameslike}
\defi
Let $X$ be a reflexive space with a 1-subsymmetric basis
$(x_n)_n$, and let $A$ be a set of ordinals. $J_X(A)$ is the
completion of $(c_{00}(A),\nrm{\cdot}_{J_X(A)})$, where for $x\in
c_{00}(A)$,
$$  \nrm{x}_{J_X(A)}=\sup \conj{\nrm{\sum_{n=1}^l \left(\sum_{i\in I_n} x(i)\right)x_n}_X}{I_1<\dots
<I_n\text{ intervals of } A}.$$
 The natural Hamel basis $(v_\al)_{\al\in A}$ of $c_{00}(A)$ is a
bimonotone 1-subsymmetric transfinite basis of $J_{X}(A)$. Also,
for every interval $I$ of $A$ the functional $I^*:J_X(A)\to \R$,
$I^*(x)=\sum_{\al\in A}x(\al)$ belongs to $J_X^*(A)$ and
$\nrm{I^*}=1$. \fdefi

\nota As we shall see next, $\ell_1$ does not embed into $J_X(A)$
and hence the basis $(v_\al)_{\al\in A}$ is not unconditional. \fnota The following two facts
are easy extensions of the corresponding results from \cite{bel-hay-od}.
\prop\label{sdfewooo1} Let $(y_n)_n$ be a semi-normalized block sequence in $J_X(A)$ with
$\sum_{\al\in A}y_n(\al)=0$ for every $n$. Then $(y_n)_n$ is equivalent to the basis $(x_n)_n$
of $X$. \fprop \prue Let $0<c<C$ be  such that $c\le \nrm{y_n}\le C$ for all $n$. It is easy
to see that:
\begin{align*}
c\nrm{\sum_n a_n x_n}_X \le &     \nrm{\sum_n a_n y_n}_{J_X(A)}\le \sup_{i_1\le i_2\le
\dots \le i_l} \nrm{\sum_{q=1}^{l-1} (|a_{i_q}|
+|a_{i_{q+1}}|)x_q}_{X} \le \\
\le &(2C K)\nrm{\sum_n a_n x_n}_X,
\end{align*}
where $K$ is the unconditional constant of $(x_n)_n$. The first inequality holds for any block
sequence and the second uses our assumptions. \fprue \cor The space $\ell_1$ does not embed
into $J_X(A)$. \fcor \prue If not, then from Proposition \ref{gliding} we could find a
semi-normalized block sequence $(y_n)_n$ equivalent to the $\ell_1$-basis. Therefore, passing
if necessary to a further block sequence, we may assume that for all $n\in \N$,  $\sum_{\al\in
A}y_n(\al)=0$. Hence Proposition \ref{sdfewooo1} yields that $(y_n)_n$ is equivalent to
$(x_n)_n$, a contradiction. \fprue
\nota Suppose that $A$ and $B$ are two sets of ordinals with the same order type. Then
the unique order-preserving mapping $f:A\to B$ defines naturally an isometry between
$\widetilde{f}:J_{X}(A)\to J_{X}(B)$ by $\widetilde{f}(\sum_{\al\in H}r_\al v_\al)=
\sum_{\al \in H}r_{\al}v_{f(\al)}$. \fnota

The next proposition also extends the corresponding result from \cite{bel-hay-od}.

\prop\label{generby} For every ordinal $\gamma$ the space $J_X^*(\ga)$ is generated in
norm by $\{[0,\al)^*\}_{\al< \ga+1}$. \fprop \prue We proceed by induction. It is clear that
the successor ordinal case follows immediately from the inductive assumption. So we assume
that $\ga$ is limit ordinal and for all $\la<\ga$ the conclusion holds. Assume to the contrary
that $Y=\overline{\langle [0,\al)^* \rangle_{\al<\ga+1}}^{\nrm{\cdot}}\cones J_X^*(\ga)$. Then there
exists $x^*\in J_X^*(\ga)$ with $\nrm{x^*}=1$  and $\vep>0$ such that $d(x^*,Y)>\vep$. Observe
also that the inductive assumption yields that for all $\al<\ga$ if $x^*_\al$ denotes the
functional defined by
\begin{align*}
x_{\al}^*(v_\be)=\left\{\begin{array}{ll} 0 & \text{if } \be<\al \\
x^*(v_\be) & \text{if } \be \ge \al,
\end{array}
\right.
\end{align*}
then $\nrm{x_\al^*}\le 1$ and $d(x_\al^*,Y)> \vep$. In particular for all $\al<\ga$,
$d(x_\al^*,\langle [\al,\ga)^* \rangle)>\vep$ and from the Hahn-Banach and Goldstine Theorems
there exists a finitely supported $\tilde{y}_\al\in J_X(\ga)$ with $\nrm{\tilde{y}_\al}\le 1$,
$\al\le \min\supp \tilde{y}_\al$, $x^*(\tilde{y}_\al)>\vep$ and
$|\sum_{\be<\ga}\tilde{y}_\al(\be)|\le \vep/4$. Assuming further that $\al$ is a successor
ordinal we consider the vector $y_\al=\tilde{y}_\al -(\sum_{\be\ge
\al}\tilde{y}_\al(\be))v_{\al^-}$. Observe that $\al^-\le \min\supp y_\al$, $x^*(y_{\al})>\vep
-\vep/4>\vep/2$ and $\sum_{\be<\ga}y_\al(\be)=0$.
 Hence we may inductively choose a block
sequence $(z_n)_n$ such that $\vep/2\le \nrm{z_n}\le 1$, $\sum_{\al<\ga}z_n(\al)=0$ and
$x^*(z_n)>\vep/2$. Observe that $(z_n)_n$ is unconditional (Proposition \ref{sdfewooo1})
therefore equivalent  to the $\ell_1$-basis which yields a contradiction. \fprue

\cor For every set of ordinals $A$ we have that $\dim J_X^*(A)= \#
A$. \qed \fcor

\subsection{ Finite interval representability of $J_{T_0}$ and the space of diagonal operators}

\defi
Let $X$ and $Y$ be Banach spaces and let $(x_\al)_{\al<\ga}$ and $(y_n)_n$ be a transfinite basis
for $X$ and a Schauder basis of $Y$ respectively. We say that $Y$ is \emph{finitely interval
representable} in $X$ if there exists a constant $C>0$ such that for every integer $n$ and
intervals $I_1\le I_2\le \dots \le I_n$ successive, not necessarily distinct, intervals of $\ga$
there exists $z_i\in \langle (x_\al)_{\al\in I_i}\rangle$ ($i=1,\dots,n$) with $\supp z_1<\supp
z_2< \dots <\supp z_n$ and such that the natural order preserving isomorphism $H:\langle
(y_i)_{i=1}^n\rangle \to \langle (z_i)_{i=1}^n\rangle $ satisfies $\nrm{H}\cdot \nrm{H^{-1}}\le
C$.
\fdefi

Recall that Maurey-Rosenthal \cite{mau-ros}, in their attempt to solve the unconditional basic
 sequence problem, have constructed a Banach space $X$ with a weakly-null normalized Schauder
basis $(e_n)_n$ having the property that every subsequence of $(e_n)_n$ finitely block
represents the James-like space $J_{c_0}$, or equivalently (and as they said it), every
subsequence of $(e_n)_n$ has a arbitrary large finite block subsequence of length $k$
equivalent to the first $k$-many members of the summing basis of $c_0$. In our attempt to
control non-strictly singular operators on $\eqs_{\ou}$, we have discovered the following
analogous result that surprised us by its powers to explain many phenomena encountered, not
only in $\eqs_{\ou}$, but in essentially any other conditional space constructed so far using
the general scheme described above in Section 2. Through all this section $\ga$ will denote a
limit  ordinal.

\teor\label{represent} Let $(y_\al)_{\al<\ga}$ be a normalized transfinite block sequence
in $\eqs_{\ou}$, and $Y$ its closed linear span. Then $J_{T_0}$ is finitely interval
representable in the space $Y$, where $T_0$ is the mixed Tsirelson space $T[(
m_{2j}^{-1}, {n_{2j}})_j]$. \fteor We will postpone the proof until Section
\ref{finblorep}. Throughout all this section $\boldsymbol{C}$ will denote the finitely
block representability constant of $J_{T_0}$ in $\eqs_{\ou}$. We will show in Section
\ref{finblorep} that $\boldsymbol{C}< 121$. \nota \noindent 1. Let us observe that since,
as we will show, the basis of $J_{T_0}$ is not unconditional and it is finitely block
representable in any block subsequence of the basis $(e_\al)_{\al<\ou}$, then
$\eqs_{\ou}$ cannot have any unconditional basic sequence. In other words the finite
interval representability of $J_{T_0}$ in the block subsequences of $\eqs_{\ou}$ must
make use of the conditional structure of $\eqs_{\ou}$. Indeed we get more. Suppose that
$\eqs$ has a transfinite basis, and suppose that a Banach space $Y$ with a conditional
basis $(y_n)_n$ is finite block representable in every block sequence of $\eqs$ Then
$\eqs$ does not contain unconditional basic sequences and from Gowers dichotomy
\cite{G1}, $\eqs$ is hereditarily indecomposable saturated.

\noindent 2. The James like space $J_{T_0}$ has the following alternative description. It is
the mixed Tsirelson space $T_G[( m_{2j}^{-1},{n_{2j}})_j]$, where $G=\conj{I^*}{I\con \N
\text{ interval}}$. The minimal set $K_0$ of $c_{00}(\N)$ which is symmetric, contains $G$,
and is closed under $( m_{2j}^{-1},{n_{2j}})$-operations norms $J_{T_0}$. \fnota

\prop\label{lkio1}
Let $x_1<\dots < x_n$ be finitely supported,   $\phi\in K_{\ou}$ and set $r_i=\phi x_i$  for
each $i=1,\dots,n$. Then $\nrm{\sum_{i=1}^n r_i v_i}_{J_{T_0}}\le \nrm{x_1+\dots +x_n}$.
\fprop \prue Fix a functional $f$ of $K_0$ with support contained in $\{1,\dots,n\}$, and a tree-analysis $(f_t)_{t\in \mc T}$ of $f$. We show by induction over the tree $\mc T$ that for every $t\in \mc T$
there is some $\phi_t\in K_{\ou}$ such that $f_t (\sum_{i=1}^n r_i v_i )=\phi_t (x_1+\dots
+x_n)$. In particular $f_0 (\sum_i r_i v_i) =\phi_0(x_1+\dots +x_n)$, and hence the desired
result holds. If $t\in \mc T$ is a terminal node, then $f_t=\pm I^*$, $I\con \{1,\dots,n\}$ an
interval. We set $\phi_t=\pm \phi|[\min \supp x_{\min I},\max \supp x_{\max I} ]$. It is clear
that $\phi_t\in K_{\ou}$, and
\begin{equation}
\phi_t (x_1+\dots + x_n)= \pm \sum_{i\in I}\phi x_i =\pm \sum_{i\in I}r_i =f_t (\sum_i r_i v_i).
\end{equation}
If $t\in \mc T$ is not a terminal node, then $f_t=(1/m_{2j})\sum_{i=1}^d f_{s_i}$, where
$S_t=\{s_1,\dots,s_d\}$ ordered by $f_{s_1}<\dots <f_{s_d}$. Then $\phi_t=(1/m_{2j})\sum_{i=1}^d
\phi_{s_i}$ clearly satisfies our inductive requirements.
\fprue

The next result shows that $J_{T_0}$ is minimal in a precise sense.
\cor Suppose that $X$ is a
Banach space with a normalized Schauder basis $(x_n)_n$ which  dominates the summing basis of
$c_0$ and is finitely block represented in $\eqs_{\ou}$. Then $(x_n)_n$ also dominates the
basis $(v_n)_n$ of $J_{T_0}$. \fcor \prue Fix scalars $(a_i)_{i=1}^n$. Choose a normalized
block sequence $(w_i)_i^n$ of $\eqs_{\ou}$  $C$-equivalent
 to $(x_i)_{i=1}^n$. Fix $f\in K_0$ with $\supp f \con\{1,\dots,n\}$ and a tree-analysis  $(f_t)_{t\in\mc T}$
 of it. We are going to find  $\phi_t\in K_{\ou}$ such that
$|f_t(\sum_{i=1}^n a_i v_i)|\le C|\phi_t(\sum a_i w_i )| $, for each $t\in \mc T$. This will show
that
\begin{equation}
\nrm{\sum_{i=1}^n a_i v_i}_{J_{T_0}}\le C\nrm{\sum_{i=1}^na_i w_i}_{\eqs_{\ou}}\le C^2 \nrm{\sum_{i=1}^na_i
w_i}_{X},
\end{equation}
as desired. If $t\in \mc T $ is a terminal node, then $f_t=\pm I^*$, $I\con [1,n]$
interval. Since $(x_n)_n$ dominates the summing basis of $c_0$, we can find $\phi_t\in
K_{\ou}$ such that
\begin{equation}
\phi_t (\sum_{i=1}^n a_i w_i)=\nrm{\sum_{i=1}^n a_i w_i}_{\eqs_{\ou}}\ge \frac1C\nrm{\sum_{i=1}^n
a_i x_i}_X\ge \frac 1 C |\sum_{i\in I}a_i|=|f_t(\sum_{i=1}^n a_i v_i)|.
\end{equation}
If $t$ is not terminal node, then we use the appropriate $(m_{2j}^{-1},
{n_{2j}},)$-operation.\fprue

\defi
Let $(x_\al)_{\al<\ga}$ be a normalized transfinite block sequence, $X$ its closed linear span. We denote by
$\mc D(X)$ the space of all bounded diagonal operators  $D:X\to X$ satisfying the property that for all
$\al<\ga$ limit there exists some $\la_\al\in \R$ such that $D(x_\be)=\la_\al x_\be$ for every $\be\in
[\al,\al+\om)$. We  also denote by $\widetilde{\mc D}(X)$ the space of all diagonal operators (not
necessarily bounded) satisfying the above condition  acting on $\langle x_\al \rangle_{\al<\ga}$.

Notice the following (linear) decomposition of $\langle x_\al \rangle_{\al<\ga}$,
\begin{equation}  \label{swqoii2}
\langle x_\al \rangle_{\al<\ga}=\bigoplus_{\al\in \La(\ga)} \langle x_\be\rangle_{\be\in
[\al,\al+\om)}.\end{equation} The \emph{canonical decomposition} of  $y\in\langle x_\al
\rangle_{\al<\ga}$ in $X$ is $y=y_1+\dots +y_n$  given by (\ref{swqoii2}). \fdefi \nota $\mc D(X)$
is a closed subalgebra of $\mc L(X)$.

\fnota For an ordinal $\mu$ we denote by $\La(\mu)$ the set of limit ordinals $<\mu$, and
by $\La(\mu)^{(0)}$ the set of limit ordinals $\al=\be+\om<\mu$ with $\be\in \La(\mu)$.
We denote this (unique) $\be$ by $\al^-$. Notice that $\La(\mu)^{(0)}$ is the set of
isolated points of $\La(\mu)$   with respect to the order-topology. For technical
reasons, 0 is considered as  limit ordinal.

\nota Notice that  for $\ga$ a limit ordinal,  $\La(\ga+1)^{(0)}$ is order isomorphic to $\La(\ga)$
via the predecessor map.  \fnota
\defi
Let $D\in \widetilde{\mc D}(X)$. We define the map $\xi_D:\La(\ga+1)^{(0)}\to \R$ by
\begin{equation}
D(x_{\al^-})=\xi_D(\al)x_{\al^-}.
\end{equation} Namely, $\xi_D(\al)$ is the eigenvalue of $D$ associated
to the eigenvectors $(x_\be)_{\be\in [\al^-,\al)}$. \fdefi

We consider the following linear map $\Xi:\widetilde{D}(X)\to
c_{00}(\La(\ga+1)^{(0)})^\#$ defined by\begin{equation} \Xi(D)(v_\al)=\xi_D(\al),
\end{equation}
where $c_{00}(\La(\ga+1)^{(0)})^\#$ denotes the algebraic conjugate of $c_{00}(\La(\ga+1)^{(0)}$. The main
goal here is to show that $\Xi$ defines an isomorphism between $\mc D(X)$ and $J_{T_0}^*(\La(\ga+1)^{(0)})$.
For $D\in \widetilde{\mc D}(X)$, let us denote
\begin{align*}
\nrm{D}=\sup\conj{\nrm{Dx}_{\eqs_{\ou}}}{x\in \langle  x_{\al} \rangle_{\al<\ga}, \, \nrm{x}_{\eqs_{\ou}}\le
1}\le \infty,
\end{align*}
and for $f\in c_{00}(\La(\ga+1)^{(0)})^{\#}$,
\begin{align*}
\nrm{f}=\sup \conj{f (x)}{x\in c_{00}(\la(\ga+1)^{(0)}), \, \nrm{x}_{J_{T_0}}\le 1}\le \infty.
\end{align*}

\prop\label{dandtheta1} $\nrm{D}\le \nrm{\Xi(D)}\le \boldsymbol{C}\nrm{D}$ for every $D\in
\widetilde{\mc D}(X)$.
 \fprop \prue
Fix $D\in \widetilde{\mc D}(X)$, and $\vep>0$. Let $y\in \langle  x_{\al} \rangle_{\al<\ga}$ with $\nrm{y}\le
1$ be such that $|\nrm{D}-\nrm{Dy}|<\vep$.  Let $y=y_1+\dots +y_n$ be the canonical decomposition of $y$ in
$X$, and $\al_1,\dots,\al_n$ be such that $y_i\in \langle x_\be \rangle_{\be\in [\al_i^-,\al_i)}$ for every
$1\le i\le n$. Let $\phi\in K$ be such that $\nrm{Dy}=\phi(Dy)$, and  set $r_i=\phi y_i$ for $i=1,\dots,n$.
By Proposition \ref{lkio1}, $\nrm{\sum_{i=1}^n r_iv_i}_{J_{T_0}}\le \nrm{x}$, and since $(v_\al)_\al$ is
1-subsymmetric we have that $\nrm{\sum_{i=1}^n r_iv_{\al_i}}_{J_{T_0}}\le \nrm{y}\le 1$. Hence
\begin{equation}
\nrm{\Xi(D)}\ge \nrm{\Xi (D)(\sum_{i=1}^n r_iv_{\al_i})}_{J_{T_0}}=\nrm{\sum_{i=1}^n
\xi_D(\al_i)r_iv_{\al_i}}_{J_{T_0}}\ge \sum_{i=1}^n \xi_D(\al_i)=\phi(Dy)\ge \nrm{D}-\vep.
\end{equation}
This shows that $\nrm{D}\le \nrm{\Xi (D)}$. Fix $v=\sum_{i=1}^n a_i v_{\al_i}\in J_{T_0}$ with
$\nrm{v}_{J_{T_0}}\le 1$, and choose a finite normalized block sequence $(w_i)_{i=1}^n$
$\boldsymbol{C}$-equivalent to $(v_{\al_i})_{i=1}^n$ with $w_i\in \langle x_\be\rangle_{\be\in
[\al_i^{-},\al_i)}$ for every $i=1,\dots,n$ (indeed we may assume that the natural isomorphism $F:\langle
w_i\rangle_{i=1}^n\to \langle v_i\rangle_{i=1}^n$ satisfies that $\nrm{F}\le 1$, $\nrm{F^{-1}}\le
\boldsymbol{C}$; see Corollary \ref{mayassme}). Then,
\begin{align}
\nrm{\Xi(D)(v)}_{J_{T_0}}& =     \nrm{\sum_{i=1}^n \xi_D(\al_i)a_iv_{\al_i}}_{J_{T_0}}\le
\nrm{\sum_{i=1}^n
\xi_D(\al_i)a_iw_i}_{\eqs_{\ou}}=\nrm{D(\sum_{i=1}^n a_iw_i)}_{\eqs_{\ou}} \le  \nonumber\\
&   \le \nrm{D}\nrm{\sum_{i=1}^n a_iw_i}_{\eqs_{\ou}}\le \boldsymbol{C}\nrm{D}.
\end{align}\fprue

%

\teor\label{isosio} The spaces $\mc D (X)$ and $J_{T_0}^*(\La(\ga+1)^{(0)})$ are isomorphic.
\fteor
\prue By  Proposition \ref{dandtheta1}, $\Xi|\mc D(X):\mc D (X)\to J_{T_0}^*(\La(\ga+1)^{(0)})$ is an isomorphism. To see that
 it is also onto consider
$f\in J_{T_0}^*(\La(\ga+1)^{(0)})$ and define $D_f\in \widetilde{\mc D}(X)$ as follows. For $\be\in
[\al^-,\al)$ set $D_f(x_\be)=f(v_\al)x_\be$. It is easy to check that $\Xi(D_f)=f$. This completes the proof.
\fprue

\cor\label{samediagonals} Let $X$ and $Y$ be the closed linear
span of two  transfinite block sequences of the same length $\ga$.
Then the natural mapping $\psi_\ga:\mc D(X)\to \mc D(Y)$ defined
by $\psi_\ga(D)=D_{\xi_D}$ is an isomorphism. \qed \fcor

Our intention now is to compare $\mc D(X)$ and $\mc D (\eqs_{\ou})$.
\defi\label{progammax}
\noindent 1. Given a closed $A\con \La(\ou+1)$, let $\widetilde{\mc D}_{A}(\eqs_{\ou})$
be the subalgebra of $\widetilde{\mc D}(\eqs_{\ou})$ consisting on all $D\in \mc
D(\eqs_{\ou})$ satisfying that for every $\al\in A^{(0)}$, there is some $\la_\al$ such
that $D |\eqs_{[\al^-,\al)} =\la_\al i_{\eqs_{[\al^-,\al)},\eqs_{\ou}}$ and $D
|\eqs_{[\max A,\ou)}=0$. Let $\mc D_A(\eqs_{\ou})$ be the subalgebra of bounded operators
of $\widetilde{\mc D}_{A}(\eqs_{\ou})$.

\noindent 2. Given a transfinite block sequence $(x_\al)_{\al<\ga}$, let $\Ga_{X} \con \La(\ou+1)$ be defined
as follows. Let
\begin{equation}
\Ga'=\conj{\sup_{n \to \infty }  \max \supp x_{\al_n}}{(\al_n)_n\uparrow, \al_n<\ga },
\end{equation}
and let $\Ga_{X}=\Ga'\cup \{0,\sup \Ga'\}$.  Another interpretation of $\Ga_X$ is to consider the
map $f_X:\La(\ga+1)\to \ou$ defined by $f_X(\al)=\sup_{\be<\al}\max \supp x_\be$ and $\Ga_X$ is
nothing else but the  image $f (\La(\ga+1))$, and hence $\Ga_X\setminus \max\{ \Ga_X\}$ and
$\La(\ga+1)^{(0)}$ are order isomorphic.

\indent 3. Given $D\in \mc D(X)$, let $E(D)\in \widetilde{\mc D}_{\Ga_X}( \eqs_{\ou})$ be the
unique extension of $D$. Notice that $D|X\in \mc D(X)$ for every $D\in {\mc D}_{\Ga_X}(
\eqs_{\ou})$. \fdefi

\teor[\emph{Extension Theorem}]\label{extensiondiagonal} For every $X\hookrightarrow \eqs_{\ou}$
generated by a transfinite block sequence the following hold:

\noindent (a) Every $D\in \mc D(X)$ is extended to a step diagonal operator $E D$ in $\mc
D(\eqs_{\ou})$.

\noindent (b) The map $D\mapsto ED$ defines a linear isomorphism from $\mc D(X)$ onto the
space $\mc D_{\Ga_X}(\eqs_{\ou})$.
%
\fteor

\prue We show that  $\nrm{E(D)}\le \boldsymbol{C}\nrm{D}$ for every $D\in \mc D(X)$.  Fix  a finitely supported $y\in \eqs_{\ou}$  such that $\nrm{y}\le 1$ and
$\nrm{E(D)}=\nrm{E(D)(y)}$. Since $\mc I= \conj{[\al^-_{\Ga_X},\al)}{\al\in \Ga_X^{(0)}}\cup \{[\max
\Ga_X,\ou)\}$ is a partition of $\ou$, $y$ has a unique decomposition $y=y_1+\dots +y_n$ for  $I_1<\dots
<I_n$ in $\mc I$ and $y_i\in \langle e_{\al} \rangle_{\al\in I_i}$. Notice that $E(D)|\eqs_{[\max
\Ga_X,\ou)}=0$, so we may assume that $I_n\neq [\max \Ga_X,\ou)$. By definition of $E(D)$ we have that
$E(D)(y)=\sum_{i=1}^n \xi_D(\be_i) y_{i}$ where $\be_i=f_X^{-1}(\al_i)$ for every $i=1,\dots,n$. Choose
$\phi\in K_{\ou}$ such that $\nrm{E(D)(y)}=  \phi (E(D)(y))$. By Proposition \ref{dandtheta1},
\begin{align}
\nrm{E(D)}= &    \phi (\sum_{i=1}^n \xi_D(\be_i) y_{i})= \sum_{i=1}^n \xi_D(\be_i)\phi
(y_i)=\Xi(D)(\sum_{i=1}^n \phi(y_i) v_{\be_i})\le \nonumber \\
\le &    \nrm{\Xi(D)}_{J^*_{T_0}(\La(\ga+1)^{(0)})}\nrm{\sum_{i=1}^n \phi(y_i)
v_i}_{J_{T_0}}\le \boldsymbol{C} \nrm{D}.
\end{align}
\fprue

\subsection{The spaces $\mc L(\eqs_{\ga})$}

\defi
A sequence $(x_1,\phi_1,\dots,x_{n_{2j+1}},\phi_{n_{2j+1}})$ is called a
$(0,j)$\emph{-dependent sequence } if the following conditions are fulfilled:
\begin{enumerate}
\item[\emph{DS0.1}] $\Phi=(\phi_1,\dots,\phi_{n_{2j+1}})$ is a
$2j+1$-special sequence and $\phi_i x_{i'}=0$ for every $1\le
i,i'\le n_{2j+1}$. \item[\emph{DS0.2}] There exists
$\{\psi_1,\dots,\psi_{n_{2j+1}}\}$ such that
$w(\psi_i)=w(\phi_i)$, $\#\supp x_i \le w(\phi_{i+1})/n_{2j+1}^2$
and $(x_i,\psi_i)$ is a $(6,2j_i)$-exact pair  for every $1\le
i\le n_{2j+1}$. \item[\emph{DS0.3}] If
$H=(h_1,\dots,h_{n_{2j+1}})$ is  an arbitrary $2j+1$-special
sequence, then
\begin{equation}
\left(\bigcup_{\ka_{\Phi,H}<i<\la_{\Phi,H}}\supp x_i\right) \cap
\left(\bigcup_{\ka_{\Phi,H}<i<\la_{\Phi,H}}\supp h_i
\right)=\emptyset.
\end{equation}
\end{enumerate}
\fdefi \prop  For every $(0,j)$-dependent sequence
$(x_1,\phi_1,\dots,x_{n_{2j+1}},\phi_{n_{2j+1}})$ we have that
\[\|  \frac{1}{n_{2j+1}}(x_1+\dots +x_{n_{2j+1}})\|\le \frac{1}{{m_{2j+1}^2}}.\] \fprop
\prue The proof is rather similar to the proof of Proposition \ref{depseqest}. One first
shows that $|\psi(1/{n_{2j+1}}\sum_{i\in E}   x_i)|\le 12(1+{\#E}/{n_{2j+1}^2})$ for
every special functional $\psi$ with $w(\psi)=m_{2j+1}$, and then the result follows from
the basic inequality, since, by condition (\emph{DS0.2}), $(x_i)_{i=1}^{n_{2j+1}}$ is a
$(12,1/n_{2j+1}^2)$-RIS. \fprue

\prop \label{keyforbridges}
 Suppose that $(y_k)_k$ is a $(C,\vep)$-RIS, and suppose that  $T:\langle y_k \rangle_k\to \eqs_{\ou}$ is
 a linear function
(not necessarily bounded)  such that $\lim_{n\to \infty} d(T y_n,\R y_n)\neq 0$. Then for
every $\vep>0$ there is some $z\in \langle y_k \rangle_k$ such that $\nrm{z}<\vep \nrm{T z}$.
\fprop \prue We may assume that there is some  $\de>0$ such that $\inf_n d(Ty_{n},\R
y_{n})>\de>0$, and also that $(T y_n)_n$ is a block sequence (hint: Consider the following
limit ordinal
\begin{equation}
\ga_0=\min\{\ga<\ou\, :\,\exists A\in [\N]^\infty \, \inf_{n\in A}d(\P_\ga Ty_n,\R y_{n})>0\},
\end{equation}
pass  to a subsequence of $(y_n)_n$ and replace $T$ by $\P_{\ga_0} T$).

\clam There exist an infinite set $A\con \N$  and a block sequence $(f_n)_{n\in A}$ of
functionals in $K_{\ou}$ such that:

\noindent \emph{(a)} For every $n\in A$, $f_n Ty_n\ge \de,\,f_n y_n=0, \,  \ran f_n\con \ran
Ty_n\text{ and } \supp f_n \cap \supp y_m=\buit \text{ for every $m\neq n$}$.

\noindent \emph{(b)} Either for every $n\in A$ $\max \supp y_n\ge \max \supp f_n$ or for every
$n\in A$ $\max \supp y_n\le \max \supp f_n$.

\fclam \prucl By the Hahn-Banach theorem, for each $n\in \N$ we can find a functional $f_n$
of norm 1 such that $f_n(T y_n)\ge \de$ and $f_n(y_n)=0$. Since the $w^*$-closure of
$K_{\ou}$ is $B_{\eqs_{\ou}^*}$ (notice that $K$ by definition is closed under rational
convex combinations) and   $K_{\ou}$ is closed under restriction over intervals, we may
assume that $f_n\in K_{\ou}$ and $\ran f_n\con \ran Ty_n$. Let $\al=\max_n\supp y_n$ and
$\be=\max_n\supp f_n$. If $\al\neq \be$, it is rather easy to achieve the desired result.
If $\al=\be$, then we can pass to a subsequence $A$ and distort $f_n$ such that for every
$n\in A$, $\max \supp f_n\ge \max \supp y_n$. \fprucl

So, we may assume that $(f_n)_{n}$ satisfies the requirements of previous Claim. Fix $j$
 with $m_{2j+1}>12/(\vep\de)$.

 \clam There is a $(0,j)$-dependent sequence $(z_1,\phi_1,\dots,z_{n_{2j+1}},\phi_{n_{2j+1}})$
such that  for every $k\le n_{2j+1}$, $z_k\in X$, $\ran \phi_k\con \ran T z_k$ and $\phi_k T
z_k>\de $.
\fclam
\prucl Choose $j_1$ even such that $m_{2j_1}>n_{2j+1}^2$, and choose $F_1\con \N $ of size
$n_{2j_1}$ such that $(y_k)_{k\in F_1}$ is a $(3,1/n_{2j_1}^2)$-RIS (going to a subsequence of
$(y_k)_k$; see Remark \ref{RRIS}). Set
\begin{align*}
 \phi_1= &
\frac1{m_{2j_1}}\sum_{i\in F_1}f_i\in K_{\ou} \text{ and }   z_1=
\frac{m_{2j_1}}{n_{2j_1}}\sum_{k\in F_1}y_k.
\end{align*}
Notice that $\phi_1 Tz_1=(1/{n_{2j_1}})\sum_{k\in F_1}f_k T y_k >\de$ and by  (a) from the
Claim, we have that $\phi_1 z_1=  (1/{n_{2j_1}})  \sum_{k\in F_1}\sum_{l\in F_1}f_k (y_l)=0$.
Pick
\begin{equation}
p_1\ge \max\{ p_\varrho (\supp z_1 \cup \supp Tz_1 \cup \supp \phi_1), \# \supp z_1 {\cdot}
n_{2j+1}^2\}
\end{equation}
and set $ 2j_2= \sig_\ro(\Phi_1, m_{2j_1}, p_1)$. Now choose $F_2>F_1$ finite of length
$n_{2j_2}$ such that $(x_k)_{k\in F_2}$ is a $(3,1/n_{2j_2}^2)$-RIS. Set
\begin{equation}
\phi_2 = \frac{1}{m_{2j_2}}\sum_{k\in F_2}f_k\in K_{\ou}\text{ and }
z_2=\frac{m_{2j_2}}{n_{2j_2}}\sum_{k\in F_2}y_k.
\end{equation}
Notice that $\phi_2>\phi_1$, $\phi_2 Tz_2>\de$ and   $\phi_2 z_2=0$. Pick
\begin{equation}
p_2\ge\max\{ p_1, p_\varrho(\supp z_1\cup \supp z_2 \cup \supp T z_1\cup \supp T z_2 \cup\supp
\Phi_1 \cup \supp \Phi_2),\# \supp z_2{\cdot} n_{2j+1}^2\}\}
\end{equation} and set
$2j_3=\sig_\ro(\phi_1,m_{2j_1},p_1,\phi_2,m_{2j_2},p_2)$, and so on. Let us check that
$(z_1,\phi_1,\dots,z_{n_{2j+1}},\phi_{n_{2j+1}})$ is a $(0,j)$-dependent sequence: Condition
(\emph{DS0.1}) and (\emph{DS0.2}) are rather easy to check from the definition of this
sequence. Let us check (\emph{DS0.3}). There are two cases: (a) Suppose that $\max \supp
z_k\le\max \supp \phi_k$  for every $1\le k\le n_{2j+1}$. Then $\supp z_k \con \supp
\cl{\phi_{\la_{\Phi,H}-1}}{p_{\la_{\Phi,H}-1}}$ for every $\ka_{\Phi,H}<k< \la_{\Phi,H}$. Then
part 2 of (TP.3) gives the desired result. (b) Suppose that $ \max \supp \phi_k\le \max \supp
z_k$ for every $1\le k\le n_{2j+1}$. Then $\supp \phi_k \con \supp
\cl{z_{\la_{\Phi,H}-1}}{p_{\la_{\Phi,H}-1}}$ for every $\ka_{\Phi,H}<k< \la_{\Phi,H}$, and we
are done by part 1 of (TP.3).
\fprucl
Fix a $(0,j)$-dependent sequence $(z_1,\phi_1,\dots,z_n,\phi_{n_{2j+1}})$ as in the Claim, and
set
\begin{align*}
z=\frac{1}{n_{2j+1}}\sum_{k=1}^{n_{2j+1}}(-1)^{k+1}z_k \text{ and
}\phi=\frac{1}{m_{2j+1}}\sum_{k=1}^{n_{2j+1}}\phi_k.
\end{align*}
Then $\phi T z=1/{n_{2j+1}}\sum_{k=1}^{n_{2j+1}}(-1)^{k+1}\phi Tz_k\ge \de/{m_{2j+1}}$
and $\nrm{z}\le 12/m_{2j+1}^2$. So, $\nrm{T(z)}\ge {\de}/{m_{2j+1}}\ge {\de
m_{2j+1}\nrm{z}}/{12}>\vep \nrm{z}$ as desired. \fprue

\cor \label{mat1} Let $(y_k)_k$ be a $(C,\vep)$-RIS, $Y$ its
closed linear span and $T:Y\to \eqs_{\ou}$ be a bounded operator.
Then $\lim_{n\to \infty}d(Ty_k,\R y_k)=0$.
\fcor
\prue If not, by the previous Proposition \ref{keyforbridges}, we can find a vector $z\in
\langle y_k\rangle_k$ such that $\nrm{z}<(1/\nrm{T})\nrm{T z}$ which is impossible if $T$ is
bounded. \fprue

\lema\label{mat2} Let $(x_n)_n$ be a $(C,\vep)$-RIS, $X$ its
closed span and  $T:X\to \eqs_{\ou}$ be a bounded operator. Then
$\la_T:\N\to \R$ defined by $d(T x_n,\R
x_n)=\nrm{Tx_n-\la_T(n)x_n}$  is a convergent sequence. \flema
\prue Fix any two strictly increasing sequences $(\al_n)_n$ and
$(\be_n)_n$ with $\sup_n \al_n=\sup_n \be_n$, and suppose that
$\la_T(\al_n)\to_n \la_1$, $\la_T(\be_n)\to_n \la_2$. By going to
a subsequences, we can assume that $x_{\al_n}<x_{\be_n}$ for every
$n$. Since the closed linear span of $\{x_{\al_n}\}_n\cup
\{x_{\be_n}\}_n$ is an H.I. space, we can find for every $\vep$
two normalized vectors $w_1\in \langle x_{\al_n}\rangle_n$ and
$w_2\in \langle x_{\be_n}\rangle_n$ such that $\nrm{T
w_1-\la_1w_1}\le \vep/3$, $\nrm{Tw_2-\la_2w_2}\le \vep/3$ and
$\nrm{w_1-w_2}\le \vep/3\nrm{T}$. Then we have that
\begin{equation}
\nrm{\la_1w_1-\la_2w_2}\le \nrm{T w_1-\la_1
w_1}+\nrm{Tw_1-Tw_2}+\nrm{Tw_2-\la_2 w_2}\le \vep,
\end{equation}
and hence,
\begin{equation}
\vep\ge \nrm{\la_1w_1-\la_2w_2}\ge
|\la_1-\la_2|\nrm{w_1}-|\la_2|\nrm{w_1-w_2}\ge
|\la_1-\la_2|-|\la_2|\vep.
\end{equation}
So, $|\la_1-\la_2|\le \vep (1+|\la_2|)$ for every $\vep$. This implies that $\la_1=\la_2$.
\fprue

\defi
 Recall that for a set $A$ of ordinals $A^{(0)}$ is the set of isolated points of $A$.
 Fix a transfinite block sequence $(x_\al)_{\al<\ga}$,  let $X$ be the closed linear span of it and  let $T:X\to
\eqs_{\ou}$ be a bounded operator. We define the \emph{step function } ${\xi_T}$ of $T$
$\xi_T:\La(\ga+1)^{(0)}  \to \R$ as follows: Let $\ga$ be a successor limit ordinal less than $\ga$. Let
$\xi_T(\ga)=\xi\in \R$ be such that $\lim_{n\to \infty}\nrm{Ty_n -\xi y_n}=0$ for every $(3,\vep)$-RIS
$(y_n)_n$  satisfying that $\sup_n\max\supp y_n=\ga$. Lemma \ref{mat2} shows that $\xi$ exists and  is
unique, and that $\xi_T$ can be extended to a continuous $\widetilde{\xi}_T:\La(\ga+1)\to \R$.

Given a mapping $\xi:\La(\ga+1)^{(0)}\to \R$ we define the diagonal, not necessarily
bounded, operator $D_\xi:X\to X$ in the natural way by $D_\xi(x_\al)=\xi(\al+\om)x_\al$.
Given a bounded $T:X\to \eqs_{\ou}$ we define the  \emph{diagonal step operator}
$D_T:\langle x_\al \rangle_{\al<\ga}\to {\eqs_{\ou}}$  of $T$  as $D_T=D_{\xi_T}$. \fdefi

\nota The function $\xi_T$ has only countable many values. This follows from the fact
that it  can be extended to a continuous function  $\widetilde{\xi}_T$ defined on
$\La(\ga+1)$. As it is well known, if $\ga=\ou$ the function $\widetilde{\xi}_T$ is
eventually constant. \fnota

\prop\label{diss0} The sequence $(\nrm{(T-D_T)(y_n)})_{n}\in c_0(\N)$ for every RIS
$(y_n)_n$ in $X$. \fprop \prue This is just a consequence of the definition of $D_T$.
\fprue

\prop\label{sscarac} A bounded  operator $T:X\to \eqs_{\ou}$     is strictly singular iff
$\xi_T=0$.

 \fprop \prue Suppose
that $T$ is not strictly singular. Then there is a block sequence $(y_n)_n$ such that $T$
is an isomorphism restricted to the closed linear span $Y$ of $(y_n)_n$. Going to a block
subsequence if necessary, we assume that $(y_n)_n$ is a RIS. Since $T|Y$ is an
isomorphism, $\lim_{n\to \infty} \nrm{Ty_n}>0$. This implies that
$\xi_T|\La(\al+1)^{(0)}\neq 0$, since otherwise $\widetilde{\xi}_T(\al)=0$ contradicting
the above inequality.

Suppose now that $\xi_T\neq 0$. Choose some successor limit $\ga$ such that $\xi_T(\ga)\neq 0$.
Then we can find a block sequence $(y_n)_n\con X_\ga$ such that $T$ is close enough to $\xi_T(\ga)
i_{Y,\eqs_{\ou}}$, where $Y$ is the closed linear span of $(y_n)_n$. Hence, $T$ is not strictly
singular. \fprue

\prop\label{disbounded} Let $(x_\al)_{\al<\ga}$ be a transfinite block sequence, $X$ its closed
linear span of $(x_\al)_{\al<\ga}$ and a bounded operator $T:X\to \eqs_{\ou}$. Then $\nrm{D_T}\le
\boldsymbol{C}\nrm{T}$ and hence $D_T\in \mc D(X)$.
\fprop
\prue Fix a normalized   $y\in \langle x_\al\rangle_{\al<\ga}$. Let  $y=y_1+\dots +y_n$ be its
decomposition  in $X$, $y_i\in \langle x_\be\rangle_{\be\in [\al_i^{-},\al_i)}$ for
$i=1,\dots,n$. Choose $\phi\in K_{\ou}$ such that $\phi(D(y))=\nrm{D(y)}$. Then,
\begin{equation}  \label{jhyu2}
\nrm{D(y)}=\sum_{i=1}^n  \xi_T(\al_i)\phi (y_i)=(\sum_{i=1}^n \xi_T(\al_i)
v_i^*)(\sum_{i=1}^n \phi(y_i)v_i)\le \nrm{\sum_{i=1}^n \xi_T(\al_i)  v_i}_{J_{T_0}},
\end{equation}
the last inequality holding because $\nrm{\sum_{i=1}^n \phi(y_i)v_{\al_i}}_{J_{T_0}}\le
\nrm{y}_{\eqs_{\ou}}\le 1$.
 We finish with the next claim. \clam $\nrm{\sum_{i=1}^n \xi_T(\al_i)  v_i^*}_{J_{T_0}^*}\le
\boldsymbol{C}\nrm{T}$. \fclam \prucl Fix $\vep>0$. By the finitely block representability of
$J_{T_0}$ in $\eqs_{\ou}$ and Proposition \ref{diss0} we can produce inductively $w_1,\dots,w_n $
such that \emph{(1)} $w_i\in \langle x_{\be}\rangle_{\be\in [\al_i^-,\al_i)}$,

\noindent \emph{(2)} the natural isomorphism
 $F:\langle  w_i \rangle_{i=1}^n\to \langle v_i
\rangle_{i=1}^n$ is such that $\nrm{F}\le 1$ and $\nrm{F^{-1}}\le \boldsymbol{C}$, and

\noindent \emph{(3)} $\sum_{i=1}^n \nrm{\xi_T(\al_i) w_i -T w_i}<\vep$.

Choose $x=\sum_{i=1}^n r_i v_i\in J_{T_0}$ of norm 1 such that $\nrm{\sum_{i=1}^n \xi_T(\al_i)
v_i^*}_{J_{T_0}^*}=\sum_{i=1}^n \xi_T(\al_i) r_i$. Then $\nrm{\sum_{i=1}^n r_i
w_i}_{\eqs_{\ou}}\le \boldsymbol{C}$ and hence
\begin{equation}
\nrm{D_T(\sum_{i=1}^n r_i w_i)}\ge  \nrm{\sum_{i=1}^n r_i \xi_T(\al_i) v_i}_{J_{T_0}}\ge \sum_{i=1}^n
\xi_T(\al_i) r_i= \nrm{\sum_{i=1}^n \xi_T(\al_i) v_i^*}_{J_{T_0}^*}.
\end{equation}
This implies that $\nrm{\sum_{i=1}^n \xi_T(\al_i) v_i^*}_{J_{T_0}^*}\le
\nrm{T(\sum_{i=1}^n r_i w_i)}+\nrm{(T-D_T)(\sum_{i=1}^n r_i w_i)}\le
\boldsymbol{C}\nrm{T}+ \vep$. \fprucl \fprue \teor\label{operblock} Let
$(x_\al)_{\al<\ga}$ be a normalized block sequence of $\eqs_{\ou}$, $X $ its closed
linear span. Then for every  bounded operator $T:X \to \eqs_{\ou}$, $D_T:X \to
\eqs_{\ou}$ is bounded and $T-D_T$ is strictly singular. \fteor

\prue This follows from Proposition \ref{sscarac} and Proposition \ref{disbounded}.
\fprue
\cor\label{decblock} Any bounded operator from the closed linear span $X$ of a
transfinite block sequence into the space $\eqs_{\ou}$ is the sum of the restriction of a
unique diagonal step operator $D\in \mc D_X(\eqs_{\ou})$  and an strictly singular
operator.
\fcor
\prue This follows from the previous theorem and Theorem \ref{extensiondiagonal}. \fprue

\cor \emph{(1)} For $T:X\to \eqs_{\ou}$ bounded TFAE: (a) $T$ is strictly singular, (b)
$\xi_T=0$, and (c) $D_T=0$.

\noindent \emph{(2)} The transformation $T\mapsto D_T$ is a projection in the operator
algebra $\mathcal{L}(X)$ of norm  $\le \boldsymbol{C}$. \qed \fcor

\prop\label{blockseverywhere} Let $X\hookrightarrow \eqs_{\ou}$, $I\con \ou$ an interval such that
$P_I|X$ is not strictly singular. Then for every $\vep>0$ there exist a normalized sequence
$(x_n)_n$ in $X$ and a normalized block sequence $(z_n)_n$ in $\eqs_{I}$ such that $\sum_n
\nrm{y_z-z_n}<\vep$. \fprop \prue Set $I=[\al,\be]$ and suppose that $P_I|X$ is not strictly
singular. Let
\begin{align*}
\ga_0=\conj{\ga\in (\al,\be]}{P_\ga |X\text{ is not strictly singular}}.
\end{align*}
We can find for every $\vep>0$, $(y_n)_n\con X$ and a block sequence $(w_n)_n\con \eqs_{\ga_0}$ such that
$P_{\ga_0}$ is an isomorphism when restricted to the closed linear span of $(y_n)_n$, $\sup_n\max\supp
w_n=\ga_0$ and $\sum_n\nrm{w_n-P_{\ga_0}y_n}\le \vep/2$. Consider $U:\overline{\langle w_n \rangle_n}\to
\eqs_{[\ga_0,\ou)}$ defined by $Uw_n=P_{[\ga_0,\ou)}y_n$. Notice that $U$ is bounded. Since $\xi_U=0$, $U$ is
strictly singular. Hence we can find a block sequence  $(z_n)_n$ of $(w_n)_n$ such that for all $n$, $\nrm{U
z_n}\le \vep/2^{n+1}$ and hence the corresponding block sequence $(x_n)_n$ of $(y_n)_n$ satisfies that
$\sum_n\nrm{z_n-x_n}\le \vep$. Finally, notice that for large enough $n_0$, $(z_{n})_{n\ge n_0}\con
\eqs_{I}$.
\fprue

\cor The space $\eqs_{\ou}$ is arbitrarily distortable. \fcor \prue For $j\in \N$, and
$x\in \eqs_{\ou}$, let $\nrm{x}_{2j}=\sup \conj{\phi(x)}{w(\phi)=m_{2j}}$. Let
$X\hookrightarrow \eqs_{\ou}$. Since for every $\vep>0$ we can find a subspace of $X$
generated by a Schauder basis $(y_n)_n$ and a normalized block sequence $(z_n)_n$ of
$\eqs_{\ou} $ such that $\sum_n \nrm{y_n-z_n}\le \vep$, without loss of generality we can
assume that $X$ is generated by a block sequence $(z_n)_n$. Now, we can find an
$(6,j)$-exact pair $(x,\phi)$, with $x\in \langle z_n\rangle_n$ and hence $1\le
\nrm{x}_{2j}\le \nrm{x}\le 6$. And for any other $j'>j$, a $(6,2j')$-exact pair
$(x',\phi')$ with $x'\in \langle z_n\rangle_n$ and hence $1\le \nrm{x'}\le 6$ and
$\nrm{x'}_{2j}\le 12/m_{2j} $. So,
\begin{equation}
\frac{ \nrm{ x/\nrm{x}}_{2j}}{ \nrm{ x'/\nrm{x'}}_{2j}}\ge \frac{1/6}{12/m_{2j+1}}=\frac{m_{2j+1}}{72}.
\end{equation}
\fprue

\defi
Two Banach spaces $X$ and $Y$ are called totally incomparable if and only if no infinite
dimensional closed $X_1\hookrightarrow X$ is isomorphic to $Y_1 \hookrightarrow Y$.
\fdefi

\cor\label{nojump} For disjoint infinite intervals $I$ and $J$, the spaces  $\eqs_I$ and
$\eqs_J$ are totally incomparable.\fcor \prue Suppose not, and let $X\hookrightarrow \eqs_I$,
and $Y \hookrightarrow \eqs_J$ such that $T:X\to Y$ is an onto isomorphism. By the previous
Proposition \ref{blockseverywhere}, we can assume that $X$ is generated by a block sequence.
But since $\xi_T=0$, $T$ cannot be isomorphism. This is a contradiction. \fprue

Another consequence of the representability of $J_{T_0}$ on each transfinite block sequence is
that we can identify  the space $\mathcal{D}(X)$ of diagonal step operators on $X$ and hence
identify $\mc L(X)/\mc S(X)$ for every closed span $X$ of a transfinite block sequence.

%
\cor\label{refofma} $\mc L(X)/\mc S(X) \cong \mc L(X,\eqs_{\ou})/\mc S(X,\eqs_{\ou}) \cong
J_{T_0}^*(\Ga_X^{(0)})$ for every $X\hookrightarrow \eqs_{\ou}$ generated by a transfinite block
sequence.
 \fcor
\prue This follows from Lemma \ref{isosio}, since
$\La(\ga+1)^{(0)}$ and $\Ga_X^{(0)}$ are order-isomorphic. \fprue

\nota\label{dfdlwerijghgg} Note that  $\mc L(X)/\mc S(X)   \cong
J_{T_0}^*(\Ga_X)$ if $\Ga_X$ is infinite. To see this, fix a
transfinite block sequence $(x_\al)_{\al<\ga}$ generating $X$ such
that $\ga\ge \om^2$. Then $\Ga_X\setminus \{\max \Ga_X\}$ and
$\La(\ga+1)^{(0)}\setminus \{\om\}$ are order-isomorphic.
 \fnota

\teor\label{fewproj} Every projection $P$ of $\eqs_{\ou}$ is of the form  $P=\P_{I_1}+\dots +
\P_{I_n}+S$, where $I_i$ are intervals of ordinals, $I_i<I_{i+1}$ and $S$ is strictly singular.
\fteor

\prue Suppose that $P:\eqs_{\ou}\to \eqs_{\ou}$ is a projection, $P=D_P+S$. Since
$P^2=P$, we obtain that $D_P^2-D_P$ is also strictly singular and therefore
$(\xi_P(\al)^2-\xi_P(\al))i_{\eqs_{[\al^-,\al)},\eqs_{\ou}}$ is strictly singular for
every successor limit $\al$. This implies that $\xi_P:\La(\ou+1)^{(0)}\to \{0,1\}$. And
since $\xi_P$ has the continuous extension property, there is no strictly increasing
sequence $\{\al_n\}_n\con \La(\ou+1)^{(0)}$ such that $\xi_P(\al_{2n})=1$ and
$\xi_P(\al_{2n+1})=0$ for every $n$. \fprue

\cor  For every
$n\in \N$ there is some $m\in \N$ such that for every projection $P$ of $\eqs_{\ou}$ with
$\nrm{P}\le n$, $P$ can be written as $P=P_{I_1}+\dots + P_{I_k}+S$ such that $k\le m$
and $I_1<<I_2<<\dots << I_k$, where    $A<<B$ denotes that the interval $(\sup A,\inf B)$
is infinite.
 \fcor
\prue  Fix $n$, and let $P:\eqs_{\ou}\to \eqs_{\ou}$ be a projection such that
$\nrm{P}\le n$. Let $j$ be the first integer such that $m_{2j}>2n\boldsymbol{C}$. We
claim that $m=n_{2j}$ works. For suppose that $P=P_{I_1}+\dots +P_{I_k}+S$ with
$I_1<<\dots <<I_k$ and $k> n_{2j}$. Fix $\vep>0$. Find a normalized block sequence
$(x_1,y_1,\dots,x_{n_{2j}/2},y_{n_{2j}/2})$ such that

\noindent (a) $x_i\in \eqs_{I_i}$, $y_i\in \eqs_{(\sup I_i,\min I_{i+1})}$ for $1\le i\le
n_{2j}/2-1$, and $y_{n_{2j}/2}> x_{n_{2j}/2}$,

\noindent (b) $(x_1,y_1,\dots,x_{n_{2j}/2},y_{n_{2j}/2})$ is $\boldsymbol{C}$-equivalent to
$(v_i)_{i=1}^{n_{2j}}$ and

\noindent (c) $\nrm{S|F}\le \vep$ where $F=\langle (x_1,,y_1,\dots,x_{n_{2j}/2},
y_{n_{2j}/2})\rangle$.

Set $x=x_1-y_1+\dots +x_{n_{2j}/2}-y_{n_{2j}/2}$. Then,
\begin{equation}  \label{equdss1}
\nrm{x}\le \boldsymbol{C}\nrm{\sum_{i=1}^{n_{2j}}(-1)^{i+1}v_i}_{J_{T_0}}\le
\boldsymbol{C}\nrm{\sum_{i=1}^{n_{2j}} t_i}_{T_0} =\boldsymbol{C} n_{2j}/m_{2j},
\end{equation}
and
\begin{equation}  \label{equdss2}
\nrm{P(x)}\ge \nrm{\sum_{i=1}^{n_{2j}/2}x_i}-\vep \ge \nrm{\sum_{i=1}^{n_{2j}/2} v_i}_{J_{T_0}}-\vep
=n_{2j}/2-\vep.
\end{equation}
 (\ref{equdss1}) and (\ref{equdss2}) imply that $\nrm{P}\ge (m_{2j}/2-\vep
m_{2j}/n_{2j} )/\boldsymbol{C}$. Hence, $\nrm{P}>n$, a contradiction.
\fprue

\subsection{Asymptotically equivalent subspaces and $\mc L(X,\eqs_{\ou})$}

Our aim here is to extend the results about operators on subspaces generated by a
transfinite block sequence to arbitrary subspace.

\defi\label{lamb}
Let $X$ be a subspace of $\eqs_{\ou}$. A subset $\Ga$ of $\ou+1$ is said to be a \emph{critical set } of $X$
if the following hold:
\begin{enumerate}
\item[\emph{(CS1)}] $\Ga$ is closed of limit ordinals, and $0\in \Ga$.
\item[\emph{(CS2)}] For all $\ga\in \Ga $, $\ga<\Om$,  $\P_{(\ga,\ga^+)}| X$ is not strictly
 singular and for all $\al \in (\ga,\ga^+)$, $\P_{(\ga,\al)}|X$ is strictly singular,
 where $\ga^+$ is the successor of $\ga$ in $\Ga$ and $\Om=\max \Ga$.
\item[\emph{(CS3)}] $\P_{[\Om,\ou)}|X$ is strictly singular (we use $P_{\buit}=0$).
\end{enumerate}
Notice that from the definition it follows easily that if $\Ga$ is a critical set of $X$, then
$\max \Ga=\min\conj{\ga\le \ou}{P_{[\ga,\ou)}|X \text{ is strictly singular}}$.
\fdefi
\prop For every
$X\hookrightarrow \eqs_{\ou}$ a critical set $\Ga$ is uniquely defined, denoted by $\Ga_X$.
\fprop \prue Fix $X\hookrightarrow \eqs_{\ou}$. We show first that a critical set $X$ exists.
We proceed by induction defining an increasing sequence $(\ga_\al)_{\al<\ou}$ as follows: We
set $\ga_0=0$. Suppose we have defined $(\ga_\be)_{\be<\al} $ satisfying conditions (CS1) and
(CS2). If $\al$ is a limit ordinal, then we set $\ga_\al=\sup_{\be<\al}\ga_\be$. Suppose now
that $\al$ is a successor ordinal. If $P_{[\ga_{\al^-},\ou)}|X$ is strictly singular, then we
set $\ga_\al=\ga_{\al^-}$. If not, let
\begin{align*}
\ga_\al=\min\conj{\ga\in (\ga_{\al^-},\ou)}{P_{[\ga_{\al^-},\ga)}|X\text{ is not strictly singular}}.
\end{align*}
Let us observe that if $X$ is separable, then the sequence $(\ga_\al)_{\al<\ou}$ is eventually
constant and we set $\Ga_X=\{\ga_\al\}_{\al<\ou}$. If $X$ is nonseparable, then the sequence
$(\ga_\al)_{\al<\ou}$ is strictly increasing  and $\Ga_X=\{\ga_\al\}_{\al<\ou}\cup \{\ou\}$.

Next we  prove the uniqueness of $\Ga_X$. Suppose on the contrary, and fix $\Ga\neq \Ga'$ two
different critical sets. Set $\ga=\max (\Ga \cap \Ga')$. First notice that $\max \Ga=\max
\Ga'$. So, either $\ga_{\Ga}^+<\ga_{\Ga'}^+$ or  $\ga_{\Ga'}^+<\ga_{\Ga}^+$. This  yields  a
contradiction using the fact that both $\Ga$ and $\Ga'$ satisfy (CS2).\fprue

\nota \noindent 1. The critical set $\Ga_X$ provides information concerning the structure of the
space $X$. For example the space $X$ is H.I. if  and only if $\Ga_X=\{0,\Om_X\}$. Also,
two subspaces $X,Y\hookrightarrow \eqs_{\ou}$ are totally incomparable  if and only if
$\Ga_X\cap \Ga_Y=\{0\}$.

\noindent 2.  For a transfinite block sequence $(x_\al)_{\al<\ga}$  its critical set is
nothing else but the set introduced  from Definition \ref{progammax} (2).
\fnota

\prop\label{subspacesubset} For every $Y\hookrightarrow X$, the corresponding critical set $\Ga_Y$
is a subset of $\Ga_X$.
\fprop
\prue This follows by an easy inductive argument. \fprue
\prop\label{blockasympto} For every separable $X\hookrightarrow\eqs_{\ou}$ and for every $\vep>0$
there exist an ordinal $\ga<\ou$,  a normalized sequence $(y_\al)_{\al<\ga}$ in $X$ and
a normalized transfinite
 block sequence $(z_\al)_{\al<\ga}$ such that (a) $\sum_{\al<\ga}\nrm{z_\al-x_\al}<\vep$ and
(b) $\Ga_X=\Ga_Z$ where $Z$ is the closed linear span of
 $(z_\al)_{\al<\ga}$.
\fprop
\prue
Use   Proposition \ref{blockseverywhere}, and a standard gliding hump argument.
\fprue

\defi Let  $X, Y\hookrightarrow \eqs_{\ou}$.
\begin{enumerate}
\item[(i)] We say that $X$ is \emph{asymptotically
finer} than $Y$, ${X\le_a Y}$, if and only if $\Ga_X\con \Ga_Y$.
\item[(ii)]We say that $X$ is \emph{asymptotically
equivalent} to $Y$, ${X\equiv_a Y}$, if and only if $\Ga_X= \Ga_Y$.
\end{enumerate}

\fdefi
It follows easily from the above definition that the relation $\le_a$ is a quasi ordering
in the class of the subspaces of $\eqs_{\ou}$ which from Proposition \ref{subspacesubset}
extends the natural inclusion. Notice also that $\equiv_a$ is an equivalence relation.

We now give  two alternative formulation of these notions.

\prop For $X,Y\hookrightarrow \eqs_{\ou}$ TFAE:

\noindent \emph{(1)} $X\le_a Y$,

\noindent \emph{(2)} if $P_I |X$ is not strictly singular, then  $P_I |Y$ is not strictly
singular, for every interval $I\con \ou$, and

\noindent \emph{(3)} $d(S_{X'},S_Y)=0$ for every $X'\hookrightarrow X$.
\fprop \prue Let us observe that for a closed infinite interval $I$, $P_I|X$ is not strictly
singular iff there is some $\ga_{\Ga_X}^+\in \Ga_X$ with $\min\Ga_X<\ga_{\Ga_X}^+\le \max I$.
The inverse direction follows immediately from the definition of the critical sets. So assume
now that $P_I |X$ is not strictly singular. Set $\ga_0=\max \conj{\ga\in \Ga_X}{\ga\le \min
I}$. Observe that $\ga_0\le \min I<\Om_X$, hence $\min\Ga_X<\ga_{\Ga_X}^+\le \max I$ by the
minimality  of $\ga_{\Ga_X}^+$ (Property (CS2)). It is easy to see that the above observation
implies easily the equivalence $(1)\Leftrightarrow (2)$. $(1) \Rightarrow (3)$: Suppose that
$X' \hookrightarrow X$. Then by Proposition \ref{subspacesubset} and our assumption,
$\Ga_{X'}\con \Ga_Y$. By Proposition \ref{blockasympto}, we can find two block sequences
$(z_n)_n$ and $(w_n)_n$ in $\eqs_{0^+_{\Ga_{X'}}}$ such that \emph{(a)} $\sup_n \max \supp z_n
= \sup_n \max \supp w_n=0^+_{\Ga_{X'}}$, and \emph{(b)} $d(S_{Z},S_X' )=d(S_{W},S_Y)=0$ where
$Z$ and $W$ are the closed linear span of $(z_n)_n$ and $(w_n)_n$ respectively. By Corollary
\ref{dis0}, $d(Z,W)=0$ and we are done. $(3) \Rightarrow (2)$: Since for every
$X'\hookrightarrow X$, $d(S_{X'},S_Y)=0$, we obtain that for every $\vep>0$, and every
$X'\hookrightarrow X$ there exists two   basic sequences $(z_n)_n$ and $(w_n)_n$ such that
$z_n\in S_{X'}$ and $w_n\in S_Y$ for all $n$ and $\sum_n \nrm{z_n-w_n}<\vep$. Assume now that
$P_I |X$ is not strictly singular. Choose $X'\hookrightarrow X$ such that $P_I |X'$ is an
isomorphism. Let $(z_n)_n\con X'$ and $(w_n)_n\con Y$ as above. Then $P_I|W$ is  isomorphism
and hence $P_I|Y$ is not strictly singular. \fprue

\prop
For $X,Y\hookrightarrow \eqs_{\ou}$ the following are equivalent: \emph{(1)} $X\equiv_a
Y$, \emph{(2)} $P_I |X$ is not strictly singular if and only if
 $P_I |Y$ is not strictly singular, for every interval $I\con \ou$, and\emph{ (3)}
$d(S_{X'},S_Y)=d(S_{Y'},S_X)=0$ for every $X'\hookrightarrow X, \,Y'\hookrightarrow Y $.
\qed\fprop

\cor \noindent (1) For every $X\hookrightarrow \eqs_{\ou}$ and every $A\con \Ga_X$ there
is $X_A\hookrightarrow X$ such that $\Ga_{X_A}=A$.

\noindent (2) For any nonseparable $X,Y \hookrightarrow \eqs_{\ou}$  there are
nonseparable $X_1\hookrightarrow X$, $Y_1\hookrightarrow Y$ such that $X_1\equiv_a Y_1$.
\qed \fcor

We shall need the following consequence of a well known result from Lindenstrauss
\cite{L}.
\lema\label{Lindes} Let  $Z\hookrightarrow X\hookrightarrow \eqs_{\ou}$ with
$Z$ separable. Then there exist a separable subspace $W$ of $X$  and $\ga<\ou$  such that $ Z
\hookrightarrow W $ and $P_\ga|X$ is a projection onto $W$. \qed
\flema
\nota Notice that for
$W$ and $X$ as in the Lemma, $\Ga_W$ is an initial part of $\Ga_X$.
\fnota
\prop\label{pretheor} Let $X$ be a subspace of $\eqs_{\ou}$  and
$T:X\to \eqs_{\ou}$ a bounded operator. Then there exists a unique $D_T \in \mc
D_{\Ga_X}(\eqs_{\ou})$ such that (a) $\nrm{D_T}\le 2\boldsymbol{C}^2 \nrm{T}$ and (b)
$T-D_T|X$ is strictly singular.
\fprop \prue Fix $X \hookrightarrow \eqs_{\ou}$ and a bounded operator $T:X\to
\eqs_{\ou}$. First suppose that $X$ is separable. Then we can find a transfinite basic
sequence $(y_\al)_{\al<\ga}\con X $ and a transfinite block sequence $(z_\al)_{\al<\ga}$
of $\eqs_{\ou}$ such that $\sum_{\al<\ga}\nrm{y_\al-z_\al}<1$ and $X\equiv_a Z$, where
$Z$ denotes the closed linear span of $(z_{\al})_{\al<\ga}$. Consider now
$T':Z\overset{U}\to Y \overset{T|Y}\to \eqs_{\ou} $ where $Y$ is the closed linear span
of $(y_\al)_{\al<\ga}$ and $U:Y\to Z$ is the isomorphism defined by $U(\sum_{\al<\ga}
a_\al z_\al)=\sum_{\al<\ga} a_\al y_\al$. Notice that $\nrm{U}\le 2$. Then there is
unique $D\in \mc D (Y)$ such that $T'-D$ is strictly singular, or equivalently there is
unique $D_{T'}\in \mc D_{\Ga_Z}(\eqs_{\ou})$ such that $T'-D_{T'}|Z$ is strictly
singular. Notice that $\nrm{D_{T'}}\le \boldsymbol{C}\nrm{D}\le
\boldsymbol{C}^2\nrm{T'}\le\boldsymbol{C}^2\nrm{U}\nrm{T}\le 2\boldsymbol{C}^2\nrm{T}$.
 Let us show that $T-D_{T'}$ is
strictly singular. Let $X'\hookrightarrow X$ and $\vep>0$.

Choose $Z'\hookrightarrow Z$ such that $(\Ga_{Z'}\setminus \{0\})\cap (\Ga_{X'}\setminus
\{0\})\neq \buit$,   $\nrm{U|Z'-i_{Z',\eqs_{\ou}}}\le\vep/(4\nrm{T})$ and
$\nrm{(T'-D_{T'})|Z'}\le \vep/4$. Pick $z'\in Z'$ and $x'\in X'$ such that $\nrm{z'-x'}\le \vep
/(2 (\nrm{D_{T'}}+\nrm{T}))$. Then
\begin{align}
\nrm{(T-D_{T'})x'}\le&  \nrm{(T-D_{T'})x' -(T'-D_{T'})z'}+\nrm{(T'-D_{T'})z'} \le
  \nrm{T}\nrm{x'-Uz'}+ \nonumber\\
+ & \nrm{D_{T'}}\nrm{x'-z'}+\frac\vep 4 \le
(\nrm{T}+\nrm{D_{T'}})\nrm{x'-z'}+\frac{\vep}2 \le \vep.
\end{align}
Now suppose that $X$ is  nonseparable. By Lemma \ref{Lindes}, we can find a sequence
$(X_\ga)_{\ga<\ou}$ of separable complemented subspaces of $X$ such that $\Ga_{X_\ga}$ is an
initial part of $\Ga_X$ for every $\ga<\ou$. Now the result for $X$  follows easily from the
result for the corresponding $T_\ga=T|X_\ga$ and the fact that $D_T\in \mc
D_{\Ga_X}(\eqs_{\ou})$ and $D_{T_\ga}\in \mc D_{\Ga_{X_\ga}}(\eqs_{\ou})$ are unique.  The
uniqueness of $D_T\in \mc D_{\Ga_X}(\eqs_{\ou})$ is clear from the analogous result for
transfinite block sequences. \fprue
\teor $\mc L(X,\eqs_{\ou}) \cong \mc D_{\Ga_X}(\eqs_{\ou}) \oplus \mc S(X,\eqs_{\ou})
\cong J_{T_0}^*(\Ga_X^{(0)} )\oplus\mc S(X,\eqs_{\ou}) $ for every $X\hookrightarrow
\eqs_{\ou}$. If in addition $\Ga_X$ is infinite, then $\mc L(X,\eqs_{\ou})  \cong
J_{T_0}^*(\Ga_X )\oplus\mc S(X,\eqs_{\ou})$. \fteor \prue Let $H:\mc D_{\Ga_X}\to \mc L
(X,\eqs_{\ou})$ be defined by $D\mapsto D|X$. Assume first that $X$ is separable.  It is
clear that $\nrm{D|X}\le \nrm{D}$. For an appropriate $\vep'>0$, we can find normalized
$(y_\al)_{\al<\ga}$ and a normalized block sequence $(z_\al)_{\al}$ such that
$\Ga_X=\Ga_Z$ and $\sum_{\al}\nrm{z_\al-y_\al}\le \vep'$ where $Z$ the closed linear span
of $(z_\al)_{\al<\ga}$. Since by Theorem \ref{extensiondiagonal} $\nrm{D|Z}\ge
\nrm{D}/\boldsymbol{C}$,  we get that
\begin{equation}
\nrm{D}/\boldsymbol{C}\le \nrm{D|Z}\le (1+\vep)\nrm{D|Y}\le (1+\vep) \nrm{D|X}=(1+\vep)\nrm{H(D)}.
\end{equation}
Hence, $H$ defines an isomorphism. To show that $H$ is an isomorphism when $X$ is nonseparable
we use a family $(X_\al)_{\al<\ou}$ of separable complemented subspaces of $X$ defined as in
the previous proof.  Proposition \ref{pretheor} shows that $\mc L(X,\eqs_{\ou})\cong \mc
D_{\Ga_X}(\eqs_{\ou})\oplus \mc S(X,\eqs_{\ou})$.

For the later isomorphism see Remark \ref{dfdlwerijghgg}. \fprue

\subsection{Examples with
$\boldsymbol{\mathcal{L}(X)/\mathcal{S}(X)
\ncong\mathcal{L}(X,\eqs_{\ou})/\mathcal{S}(X,\eqs_{\ou})}$} We present a family $\{{\mc
Z}_\zeta\}_{\zeta<\ou}$ of separable subspaces of $\eqs_{\ou}$ such that each ${\mc
Z}_\zeta$  is indecomposable but  has a $\zeta$-closed direct sum as a subspace.
\defi
For given $\al\le \be<\ou$, let $d_{\al,\be}=e_\al+e_\be$. Given $A=\{\al_n\}_n
\uparrow,B=\{\be_n\}_n\uparrow\con \ou$   such that $A<B$, let ${\mc Z}_{A,B}$ be the closed
linear span generated by $\{d_{\al_n,\be_n}\}_n$.
\fdefi
\prop \label{lots} $\Ga_{{\mc Z}_{A,B}}=\{0,\al,\be\}$ where $\al=\sup A, \, \be=\sup B$.
\fprop
\prue We get the direct inclusion above, since $\mc Z_{A,B}\con \eqs_{A\cup B}$. It remains to
show that $P_{\al}|\mc Z_{A,B}$ and $P_{(\al,\be)}|\mc Z_{A,B}$ are not strictly singular.  We
check the case of $P_{\al}|\mc Z_{A,B}$ since the other is similar. Let $U:\eqs_A\to \mc
Z_{A,B}$ be the linear map defined by $e_{\al_n}\mapsto d_n$. Since $\lim_n d(U e_{\al_n}, \R
e_{\al_n})\ge 1$, we can apply Proposition \ref{keyforbridges}, and we can obtain a block
sequence $(x_n)_n$ such that $\nrm{U x_n}=1$ and $\nrm{x_n}<1/2^n$ for every $n$. Now
$\nrm{P_\al Ux_n}\ge \nrm{U x_n }-\nrm{x_n}\ge 1/2$ for every $n$. Hence, $P_\al | X$ is an
isomorphism where $X$ is the closed linear span of the Schauder basic sequence $(Ux_n)_n$.
\fprue
\nota\label{subsequequiv} Note that this shows that $\mc Z_{A',B'}\equiv_a \mc Z_{A, B}$ for every
infinite $A'\con A$, $B'\con B$. \fnota

\prop \label{lidforz}
Suppose that  $T:\mc Z_{A,B}\to \eqs_{\ou}$  is  bounded and satisfies  for every $n,m$
\begin{equation}  \label{symmd}
e_{\al_n}^*T d_m =   e_{\be_n}^*T d_m.
\end{equation}
Then there is some scalar $\la$ such that $T-\la i_{\mc Z_{A,B},\eqs_{\ou}}$ is strictly
singular.  Consequently, every bounded operator $T:Z\to Z$ is of the form $T=\la Id_Z+S$,
where $S$ is strictly singular.  Hence, $Z$ is indecomposable. \fprop \prue Let $T:\mc
Z_{A,B}\to \eqs_{A,B}$ be bounded and satisfying (\ref{symmd}). Let $d_n=d_{\al_n,\be_n}$
for every $n$. \clam $\lim_{n\to \infty} d(Td_n, \R d_n)=0$. \fclam \prucl Condition
$(\ref{symmd})$ implies that
\begin{equation}
\max\{\inf_{n}d(P_\al T d_n,\R e_{\al_n}),\inf_{n}d(P_{[\al,\ou)} T d_n,\R e_{\be_n})\}>0.
\end{equation}
Without loss of generality we may assume that $\inf_n d(P_\al T d_n,\R e_{\al_n})> 0$.
Applying Proposition \ref{keyforbridges} to $U=P_\al T :\langle e_{\al_n}\rangle_n\to
\eqs_{\ou}$, we can find $x=\sum_{k\in F} a_k e_{\al_k}\in \langle e_{\al_n}\rangle_n$
such that $\nrm{x}<(1/3\nrm{T})\nrm{U x}$ and  $\nrm{\sum_{k\in F} a_k e_{\be_k}}\le
(1/3\nrm{T})\nrm{U x}$. This implies that $\nrm{\sum_{k\in F}a_k d_k} \le
(2/{3\nrm{T}})\nrm{T (\sum_{k\in F}a_k d_k)}\le  (2/3)\nrm{\sum_{k\in F}a_k d_k}$,  a
contradiction. \fprucl Now for each $n$, let $\la_n\in \R$ realizing $d(T d_n,\R
d_n)=\nrm{T d_n-\la_n d_n}$, and choose any accumulation point $\la$ of $(\la_n)_n$. Let
us show that $S=T-\la_{\mc Z_{A,B}}$ is strictly singular. Fix $\vep>0$, and let $N\con
\N$ be infinite such that $\la_n \to_{n\to\infty, n\in N} \la$ and $\nrm{T-\la i_{\mc
Z_{A',B'}}}\le \vep/2$, where $A'=\{\al_n\}_{n\in N}$, $B'=\{\be_n\}_{n\in N}$. Notice
that from Remark \ref{subsequequiv} we know also that $\mc Z_{A',B'}\equiv_a \mc
Z_{A,B}$. So, given any $X \hookrightarrow \mc Z_{A,B}$, we can find normalized $x\in X$,
$y\in \mc Z_{A',B'}$ with $\nrm {x-y}\le \vep/2 \nrm{S}$. Hence, $\nrm{S x }\le \nrm{S y}
+\nrm{S (x-y)}\le \vep$ as desired.
\fprue

We generalize the previous ideas and we present a family $\mc Z_\ga$ ($\ga<\ou$, $\ga$ limit)
of infinite dimensional closed subspaces of $\eqs_{\ou}$ such that for every limit ordinal
$\ga$, $\mc Z_\ga\equiv_a \eqs_\ga$ and such that $\dim \mathcal{L}(\mc Z_\ga)/\mathcal{S}(\mc
Z_\ga)=1$. In particular, each $\mc Z_\ga$ is an indecomposable space.

\defi
Fix a limit ordinal $\ga<\ou$. Let $\mc I_\ga$ be the family of minimal infinite intervals of
$\ga$, i.e.,  $\mc I_\ga=\conj{[\al,\al+\om)}{\al\text{ is a limit ordinal }, \, \al+\om\le \ga}$.
For each $I\in \mc I_\ga$, we choose a partition $\conj{L^I_J\con I}{J\in \mc I_\ga}$ into
infinite sets. Notice that since $I=[\al,\al+\om)$ for some limit $\al$, all infinite sets $L^I_J$
have order type $\om$. Now for each $n\in \N$ we  consider the vectors
$d_n^{I,J}=e_{\al_n}+e_{\be_n}$,  where $\{\al_n\}_n$ and $\{\be_n\}_n$ is the increasing
enumeration of the sets $L_J^I$ and $L_I ^J$ respectively. Finally, let $\mc Z_\ga$ be the closed
linear span of $(d_n^{I,J})_{I,J\in \mc I_\ga,n\in \N}$.
\fdefi
\teor
$\mc L(\mc Z_\ga,\eqs_{\ou})/\mc S(\mc Z_\ga,\eqs_{\ou})\cong J_{T_0}^*(\La(\ga))$ and $\dim
\mathcal{L}(\mc Z_\ga)/\mathcal{S}(\mc Z_\ga)=1$, for a limit ordinal  $\ga$.
\fteor
\prue Notice that for every limit ordinal $\al$ such that $\al+\om\le\ga$ we have that
$d^{[\al,\al+\om),[\al,\al+\om)}_n=2 e_{\al_n}\in \mc Z_{\xi}$, where
$L_I^I=\{\al_n\}_n\uparrow$. This together with the fact that $\mc Z_\ga \hookrightarrow
\eqs_{\ga}$ gives that $\mc Z_\ga \equiv_a  \eqs_\ga$ (i.e, $\Ga_{\mc Z_\ga}=\La(\ga+1)$) and
hence $\mc L(\mc Z,\eqs_{\ou})/\mc S(\mc Z,\eqs_{\ou})\cong J_{T_0}^*(\La(\ga+1)^{(0)})\cong
J_{T_0}^*(\La(\ga))$.

Fix now a bounded operator $T:\mc Z_{\zeta}\to \mc Z_{\zeta}$. By Proposition \ref{pretheor} there is $D\in
\mc D(\eqs_{\ga})$ such that $T-D |\mc Z_\ga$ is strictly singular. The proof will finish by proving that
$\xi_T$ is constant. We observe that given $I<J\in \mc I_\ga$, $I=\{\al_n\}_n$, $J=\{\be_n\}_n$ increasing
enumeration, we have that $e_{\al_n}^*T d_m^{I,J}= e_{\be_n}^*T d_m^{I,J}$ for every $n, m$. So from
Proposition \ref{lidforz} we obtain that  for every pair $I<J$ in $\mc I_\ga$ there is some $\la_{I,J}$ such
that $T|\mc Z_{I,J}- \la_{I,J} i_{\mc  Z_{I,J},\eqs_{\ou}}$  is strictly singular, and this clearly implies
that $\xi_T$ is constant.
\fprue

\nota Notice that it is not possible to improve the previous result to a nonseparable subspace
of $\eqs_{\ou}$, since every nonseparable reflexive space admits non trivial projections \cite{L}.
\fnota

\subsection{Further results on Operators}

\cor No closed linear span $X$ of a transfinite block sequence of $\eqs_{\ou}$ is isomorphic to
finite cartesian power of a Banach space. \fcor
\prue
This is so since $\mathcal{L}(X)$ admits a non trivial linear multiplicative functional.
\fprue

Recall the following  facts about semi-Fredholm operators (see \cite{go-ma2}, \cite{li-tza})
\prop Suppose that $T:X\to Y$ bounded  such that $TX$ is closed and $\al(T)<\infty$. Then there is
some number $\vep(T)>0$ such that if $S:X\to Y$ is bounded and satisfying that for any $X_1
\hookrightarrow X$ there is some $x\in S_{X_1} $ with $\nrm{S(x)}<\vep$, then $T+S$ has closed
range and $\al(T+S)<\infty$. \fprop \prue Since $\ker T$ is finite dimensional, $X=X_1\oplus \ker
T$. Let $T_1=T|X_1$. Notice that $T_1|X_1= TX_1=TX \hookrightarrow Y$ is closed, and  therefore
$T_1:X_1\to TX_1$ is an isomorphism. Let $\vep=\vep(T)=(1/2)\nrm{T_1^{-1}}^{-1}$. Fix now $S$
satisfying the condition about $\vep$. Suppose that $T+S$ has $\al(T+S)=\infty$. Then, $T_1+S|X_1$
has infinite dimensional kernel.\fprue

\prop\label{2c9}
Suppose that $T:X\to Y$ is semi-Fredholm.
\begin{enumerate}
\item Then there is some number $\vep=\vep(T)>0$ such that for all $S:X\to Y$ with $\nrm{S}<\vep$, $T+S$ is
semi-Fredholm and $i(T+S)=i(T)$.
\item If $\al(T)$ finite, and  $S:X\to
Y$ is a strictly singular operator, then $T+S$ is semi-Fredholm, $\al(T+S)$ is finite and
$i(T+S)=i(T)$.\qed
\end{enumerate}
\fprop In the next results $X$  denotes the closed linear span of a transfinite block
sequence $(x_\al)_{\al<\ga}$ of $\eqs_{\ou}$.
\prop \label{ontoclosed} Suppose that $D\in
\mc D(X)$ is such that $\inf\conj{\xi_D(\al)}{\al\in \La(\ga+1)^{(0)}}>0$. Then $D$ is a
Fredholm operator with index 0, and hence it is an onto isomorphism.
\fprop
\prue Let $\widetilde{\xi}_D:\La(\ga+1)\to \R$ be the unique continuous extension of $\xi_D$.
Notice that the above inequality is equivalent to say that $\widetilde{\xi}_D$ is never zero.
In order to show that $D$ is an onto isomorphism it is enough to show that $DX$ is closed. If
not, for every $n$ we can find an
 block sequence $X_n \hookrightarrow \eqs_{\ga}$ such that $\nrm{D|X_n}\le 2^{-n}$.
  Notice that  for every $n$,  $D
|X_n -\widetilde{\xi}_D(\ga_n) i_{X_n,\eqs_\ga} \text{ is strictly singular}$, where $\ga_n=\max_n
\conj{\supp x}{x\in \langle X_n\rangle}$. Now for all $n$, we can find a norm 1 normalized vector
$x_n\in X_n$ such that $\nrm{Dx_n -\widetilde{\xi}_D(\ga_n)x_n}<2^{-n}$, and hence
$|\widetilde{\xi}_D(\ga_n)|\le 2^{1-n}$ for every $n$. Continuity of $\xi^e_D$  implies that there
is some limit $\al\le \ga$ such that $\widetilde{\xi}_D(\al)=0$, a contradiction.\fprue

\teor\label{semi=full}
$T\in \mathcal{L}(X)$ is Fredholm with index $i(T)=0$ iff it is semi-Fredholm.
\fteor

\prue
 Suppose that
$T:X\to X$ is Semi-Fredholm. Let us take the decomposition $T=D_T+S$. If $\al(T)$ is finite,
then since $S$ is strictly singular, by Proposition \ref{2c9} (2), $D_T$ is semi-Fredholm with
finite dimensional kernel and with the same index as $T$. This implies that $\xi_T$ is never
zero (otherwise, the kernel is infinite dimensional), and hence $D_T$ is indeed 1-1. Since for
all $\al<\ga$, $x_\al\in D_T \eqs_{\ga}$, and $D_T\eqs_{\ga}$ is closed, we get that $D_T$ is
an onto isomorphism. Hence $T$ is Fredholm with index 0.

Suppose now that $\be(T)$ is finite. Let $\vep>0$ be given by Proposition \ref{2c9} (1), and
let $\la\in (-\vep,\vep)\setminus \widetilde{\xi}_T (\La(\ga +1))$. Notice that this is
possible since $\widetilde{\xi}_T (\La(\ga +1))$ is countable by the fact that
$\widetilde{\xi}_T :\La(\ga+1)\to \R$ is continuous. Then $T'=T-\la Id_X$ is Semi-Fredholm
with the same index as $T$ and $\widetilde{\xi}_{T'}$ is never zero. So, $D_{T'}$ satisfies
that $\widetilde{\xi}_{D_{T'}}$ is never zero. By the previous Proposition \ref{ontoclosed},
$D_{T'}$ is an isomorphism onto. Hence $T'$ is Fredholm with index 0 and $i(T)=i(T')=0$.
\fprue

\cor $X$ is not isomorphic  to  either a
 proper  subspace of it, or to a non trivial quotient.
\fcor
\prue Let $Y \hookrightarrow X$. Suppose first that $T: X\to Y$ is an onto isomorphism. Then
the composition $U=i_{Y,X}\circ T:X\to X$ is a semi-Fredholm operator, with $\al(T)=0$. By
Theorem \ref{semi=full}, $U$ is indeed Fredholm with index 0, hence $U$ is onto, i.e, $X=U
X=Y$.

Suppose now that $T:X/Y\to X$ is an onto isomorphism. Now the composition $U=T\circ \pi_Y:X\to
X$ is semi-Fredholm and onto, where $\pi_Y:X\to X/Y $ is the quotient mapping. Again $U$ has
to be Fredholm with index 0, hence  $U$ is 1-1, i.e.,  $Y=\ker U=\{0\}$. \fprue

\prop
$\overline{ \langle v_{2n}^*\rangle_n }\cong T_0^*$.
\fprop
\prue
Fix a sequence of scalars $(b_n)_n$. Let $\sum_n a_n v_n\in J_{T_0}$ be of norm 1 such that
\begin{equation}
\nrm{\sum_n b_n v_{2n}^*}_{J_{T_0}^*}=(\sum_n b_n v_{2n}^*)(\sum_n a_n v_n)= \sum_n b_n a_{2n}.
\end{equation}
Since $\nrm{\sum_n a_{2n}t_{n}}_{T_0}=\nrm{\sum_n a_{2n}t_{2n}}_{T_0}\le \nrm{\sum_n a_n
v_n}_{J_{T_0}}=1 $, it follows that $\nrm{\sum_n b_n t_n^*}_{T_0^*}\ge (\sum_n b_n t_n^*)
(\sum_n a_{2n} t_n )= \sum_n b_n a_{2n}= \nrm{\sum_n b_n v_{2n}^*}_{J_{T_0}^*}$.  The
other inequality follows from the fact that $(-v_{2n-1}+v_{2n})_n$ is equivalent to
$(t_n)_n$ (by Proposition \ref{sdfewooo1}). \fprue

\prop There are $X,Y\hookrightarrow \eqs_{\ou}$ such that  $\mc
L(X)/\mc S(X)\cong \mc L(Y)/\mc S(Y) \cong J_{T_0}^*$ and such that  $\mc L(X,Y)/\mc S (X,Y) \cong
T_0^*$.
\fprop

\prue

Let $X=\eqs_{\om^2}$, and let $Y={\eqs_A}$, where $A=\bigcup_n [\om (2n),\om(2n+1))$. The result
follows from the fact that if $T:X\to Y$ is  a bounded operator, then  necessarily $\xi_{T'}(\om
(2n))=0$, for $T'=i_{Y,X}\circ T= D+S$. And $\mc L(X,Y)/\mc S (X,Y)\cong \conj{\xi\in
J_{T_0}^*}{\xi(v_{2n+1})=0 }=\overline{ \langle v_{2n}^*\rangle_n }\cong T_0^*$.\fprue

 \prop Every  $T\in \mathcal{L}(\eqs_{\ou})$ is of the form $T=\la Id +R$, where $R$ has separable range. \fprop \prue We know
that $T=D_T+S$ with $S$ strictly singular. Corollary \ref{mat1} shows that $S\in
c_0(\ou)$. Since $\xi_T:\La(\ou)^{(0)}\to \R$ has the continuous extension property,
there is some $\la\in \R$ and some ordinal $\al<\ou$ such that $\xi_T(\be)=\la$ for every
$\be \in \La(\ou)^{(0)}\cap [\al,\ou)$. Hence $D_T-\la Id$ has separable range. So $T-\la
Id=(D_T-\la Id)+S$ has separable range, as desired. \fprue Nonseparable spaces with this
property of operators have   been constructed before in  \cite{sh},  \cite{sh-st} and
\cite{wark}. These constructions however give no control on separable subspaces.

The following theorems summarize our results for the structure of
$\eqs_{\omega_1}$, its subspaces and the spaces of operators.

\teor \noindent There exists a reflexive space $\eqs_{\omega_1}$
with a transfinite basis $(e_{\al})_{\al<\omega_1}$ such that

\noindent(1) It does not contain an unconditional basic sequence.

\noindent(2) It is arbitrarily distortable.

\noindent(3) $\eqs_I$ and $\eqs_J$ are totally incomparable for
disjoint infinite intervals $I$ and $J$.

\noindent(4) It is $\omega_1$ hereditarily indecomposable (i.e., for every nonseparable
$X,Y \hookrightarrow \eqs_{\omega_1}$, $dist(S_X,S_Y)=0$).

\noindent(5) Every subspace $X \hookrightarrow \eqs_{\ou}$
generated by a transfinite block sequence is, neither isomorphic
to a proper subspace, nor to a non-trivial of its quotients.
\fteor

To each infinite dimensional subspace $X$ of $\eqs_{\omega_1}$ we
assign a closed subset $\Gamma_X$ of $\omega_1$, called the
critical set of $X$. The following theorem describes the
interference of $X$ and $\Gamma_X$.

 \teor For $X,Y$ subspaces of $\eqs_{\omega_1}$ the following holds

\noindent(1) If $Y\hookrightarrow X$ then
$\Gamma_Y\subset\Gamma_X$.

\noindent(2) The subspaces $X,Y$ are totally incomparable iff
$\Gamma_X\cap\Gamma_Y=\{0\}$.

\noindent(3) The subspace $X$ is hereditarily indecomposable iff
$\# \Gamma_X=2$.

\noindent(4) For every subspace $X$ of $\eqs_{\omega_1}$ there
exists $Y$ generated by a block sequence $(y_{\al})_{\al<\gamma}$
such that $\Gamma_X=\Gamma_Y$.
 \fteor

 Finally the structure of the spaces of  operators is described
 by the next theorem. Recall that
 $\Gamma_X^{(0)}$ is the set of isolated ordinals of $\Gamma_X$.

 \teor \noindent (1)  For every $X\hookrightarrow
\eqs_{\ou}$, $\mc L(X,\eqs_{\ou}) \cong \mc D_{\Ga_X}(\eqs_{\ou})
\oplus \mc S(X,\eqs_{\ou}) \cong J_{T_0}^*(\Ga_X^{(0)} )\oplus\mc
S(X,\eqs_{\ou}) $. If in addition $\Ga_X$ is infinite, then $\mc
L(X,\eqs_{\ou})\cong J_{T_0}^*(\Ga_X )\oplus\mc S(X,\eqs_{\ou}) $.

\noindent (2)  For every $X \hookrightarrow \eqs_{\ou}$ generated
by a transfinite block sequence, $\mc L(X)  \cong
J_{T_0}^*(\Ga_X^{(0)})\oplus \mc S(X) $. If in addition $\Ga_X$ is
infinite, $\mc L(X) \cong J_{T_0}^*(\Ga_X)\oplus \mc S(X)$.

\noindent (3) For every $\gamma<\omega_1$ there exists a subspace $Y_{\gamma}$ of
$\eqs_{\omega_1}$ such that $\mc L(Y_{\gamma}) \cong \langle \mathrm{Id}_{Y_{\gamma}}
\rangle \oplus \mc S(Y_{\gamma})$ and $\mc L(Y_{\gamma},\eqs_{\omega_1})\cong
J_{T_0}^*(\gamma) \oplus \mc S(Y_{\gamma},\eqs_{\omega_1})$.
 \fteor

\section{ Universal and smooth $\ro$-functions}\label{difros}
In this section we present two new properties which a $\ro
-$function can have. In this and in
the subsequent section we show how these properties of $\ro$
influence the corresponding space  $\eqs_{\ou}$ based on $\sigma_{\rho}$.
\subsection{The construction of a universal $\ro$-function}\label{universal}
In this subsection we  show how the construction of the $\varrho$-function of \cite{tod1} can be
adjusted in order to give us a function $\varrho:[\ou]^2\to \omega$ with the following properties:
\begin{enumerate}
\item $\varrho(\al,\ga)\le \max\{\varrho(\al,\be),\varrho(\be,\ga)\}$ for all
$\al<\be<\ga<\ou$.
\item $\varrho(\al,\be)\le \max\{\varrho(\al,\ga),\varrho(\be,\ga)\}$ for all
$\al<\be<\ga<\ou$.
\item $\conj{\al<\be}{\varrho(\al,\be)\le n}$ is finite for all $\be<\ou$ and $n\in \N$.
\item $(\ou,\varrho)$ is universal.
\end{enumerate}
To describe what we mean by ``$(\ou,\varrho)$ is universal" we need some more definitions.
\defi
A finite $\varrho$-model is a model of the form $(M,<,\varrho_M,p_M)$ where $M$ is a set,
$<$ is a total ordering on $M$, $p_M$ is an integer and $\varrho_M:[M]^2\to p_M$ is a
function with properties 1. and 2. listed above. We also assume that there
exist $x<y$ in $M$ such that $\varrho_M(x,y)=p_M$.
\fdefi
\defi
Suppose that $\varrho:[\la]^2\to \omega$ satisfies 1., 2., and 3.
above. For $M\subseteq \la$, let $p_M=\max\{\varrho(\al,\be)\, : \, \al,\be\in M \}$. Such an $M$ is
$\varrho$-closed if
$$M=\{ \al<\la \, :\, \exists \be\in M \,(\al\le \be \, \& \varrho(\al,\be)\le p_M\,)\}.$$
\fdefi

We use the convention of $\varrho(\al,\al)=0$ for all $\al$. Note that for a
$\ro$-closed subset $M$ of $\la$, $(M,<,\varrho| [M]^2,p_M)$ is an example of a
$\ro$-model. Note also that an initial part $M_0$ of a $\ro$-closed set $M$ is a
$\ro$-closed set and that its integer $p_{M_0}$ might be smaller than $p_M$. Similarly,
an initial part of a $\ro$-model is also a $\ro$-model with a possibly smaller
integer $p_M$.

\defi\label{rhoisomophism}
Two $\ro$-models $(M_1,<_1,\varrho_1,p_1)$ and $(M_2,<_2,\varrho_2,p_2)$ are isomorphic if $p_1=p_2$
and if there is order-isomorphism $\pi:(M_1, <_1)\to (M_2,<_2)$ such that $\ro_1(a,b)=\ro_2(\pi(a),\pi(b))$
for all $a,b\in M_1$.
\fdefi

\defi
A function $\ro:[\la]^2\to \omega$ defined on some limit ordinal $\la\le \ou$ and satisfying 1., 2. and 3. is
said to be \emph{universal} if for every finite $\ro$-model $(M,<,\ro_M,p_M)$, every $\ro$-closed subset
$M_0$ of $\la$ such that $(M_0,<,\ro |[M_0]^2,p_{M_0})$ is isomorphic to an initial segment of
$(M,<,\ro_M,p_M)$, and every ordinal $\de$ such that $\de+\omega \le \la$, there is a $\ro$-closed subset
$M_1$ of $\de+\omega$ such that
\begin{enumerate}\addtocounter{enumi}{4}
\item $(M_1,<,\ro |[M_1]^2,p_{M_1})\cong (M,<,\ro_M,p_M)$
\item $M_0$ is an initial segment of $M_1$.
\item $M_1\setminus M_0\con [\de,\de+\omega)$.
\end{enumerate}
\fdefi
The existence of a  universal $\ro:[\ou]^2\to \omega$ is established by recursively
constructing an increasing sequence $\ro_\la:[\la]^2\to \omega$  ($\la\in \La$). Let
$\ro_0=\buit$, and suppose $\ro_\la: [\la]^2\to \omega$ has been determined for some
countable limit ordinal $\la$. Let $C$ be a subset of $\la$ of order-type $\omega$ such
that $\la=\sup C$. Define
$$  \ro_{\la+\omega}(\al,\la)=\max\{\# (C\cap \al), \ro_\la(\al,\min (C\setminus \al)) ,
\ro_\la(\xi,\al) \, : \, \xi\in C\cap \al \}.$$

It can be checked (see e.g. \cite{tod1}) that this defines a function
$\ro_{\la+\omega}:[\la+1]^2\to \om$ satisfying the conditions 1., 2. and 3. Starting with this
extension of $\ro_\la$ and the assumption that $\ro_\la$ is  universal we build extensions
$\ro_{\la+\om}:[\de]^2\to \omega$ ($\la+1\le \de <\la+\om$) in such a way that at a given
stage $\de$ we take care about a particular instance of universality of $\ro_{\la+\om}$. Thus,
modulo some book keeping device, it suffices to show how one deals with the following task: We
have already defined an extension $\ro_{\la+\om}:[\de]^2\to \om$, we are given a finite
$\ro$-model $(M, <, \ro_M,p_M)$ and a $\ro$-closed  subset $M$ of $\de$ such that $(M_0, <,
\ro_{\la+\om}|[m_0]^2 ,p_{M_0})$ is isomorphic to a proper initial segment of $(M, <,
\ro_M,p_M)$. Let $l=\# M-\# M_0$ and extend $\ro_{\la+\om}$ from $[\de]^2$ to $[\de+l]^2$ as
follows: First of all define $\ro_{\la+\om}$ on
$$  [M_0\cup [\de,\de+l)]^2\setminus [M_0]^2$$
in such a way that we have the isomorphism
$$  (M_0\cup [\de,\de+l),<,\ro_{\la+\om},p_M)\cong (M,<,\ro_M,p_M).$$
Thus, it remains to define $\ro(\al,\ga)$ for $\al\in \de\setminus M_0$ and $\ga\in (\de,\de+l)$.
If $\al<\de$ and $\al>\max M_0$, then set
$$  \ro (\al,\ga)=\max\{p_M+1,\ro(\al,\de-1),\ro(\xi,\al)\,:\, \xi\in M_0\}.$$
If $\al\le \max M_0$, then set
$$  \ro (\al,\ga)=\max\{\ro(\al,\min (M_0\setminus \al)), \ro(\xi,\al)\,:\, \xi\in M_0\cap \al\}.$$

It remains to show that $\ro_{\la+\om}|[\de+l]^2$ satisfies the properties 1. and 2. So
let $\al<\be<\ga<\de+l$ be given. We simplify the notation by writing $\al\be$ instead of
$\ro(\al,\be)$, $\xi \eta \vee \al\be$ instead $\max\{\ro(\xi,\eta),\ro(\al,\be)\}$, and
$\de^-$ in place of $\de-1$.

\noindent{\emph{Case 1.}} $\al<\de\le \be <\ga <\de+l$. If $\al\in M_0$, then  properties 1.
and 2. for $\al<\be<\ga$ follow from the fact that in the definitions of $\al\be$, $\be\ga$
and $\al\ga$ we have copied the $\ro$-model $(M,<,\ro_M,p_M)$ which satisfies 1. and 2.

If $\al\notin M_0$, then in both the case $\al>\max M_0$ and $\al\le \max M_0$ we conclude
that $\al\be=\al\ga$, so 1. and 2. for $\al<\be<\ga$ follow immediately.

\noindent{\emph{Case 2.}} $\al<\be<\de<\ga<\de+l$. \noindent \emph{Subcase 2.1.} $\max
M_0<\de<\be$. Consider first the inequality $\al\ga\le \al\be \vee \be\ga$. The quantities $p_M+1$
and $\xi \al$ ($\xi\in M$) from the definition of $\al\ga$ are all present in the definition of
$\be\ga$, so it remains only to show that the quantity $\al\de^-$ is bounded by $\al\be\vee
\be\ga$. Applying 1. for $\ro_{\la+\om} |[\de]^2$ we get
$$  \al\de^- \le \al\be \vee \be\de^-,$$
and so we are done as $\be\de^-$  shows up in the definition of $\be\ga$. Consider now
the inequality $\al\be\le \al\ga\vee\be\ga$. Applying 2. for $\ro_{\la+\om} |[\de]^2$ to
the triple $\al<\be<\de^-$ we get
$$  \al\be\le \al\de^- \vee \be\de^-,$$
so we are done also in this case since the quantity on the right-hand side is bounded by
$\al\ga\vee \be\ga$.

\noindent \emph{Subcase 2.2.} $\al\le \max M_0 <\be <\de$. Consider first the subcase when $\al\in
M_0$. To see that $\al\ga\le \al\be \vee \be\ga$ observe that $\al\ga\le p_M <p_M+1\le \be\ga$. To
see that $\al\be\le \al\ga \vee \be\ga$ observe that $\al\be$ appear as a quantity int the
definition of $\be\ga$. Let us consider the case $\al\notin M_0$ and let
$$  \al'=\min(M_0\setminus \al).$$
The quantity $\al\al'$ from the definition of $\al\ga$ is bounded by $\al\be\vee \be\ga$
since by 2. for $\ro_{\la+\om}|[\de]^2$ we have that
$$  \al\al'\le \al\be \vee \al'\be,$$
and $\al'\be$ appears in the definition of $\be\ga$. Since the quantities $\xi\al$ ($\xi
\in M_0\cap \al$ ) appear also in the definition of $\be\ga$, this establishes the
inequality $\al\ga\le \al\be \vee \be\ga$. Consider now the inequality $\al\be\le
\al\ga\vee \be\ga$. Apply 1. of $\ro_{\la+\om}|[\de]^2$ to $\al<\al'<\be$ and get
$$  \al\be\le \al\al'\vee \al'\be,$$
and this finishes the proof since $\al\al'\le \al \ga$ and $\al'\be\le \be\ga$.

\noindent \emph{Subcase 2.3.}  $\al<\be\max M_0$. If $\al,\be\in M_0$, then the inequalities 1.
and 2. for $\al<\be<\ga$ follow from the fact that in the definitions of $\al\be$, $\be\ga$ and
$\al\ga$ we copied the $\ro$-model $(M,<,\ro_M,p_M)$.

\noindent \emph{Subcase 2.3.1} $\al\in M_0$ and $\be\notin M_0$. Consider the inequality
$\al\ga\le \al\be\vee \be\ga$. This follows from the fact that
$$  \al\ga\le p_M<\be\be', \, \text{where } \be'=\min(M_0\setminus\be)$$
and the fact that in the definition of $\be\ga$ the quantity  $\be\be'$ appears. The inequality
$\rod$ in this subcase follows from the fact that the quantity $\al\be$ appears in the definition
of $\be\ga$.

\noindent \emph{Subcase 2.3.2} $\al\notin M_0$ and $\be\in M_0$. Consider the inequality
$\rou$. Let $$  \al'=\min (M_0\setminus \al).$$ We need to bound the quantities $\al\al'$
and $\xi\al$ ($\xi \in M_0\cap \al$) by $\al\be\vee \be\ga$. Apply 2. of
$\ro_{\la+\om}|[\de]^2$ to $\al<\al'\le \be$ and get
$$  \al\al'\le \al\be\vee \al'\be.$$
Since $\al'\be\le p_M$ and $\al\al'>p_M$ we conclude that $\al\al'\le \al\be$ as
required. Similarly note that
$$  \xi\al\le \xi\be\vee \al\be=\al\be,$$
since $\xi\be\le p_M$ while $\al\be>p_M$. It remains to check the inequality $\rod$ in this
subcase. As before note that
$$  \al\be\le \al\al' \vee \al'\be,$$
and that $\al\be>\al'\be$ since $\al'\be\le p_M$ while $\al\be>p_M$. It follows that
$\al\be\le \al\al'\le \al\ga$ as required.

\noindent \emph{Subcase 2.3.3} $\al\be\notin M_0$ (and $\al,\be\le \max M_0$). So in this subcase
both quantities $\al\ga$ and $\be\ga$ are defined according to the second definition. Let
$$  \al'=\min(M_0\setminus \al) \text{ and } \be'=\min(M_0\setminus \be). $$
Note that $\al'\le \be'$. We first check the inequality $\al\ga\le \al\be\vee \be\ga$. Consider
first the quantity $\al\al'$ that appears in the definition of $\al\ga$. If $\al'=\be'$, then
$$  \al\al'\le \al\be\vee \be\be'\le \al\be\vee \be\ga,$$
as $\be\be'$ appears in the definition of $\be\ga$. Suppose that $\al'<\be'$ i.e., that
$\al'<\be$. Then
$$  \al\al'\le \al\be\vee\al'\le \al\be\vee \be\ga$$
as $\al'\be$ appears as a quantity in the definition of $\be\ga$. Consider  the quantity
$\xi\al$ for $\xi\in M_0\cap \al$. Note that
$$  \xi\al\le \al\be\vee \xi\be\le \al\be\vee \be\ga,$$
as $\xi\be$ appears in the definition of $\be\ga$. It remains to check the inequality $\al\be\le
\al\ga\vee\be\ga$ in this subcase. If $\al'=\be'$, then we get that
$$  \al\be\le \al\be'\vee \be\be'\le \al\ga\vee\be\ga, $$
as the quantity $\al\be'=\al\al'$ appears in $\al\ga$ while $\be\be'$ appears in $\be\ga$. If
$\al'<\be'$ i.e., $\al'<\be\le \be'$, then we get that
$$  \al\be\le  \al\al' \vee \al'\be\le \al \ga \vee \be\ga,$$
since $\al\al'$ appears in $\al\ga$ and $\al'\be$ appears in $\be\ga$. This finishes checking
that the extension $\ro_{\la+\om}|[\de+l]^2$  satisfies the conditions 1. 2. and 3. Note that
$$  \ro(\al,\ga)>p_M \text{ for all } \al\in \de\setminus M_0 \text{ and } \ga\in [\de,\de+l),$$
we conclude that the set $M_0\cup [\de,\de+l)$ is $\ro_{\la+\om}$-closed. So the
extension $\ro_{\la+\om}|[\de,\de+l]^2$ has a set
$$  M_1=M_0\cup [\de,\de+1)\con \de+l,$$
which is $\ro_{\la+\om}$-closed while the corresponding model
$(M_1,<,\ro_{\la+\om}|[M_1]^2,p_M)$ is isomorphic to the given $\ro$-model
$(M,<,\ro_M,p_M)$. This finishes the recursive construction of a universal $\ro$-function
$\ro:[\ou]^2\to \om$. The reader is referred to \cite{tod1} for more on $\ro$-functions
and their uses. Some of the applications need the following {\it unboundedness property},
stronger than 3.
\begin{enumerate}
\item[]3'. For every $n<\om$ and infinite $M\con \ou$ there exist $\al<\be$ in $M$ such
that $\ro(\al,\be)>n$.
\end{enumerate}
As there is no reason to suspect that the universal $\ro:[\ou]^2\to \om$ just produced
satisfies 3'. we offer the following derived function $\bar{\ro}:[\ou]^2\to \om$:
$$  \bar{\ro}(\al,\be)=\max\{\ro(\al,\be), \#(\conj{\xi\le\al}{\ro(\xi,\al)\le \ro(\al,\be)})\}.$$
It may be checked that $\bar{\ro}$ has the properties 1., 2. and 3'.
\subsection{A smooth $\ro$-function}
We construct a $\ro$-function such that the  corresponding coding $\sig_\ro$ yields that
$\eqs_\xi$ has a Schauder basis  for every $\xi<\ou$. These bases will be a reordering of the
transfinite basis $(e_\al)_{\al<\xi}$ in order type $\om$.

For a given $\ro$-function and an ordinal  $\al<\ou$, let
$F_n^\al=\conj{\be<\al}{\ro(\be,\al)\le n}$, which are $n$-closed.
\defi
A $\ro$-function is called \emph{smooth} if for every limit ordinal $\la<\ou$, the
numerical sequence $(\# F_n^\la/n)_n$ is bounded. \fdefi

\prop
There exists a smooth $\ro$-function.
\fprop
\prue Let us show that such smooth $\ro$-function exists. The definition will be again
inductive,  i.e., for each  limit ordinal $\la$ we are going to define $\ro_{\la}:[\la]^2\to
\om$. Suppose we have defined $\ro_\la$, and set
\begin{equation}
\ro_{\la+\om}(\al,\la)=\max\{{g_{\la}(\# (C_\la\cap \al))}, \ro_\la(\al,\min (C_\la\setminus \al)) ,
\ro_\la(\xi,\al) \, : \, \xi\in C_\la\cap \al \}
\end{equation}
where  $g_\la:\N\to \N$ is increasing and $C_\la$ is cofinal in $\la$ and they are defined as  follows: For a
given $\al<\la$, let $i(\al)=\# C_\la \cap \al$. Suppose we have constructed $\ro_\la$ such that for all
limit $\ga<\la$ with the smooth property $ \lim_{n\to \infty}\# F_n^\ga/n=0 $. There are two cases.

\noindent (a) Suppose that $\la=\ga+\om$ is a successor limit.
 For each integer $i$, let
$g_{\la}(i)=  2^{i} \text{ and } C_\la =  \{\ga+n\}_n$. Notice that for $\al<\la$, if
$\ro(\al,\la)\le n$, then either $\al<\ga$ and $\al\in F^\ga_{n}$ or if $\al=\ga+l$, then
since $i(\al)=\# C_\la\cap \al= l+1$, we have that $l\le \log_2(n)$.  So, $\# F^\la_n\le \#
F^\ga_n+ 1+\log_2 n$, which certainly implies that $\# F^\la_n /n\to_n 0$.

\noindent (b) Suppose that $\la$ is a limit of limit ordinals.   Let $C_\la=\{\la_n\}_n\con
\la$ cofinal in $\la$, with  each $\la_n$ a limit ordinal.  Let $(n_i)_i$ be a strictly
increasing sequence of integers such that
\begin{equation}
\# F^{\la_0}_n/n +\dots + {\# F^{\la_i}_n}/n\le 2^{-i}
\end{equation}
for every $i$ and every $n\ge n_k$. Let $g_\la(i)=n_i$ for all $i$. Fix $\vep>0$, and let $j$
be such that $1/2^{j}\le \vep$. We show that for all $n\ge n_j$, $\# F^\la_n/n \le \vep$: Fix
$n\ge n_j$, and let $i_0$ be the maximal integer such that $n_j\le n_{i_0}\le n$. Notice that
then
\begin{equation}
F^\la_n\con \conj{\al<\la}{i(\al)\le i_0 \text{ and } \al \in F^{\la_{i(\al)} } _n }=F^{\la_0}_n\cup \dots
\cup F^{\la_{i_0} }_n.
\end{equation}
So, ${\# F^{\la}}/n \le {\# F^{\la_0}_n}/n+\dots + {\# F^{\la_{i_0}}_n}/n \le {2^{-i_0}}\le {2^{-j}}\le
\vep$.
\fprue

\lema \label{decomp1} Let $G$ be a $p$-closed set, and let
$\phi=(1/w(\Phi))\sum_{i=1}^d\phi_i\in K$ be of odd weight $w(\phi)<p$, and such that
$\supp\phi\cap G\neq \buit$. If defined, set
\begin{align*}
d_1= &\max\conj{i\in [1,d]}{w(\Phi_{i})<p} \text{ and } d_0=  \max\conj{i<d_1}{\supp \Phi_i\cap G\neq \buit},
\end{align*}
Then

\noindent(1)  $(1/w(\phi))\sum_{i=1}^{d_0}\Phi_i|[0,\al]$ has support contained in $G$, where
$\al=\max (\supp \Phi_{d_0}\cap G)$.

\noindent(2) $\supp \phi_i \cap G=\buit$ for every $d_0<i <d_1$.

\noindent (3)  $w(\phi_i)\ge p$ for every $d_1<i\le d$.

\flema
\prue
If $d_1$ is not well defined, then for all $k\le d$, $w(\phi_k)\ge p$. If $d_1$ is well defined, but $d_0$ is
not, then for all $k<d_1$ we have that $\supp \phi_k \cap G=\buit$. Suppose that both are well defined.

Finally, since $p\ge w(\phi_{d_1})\ge \max\{ p_\ro(\bigcup_{i=1}^{d_0}\supp \phi_i),w(\phi_i)\}$
and $\cl{\bigcup_{i=1}^{d_0}\supp \phi_i \cap (\al+1)}{p}\con G$ ($G$ is $p$-closed), it follows
that the support of $(1/w(\phi))\sum_{i=1}^{d_0}\Phi_i|[0,\al]$ is included in $G$.\fprue

\lema\label{f0f1}
 Let  $G\con \ou$  be  $p$-closed. Then for all $\phi\in K_{\ou}$  there are
 some $f_0$ and $f_1$ such that
\begin{enumerate}
 \item $\supp f_0,\supp f_1\con G$, $f_0+f_1=\phi | G$
 \item $\nrm{f_0}_\infty \le 1/p$,
 \item $f_1\in 2 K_{\ou}(G)$, where $K_{\ou}(G)$ is the subset of $K_{\ou}$ consisting on the functionals
 $\phi$ with support contained in $G$.
\end{enumerate}
\flema
\prue
Let $(\phi_t)_{t\in \mc T}$ be a tree-analysis of $\phi$. Let
\begin{align*}
\mc T_0= &\{t\in \mc T\, : \, \text{ there is some } u \preceq t \text{ with } w(\phi_u)\ge p\} \text{ and }
\mc T_1 = \mc T \setminus \mc T_0.
\end{align*}
Notice that $\mc T_1$ is a downwards closed subtree of $\mc T$, and hence for a given $t\in \mc T_1$, the set
$S_t^{1}$ of immediate successor of $t$ in $\mc T_1$ is exactly equal to $S_t^1 =S_t\cap \mc T_1.$ If $\mc
T_1=\buit$, then $\phi_0=\phi$ has to be of type I and $w(\phi_0)\ge p$. In this case, let $f_0=\phi$ and
$f_1=0$, that clearly satisfies what we want. Suppose now that $\mc T_1\neq \buit$. We are going to find for
all $t\in \mc T_1$, $f_0^t,f_1^t$ such that

\begin{enumerate}
 \item $\supp f_0^t,\supp f_1^t\con G$, $f_0^t+f_1^t=\phi | G$
 \item $\nrm{f_0^t}_\infty \le 1/p$,
 \item $f_1^t\in 2 K_{\ou}$.
\end{enumerate}It is clear that the pair $f_0=f_0^0$, $f_1=f_1^0$ satisfies our requirements.  The proof goes by downwards
induction over $t\in \mc T_1$ on the tree $\mc T_1$. Suppose that $t\in \mc T_1$ is a
terminal node of $\mc T_1$.

\noindent \emph{(1)} If $t$ is terminal node of $\mc T$, then we set $f_0^t=0$ and $f_1^t=\phi_t$
if $\supp \phi_t\con G$, and $f_0^t=f_1^t=0$ otherwise.

\noindent{\emph{(2)}} If $t$ is not terminal node of $\mc T$, then this means that for all $s\in
S_t$, $\phi_s$ is of type I and $w(\phi_s)\ge p$. Set $f_0^t=\phi_t|G$, $f_1^t=0$.

Suppose now that $t\in \mc T_1$ is not terminal in $\mc T_1$. Clearly this implies that $t$ is not
terminal in $\mc T$. There are three cases: \emph{Case 1}. $\phi_t$ is of type II,
$\phi_t=\sum_{s\in S_t } r_s \phi_s$. Then we set
\begin{equation*}
  f_0^t=  \sum_{s\in S_t^1} r_s f_0^s+\sum_{s\in S_1\setminus \mc T_1} r_s\phi_s|G \text{
and } f_1^t = \sum_{s\in S_t^1}r_s f_1^s.
\end{equation*}
Since for $s\in S_t$, $s\notin \mc T_1 $ iff  $\phi_s$ is of type I and  $w(\phi_s)\ge p$, this gives that
$\nrm{f_0^t}_\infty \le 1/p$. The rest our inductive promises for $f_0^t$ and $f_1^t$ are clearly satisfied.

\noindent \emph{Case 2}. $\phi_t$ is of type I, and $w(\phi_t)$ is even. We set
\begin{equation*}
  f_0^t=  \frac1{w(\phi_t)}\sum_{s\in S_t^1}  f_0^s+ \frac1{w(\phi_t)} \sum_{s\in
S_t\setminus \mc T_1} \phi_s|G \text{ and } f_1^t =  \frac1{w(\phi_t)}\sum_{s\in S_t^1}  f_1^s.
\end{equation*}
The condition $\nrm{f_0^t}_\infty \le 1/p$ is satisfied by the same reason as in the previous
case.

\noindent \emph{Case 3}. $\phi_t$ is of type I, and $w(\phi_t)$ is odd,
$\phi_t=(1/w(\phi_t))\sum_{i=1}^d \phi_{s_i}$, where $\{s_1,\dots,s_d\}=S_t$. Find $d_0<d_1\le
d$ as in the previous Lemma \ref{decomp1}. If $d_1$ is not well defined, this implies that
$w(\phi_s)\ge p$ for every $s\in S_t$. Then $S_t^1=\buit$ and we set $f_0^t=  \phi_t| G \text{
and } f_1^t = 0$. Suppose that $d_1$ is well defined but $d_0$ is not. This means that $\supp
\phi_k\cap G=\buit$ for every $k<d_1$. Then we set
\begin{equation*}
f_0^t=  \frac1{w(\phi_t)}\left( f_0^{s_{d_1}}+\sum_{i=d_1+1}^d \phi_{s_i}|G\right) \text{ and }
f_1^t = \frac1{w(\phi_t)} f_1^{s_{\la_{\phi,\phi'}}}.
\end{equation*}
Suppose now that both $d_0$ and $d_1$ are well defined, then we set
\begin{equation*}
f_0^t=  \frac1{w(\phi_t)}\left( f_0^{s_{d_1}}+\sum_{i=d_1+1}^d \phi_{s_i}|G\right) \text{ and }
f_1^t = \frac1{w(\phi_t)}\left(\sum_{i=1}^{d_0} \phi_i|[0,\al]\right)+ \frac1{w(\phi_t)}
f_1^{s_{\la_{\phi,\phi'}}},
\end{equation*}
where $\al=\max (\supp \phi_{d_0}\cap G)$. Notice that $1/{w(\phi_t)}(\sum_{i=1}^{d_0}
\phi_i|[0,\al])\in K_{\ou}(G)$. Therefore, using the induction hypothesis we conclude that
$f_1^t\in 2 K_{\ou}$.\fprue

\lema Assume that $\eqs_{\ou}$ is built upon a smooth
$\ro$-function and fix a limit ordinal $\la<\ou$. Then the
projections $(\P_{F^\la_n})_n$ are uniformly bounded by $2+D_\la$,
where $D_\la=\sup_n \# F^\la_n/n<\infty$.
\flema
\prue Fix a limit ordinal $\la$, let $x\in \eqs_\la$ be of norm 1, and  $\phi\in K_{\ou}$.Take
the decomposition $\phi=f_0+f_1$  from the previous Lemma \ref{f0f1} applied to the $n$-closed
set $F^\la_n$. Then,
\begin{equation}
|\phi \P_{F_n^\la}x|=|\langle f_0 ,  \P_{F_n^\la}x \rangle +\langle f_1 , \P_{F_n^\la}x \rangle|\le \frac{\#
F_n^\la}n+|\langle f_1 , \P_{F_n^\la}x \rangle|\le D_\la+ |\langle f_1 , \P_{F_n^\la}x \rangle|.
\end{equation}
Now using that  $f_1\in 2K_{\ou}({F_n^\la})$, we can write $f_1=\sum_i \la_i \phi_i$,
$\sum_i \la_i\le 2$, $\la_i\ge 0$, and $\phi_i\in K_{\ou}({F_n^\la})$. Therefore,
$\langle \phi_i, \P_{F_n^\la}x\rangle = \langle \phi_i, x\rangle\le 1$. So, $|\langle f_1
, \P_{F_n^\la}x \rangle|\le 2$, and we are done. \fprue

Let $Q_n^\la=\P_{F_n^\la}$.  Notice that $Q_n^\la Q_m^\la= Q_{\min\{n,m\}}^\la$.

\teor
For every $x\in \eqs_\la$, $\lim_{n\to \infty} Q_n^\la (x)=x$.
\fteor
\prue Let us show that for all limit $\be\le \la$ and all $x\in \eqs_\be$, $\lim_{n\to \infty}
\P_{F_n^\la}x=x$. The proof is by the induction over the set of limit ordinals $\le\la$. Fix
$x\in \eqs_\be$.

\noindent (a) $\be=\om$.  We know that $\lim_{n\to \infty}\P_n x=x$. Fix $\vep>0$, and let
$n_0$ be such that for all $n\ge n_0$, $\nrm{x- \P_{n}x}\le \vep/(3+D_\la)$. Let $n_1\ge n_0$
be such that for all $n\ge n_1$, $[0,n_{0}]\con F_{n}^\la$. Hence
$\nrm{x-\P_{F_n^\la}x}\le\nrm{x- \P_{n_0}{x}} +\nrm{ \P_{F_n^\la}(x- \P_{n_0}{x})}\le
(1+\nrm{\P_{F_n^\la}})\nrm{x- \P_{n_0}{x}}\le  \vep$ for every $n\ge n_1$.

\noindent (b) $\be=\ga+\om$. Then, $x=y+z$, where $y\in \eqs_\ga$, and $z\in
\eqs_{[\ga,\ga+\om)}$. By the induction hypothesis, $\lim_{n\to \infty} \P_{F_n^\la}y=y$.  Now
use the projections $(\P_{[\ga,\ga+n)})_n$ to approximate $z$ and follow the ideas of the case
$\be=\om$.

\noindent (c) $\be$ is limit of limit ordinals. Fix a strictly increasing  sequence
$(\be_n)_n$ with limit $\be$, and let $x_n=\P_{\be_n}x$. We know that $\lim_{n\to \infty}x_n=
x$. Fix $\vep>0$, and let $n_0$ be such that $\nrm{x-x_{n_0}}\le \vep/2(3+D_\la)$. Let $n_1\ge
n_0$ be such that $\nrm{x_{n_0}-\P_{F_n^\la}x_{n_0}}\le \vep/2$ for all $n\ge n_1$, that we
know that it is possible by  the induction hypothesis since $x_{n_0}\in \eqs_{\be_{n_0}}$.
Then
 for all $n\ge n_1$,
\begin{align*}
\nrm{x-\P_{F_n^\la}x_{n}}\le &  \nrm{x-x_{n_0}}+
 \nrm{x_{n_0}-\P_{F_n^\la}x_{n_0}}+
 \nrm{\P_{F_n^\la}}\nrm{x-x_{n_0}}\le  \\
\le & (3+D_\la) \nrm{x-x_{n_0}}+\nrm{x_{n_0}-\P_{F_n^\la}x_{n_0}}\le \vep.
\end{align*}
\fprue

\cor The space $\eqs_\al$ has a Schauder basis for every ordinal $\al<\ou$. Moreover, for
every $\al<\ou$ there exists a reordering  $(e_{\be_n})_n$ of $(e_\be)_{\be<\al}$ such
that $(e_{\be_n})_n$ is a Schauder basis of the space $\eqs_{\al}$.
 \fcor
  \prue It is
enough to show the result for a limit ordinal $\la$. By the previous Theorem, the projections
$(Q^\la_n)_n$ define a finite dimensional Schauder decomposition of $\eqs_\la$. Notice
 that the natural ordering $<_\la$ on $\la$   defined by
$$  \al<_\la \be \text{ iff }\left\{\begin{array}{ll}
\ro(\al,\la)<\ro(\be,\la) & \text{ or} \\
\ro(\al,\la)=\ro(\be,\la) & \text{ and } \al<\be
\end{array}\right.$$
 has order type $\om$.
Let  $\{\la_n\}_n$ be an enumeration of $(\la,<_\la)$ in order type $\om$,   and let us show that
$(x_n=e_{\la_n})_n$ is a basis of $\eqs_\la$:  Let $(R_n)_n$  be the projections $R_n:\eqs_\la\to
\eqs_\la$ associated to $(x_n)_n$.  For a given $k$, let $n_k=\ro(\la_k,\la) $. Then,
$R_k=Q^\la_{n_k-1}+ \P_{\la_k}\circ (Q_{n_k}^\la -Q_{n_k-1}^\la)$. This clearly shows that
$(x_n)_n$ is a Schauder basis of $\eqs_\la$.\fprue

\nota It is unclear whether there is a variation on $\ro$ such that some of the resulting
spaces $\eqs_\la$ $(\la<\omega_1)$ do not admit Schauder basis. \fnota

\section{Universality of  $\ro$  and nearly subsymmetric bases}

Throughout  this section we assume that the coding $\sig_\ro$ is based on an universal
$\ro$ function discussed  in the previous section.

\nota For the sequel we need a slight modification of the definition of special sequences. More
precisely, we assume that for each
$(\phi_1,w_1,p_1,\dots,\phi_{n_{2j+1}},w_{n_{2j+1}},p_{n_{2j+1}})$ every $p_i$ satisfies the
additional property that  for all $l\le i$,  $\phi_l$  admits a tree-analysis with supports in the
set
\begin{align*}
G_i=\cl{\bigcup_{k=1}^{l}\supp \phi_k}{p_i}.
\end{align*}
Note that the definition of the special functionals and the fact that $K_{\ou}$ is rationally
closed does not allow one to conclude that every functional $\phi\in K_{\ou}$ admits a
tree-analysis $(\phi_t)_{t\in \mc T}$ such that for every $t\in \mc T$, $\supp \phi_t\con
\supp \phi$. However there will always be large enough $p$ such that $\cl{\supp \phi}p$
contains a tree-analysis of $\phi$. This follows from the fact that there is a tree-analysis
$(\phi_t)_{t\in \mc T}$  of $\phi$ such that $\max\bigcup_{t\in \mc T}\supp \phi_t=\max \phi$.
 \fnota

\defi
For a given $p$, and a subset $G\con \ou$
 let
\begin{align*}
K_p(G)=  \{ & \phi\in K \, :\,\phi \text{ has some tree-analysis }(\phi_t)_{t\in \mc T}
\text{ such that } \ptot t\in \mc T \supp \phi_t\con G \text{ and}\\
& \text{if } \phi_t \text{ has type I, then } w(\phi_t)<p \}.
\end{align*}
We will call such  tree-analysis $\mc F=(\phi_t)_{t\in \mc T}$ of $\phi\in K_p(G)$ a
\emph{$(p,G)$-tree-analysis} of $\phi$. Notice that if $\mc F=$ is a $(p,G)$-tree-analysis of
$\phi$, then for all interval $E$, $(\phi_t|E)_{t\in \mc T}$ is a $(p,G)$-tree-analysis of
$\phi|E$.
\fdefi

\prop\label{moving}
Suppose that $G,G'\con \ou$ are $p$-complete, and ($\ro$-)isomorphic (see Definition
\ref{rhoisomophism}). Then the unique order preserving mapping $\pi:G\to G'$ defines a bijection
$$  \tilde{\pi}:K_p(G)\to K_p(G')$$ such that for every $\al\in G$ $\tilde{\pi}(e_\al^*)=e_\be^*$,
preserves $(p,G)$-tree-analysis in $K_p(G)$ and weights.
\fprop
\prue
The proof is an easy use of downwards induction over a $(p,G)$-tree-analysis.
\fprue

Using the properties of our new coding $\sig_\ro$ we can improve Lemma \ref{f0f1} as follows.
\lema\label{f0f11}
 Let  $G\con \ou$  be  $p$-closed. Then for all $\phi\in K_{\ou}$  there are
 some $f_0$ and $f_1$ such that
\begin{enumerate}
 \item $\supp f_0,\supp f_1\con G$, $f_0+f_1=\phi | G$,
 \item $\nrm{f_0}_\infty \le 1/p$, and
 \item $f_1\in 2 K_{p}(G)$.
\end{enumerate}
\flema
\prue The Proof follows exactly the same steps than the proof of Lemma \ref{f0f1} with the
exception that the inductive premise 3. $\supp f_1^t\in 2 K_{\ou}(G)$ now is replaced by
$f_1^t\in 2 K_p(G)$. To check that one can find the corresponding decomposition when one deals
with the case of odd weight, we notice that the premise 3. will be fulfilled since the new
coding $\sig_\ro$ will guarantee that in Lemma \ref{decomp1}, the corresponding
$(1/w(\phi))\sum_{i=1}^{d_0}\phi_i|[0,\al]\in K_p(G)$. \fprue

\defi
A transfinite basis $(e_{\al})_{\al<\ga}$ is said to
 be $C$-\emph{{nearly subsymmetric}} if for every $\vep>0$ and for every
 family of finite successive subsets  $\{F_i\}_{i=1}^d$ of
 $\ga$ and every  family  $\{I_i\}_{i=1}^d$ of successive
 infinite intervals there exists $\{G_i\}_{i=1}^d$ with
 $G_i\subseteq I_i$, $\#G_i=\#F_i$ such that the natural isomorphism
 $T:\langle(e_{\alpha})_{\alpha\in \cup_{i=1}^dF_i}\rangle\to
 \langle(e_\beta)_{\beta \in \cup_{i=1}^dG_i}\rangle$
 satisfies $\|T\| \cdot \|T^{-1}\|\le C+\vep$.
\fdefi

The purpose of this section is to prove the following result.
\teor\label{movgen} The transfinite basis $(e_\al)_{\al<\ou}$ is $4$-nearly subsymmetric.

\fteor

\prue
We want to show  that for every sequence of finite sets  $F_1<F_2<\dots <F_n$, infinite
intervals $I_1\le I_2\le \dots \le I_n$ (with possible repetitions) and $\vep>0$, there
is some $G_1<G_2<\dots <G_n$ such that

\noindent \emph{(a)} $G_i\con I_i $, $i=1,\dots,n$,

\noindent\emph{(b)} $\# F_i=\# G_i$, $i=1,\dots,n$,

\noindent \emph{(c)} The natural isomorphism $T$ between $\eqs_{F}$ and $\eqs_G$ satisfies
that $\nrm{T},\nrm{T^{-1}}\le 2+\vep$, where $\eqs_F=\langle e_\al \rangle_{\al \in F}$ and
$\eqs_G=\langle e_\al \rangle_{\al \in G}$ and $T$ is defined for $\al\in F$ and $\be \in G$
such that  $T(e_{\al})=e_\be$ satisfies that  $\al\mapsto \be$ is order preserving and onto,
and $F=\bigcup_{i=1}^n F_i$ and $G=\bigcup_{i=1}^n G_i$.

Let $p\ge \max\{p_\ro(\bigcup_{i=1}^n F_i),  \# F/\vep\}$, and let $\tilde{F}=\cl{F}p$. For each
$i=1,\dots,n$, let $\al_i=\max F_i+1$, and let $F_i'=G\cap \al_i$. Notice that

\noindent (1) $F_i\con F_i'$ for every $i=1,\dots,n$,

\noindent  (2) each  $F_i'$ is a $p$-closed set for $i=1,\dots,n$, and

\noindent (3)  $F_i'$ is an initial segment of $F_j'$ for $i\le j\le n$.

By the universality of  $\ro$, there is some $G_1'\con I_1$  which is $\ro$-isomorphic to
$F_1'$. Since $(F_2',<,\ro|[F_2']^2,p)$ and  $(F_1,',<,\ro|[F_1']^2,p)$ are  $\ro$-models,
$F_1'$ is an initial segment of $F_2'$, and $(F_1',<,\ro|[F_1']^2,p)\cong
(G_1',<,\ro|[G_1']^2,p)$, the universality of $\ro$ gives a set $H_2\con I_2\setminus G_1'$
such that $G_2'=G_1'\cup H_2$ satisfies that

\noindent (1) $G_2'$ is $p$-closed, and

\noindent (2)$(G_2',<,\ro|[G_2']^2,p)\cong (F_2',<,\ro|[F_2']^2,p)$.

And so on. At the end we get $n$ many $\ro$ models $(G_i',<,\ro|[G_i']^2,p)$ for $i=1,\dots,n$
such that

\noindent (1) For $i\le j\le n$,  $G_i'$ is an initial segment of $G_j'$,

\noindent (2) $G_{i}'\setminus G_{i-1}'\con I_i$, for $i=2,\dots,n$, and $G_1'\con I_1$,

\noindent (3) $(G_i',<,\ro|[G_i']^2,p)\cong (F_i',<,\ro|[F_i']^2,p)$ for $i=1,\dots,n$.

Therefore $\tilde{G}=G_n'$ satisfies that $(\tilde{G},<,\ro|[G']^2,p)\cong (\tilde{F},<,\ro|[F']^2,p)$. Let
$\pi$ be the isomorphism between them, and for each $i=1,\dots,n$, let $G_i =\pi F_i$.  Let us show that
$T:\eqs_F\to \eqs_G$ satisfies what we wanted:   Fix one vector $x=\sum_{i\in F} \la_i e_i$ such that
$\nrm{x}=1$. Take $\phi\in K$ such that $\phi x =1$.  Then take the decomposition of  $\phi=f_0+f_1$ as in
previous Lemma \ref{f0f11}. Using Proposition \ref{moving},   we can take a copy $g_1$ of $f_1$ in $2
K_{p}(\tilde{G})$. Since $1=\phi x =|f_0x+f_1x|\le |f_0 x|+ N/p<|g_1Tx| +\vep $, $|g_1 Tx|>1-\vep$. This
implies that there is some $\psi\in K_p({\tilde{G}})$ such that $|{\psi}Tx|> (1-\vep)/2$.
 So, $\nrm{Tx}\ge (1-\vep)/2$.

Now suppose that $\nrm{Tx}> 2+\vep$. Then, let $\phi\in K$ be such that $\phi Tx
>2+\vep$. Take the decomposition $\phi=g_0+g_1$ as in the previous Lemma \ref{f0f11}, now in
$K_{p}(\tilde{G})$. This implies that $g_1 Tx >2$, and hence, there is some ${\psi}\in K_{p}(\tilde{G})$ such
that ${\psi} Tx
>1$. Hence, the copy $\phi$ of ${\psi}$ in $K_{\tilde{F}}(p)$ is such that $\phi x = {\psi}
Tx>1$, a  contradiction. So, $\nrm{T}\le 2+\vep$ and $\nrm{T^{-1}}\le  2/(1-\vep)\le  2+\vep$.
\fprue

\defi
Recall the following (modified) notion from \cite{Ma-Mi-To}. Let $X$ be a Banach space with a
Schauder basis $(u_n)_n$, fix $n\in \N$ and $C\ge 1$. A finite $n$-dimensional space $E$ with
a basis $(e_i)_{i=1}^n$ is called a \emph{$C$-asymptotic space} of $X$ iff
\begin{equation}
\sup_{X_1} \inf_{x_1\in S(X_1)} \sup_{X_2}\ldots \inf_{x_n\in S(X_n)} d_b(\langle
x_1,\ldots,x_n\rangle, E) \le C,
\end{equation}
where $d_b=\nrm{T}\cdot \nrm{T^{-1}}$ for $T:\langle x_1,\ldots,x_n\rangle\to E$ to be the natural
isomorphism defined by $T(x_i)=e_i$, and $X_i$ runs over all tail subspaces, i.e,
$X_i=\overline{\langle u_i \rangle_{i>k}}$ for some $k$. A space $Y$ with a monotone basis
$(y_n)_n$ is called a \emph{$C$-asymptotic version} of $X$ iff for every $n$, $\langle y_i
\rangle_{i=1}^n$ is an asymptotic space of $X$.
\fdefi

\cor
There exists a family $\{X_{\ga}\}_{\ga<\ou}$ of reflexive totally incomparable hereditarily
indecomposable spaces with Schauder bases such that $X_\ga$ is an asymptotic version of
$X_{\ga'}$ for every $\ga,\ga'<\ou$.\qed\fcor

\defi
Two transfinite basis $(x_\al)_{\al<\ga_0}$ of $X_0$ and $(x_\al)_{\al<\ga_1}$ of $X_1$ are called
\emph{finitely equivalent} iff there is some constant $C>0$ such that for all finite set $F_0\con \ga_0$
there is some finite set $F_1\con \ga_1$ with the same cardinality such that $(x_\al)_{\al\in F_0}$ and
$(y_\al)_{\al\in F_1}$ are $C$-equivalent.
\fdefi

\nota There are finitely equivalent subspaces of $\eqs_{\ou}$
which are incomparable. \fnota

\nota Using the fact that $\ro$ is universal  it can be shown that
for any bounded $T:\eqs_{\ou}\to \eqs_{\ou}$, $\nrm{D_T}\le 4\nrm{T}$. The proof goes as
follows. We assume that $\nrm{T}=1$. Fix a normalized finitely supported vector $x$; let
$x=x_1+\dots +x_n$ be its decomposition  in $\eqs_{\ou}$, and $\vep>0$. Let $F_i=\supp x_i$
for each $i=1,\dots,n$, and consider infinite intervals $I_1\le I_2\le \dots \le I_n$ of $\ou$
such that $\nrm{D_T(y)-T(y)}\le \vep\nrm{y}$ for every $y\in \eqs_{I_1\cup \dots \cup I_n}$.
By Theorem \ref{movgen} we can find for every $i=1,\dots,n$ $G_i\con I_i$ such that $\# G_i=\#
F_i$ and the order isomorphism between $F=\bigcup_{i=1}^n F_i$ and $G=\bigcup_{i=1}^n G_i$
defines an isomorphism  $H$ between $\langle e_\al \rangle_{\al\in F}$ and $\langle e_\al
\rangle_{\al\in G}$ with $\nrm{H},\nrm{H^{-1}}\le 2+\vep$. Set $y=F (x)$ and then  $\nrm{y}\le
2+\vep$ and  $ \nrm{(T-D_T)(y)\le \vep}$. Since $H(D_T(x))=D_T(y)$ we have that  $
\nrm{D_T(x)}\le ({2+\vep})\nrm{D_T(y)}$. So,
\begin{equation}
\nrm{D_T(x)}\le  ({2+\vep})\nrm{D_T(y)} \le  ({2+\vep})(\nrm{D_T(y)-T(y)}+\nrm{T y})\le ({2+\vep})\vep +
({2+\vep})^2 \nrm{T}.
\end{equation}
\fnota

\section{Tree-analysis of functionals: Basic inequality and Finite interval representability of $J_{T_0}$}\label{finblorep}
The goal of this section is to prove the basic inequality (Lemma \ref{bin}) and show the
finite interval representability of the James-like space $J_{T_0}$ in $\eqs_{\ou}$ (Theorem
\ref{represent}). Reaching these two goals involve similar sort of problems and for this
reason we introduce a  general theory applicable to both cases and hopefully to many other
cases to come.

\subsection{General theory}
The theory deals with a block sequence of vectors $(x_k)_{k=1}^n$, a sequence of scalars
$(b_k)_{k=1}^n$, and a functional $f\in K_{\ou}$, and tries to estimate $|f(\sum_{k=1}^n b_k
x_k)|$ in terms of $|g(\sum_{k=1}^n b_k e_k)|$ for an appropriately chosen functional $g$ of
an auxiliary Tsirelson-like space $X$ with basis $(e_i)_i$. The natural approach is to start
with a tree-analysis $(f_t)_{t\in \mc T}$ of $f$, and try to replace the functional $f_t$ at
each node $t\in \mc T$  by a functional $g_t$ in the norming set of the auxiliary space, and
in doing this try to copy, as much as possible,  the given tree-analysis $(f_t)_{t\in \mc T}$.
Not all nodes $t\in \mc T$ have the same importance in this process.  It turns out that the
crucial replacements $f_t \mapsto g_t$ are made for $t$ belonging to some sets $\mc A\con \mc
T$ such that $(f_t)_{t\in \mc A}$ is in some sense responsible for the estimation of the
action of the whole functional $f$ on each of the vectors $x_k$. These are  the \emph{maximal
antichains} of $\mc T$ defined below. Observe that some of the replacements $f_t\mapsto g_t$
are necessary before this procedure has a chance to work. Suppose for example the replacements
are  made in an auxiliary mixed Tsirelson space $X$ where a particular $(m_{j_0}^{-1},
{n_{j_0}})$-operation is not allowed. Then, every time we find a node $t\in \mc T$ such that
the corresponding $f_t$ has weight $w(f_t)=m_{j_0}$ the replacement $g_t$ has to be something
avoiding this operation, i.e., we cannot put  the combination $g_t=(1/w(f_t))\sum_{s\in S_t}
g_s$. These sorts of nodes are the ones that we call ``catchers" below, because their own
tree-analysis  $(f_s)_{s \succeq t}$ cannot be taken into account.

\subsubsection{Antichains and arrays of antichains}
Recall that every $f\in K_{\ou}$ has a tree-analysis $(f_t)_{t\in \mc T}$ such that: For every
$t\in \mc T$ (a) if $u \succeq t$, then $\ran f_u \con \ran f_t$, and (b) if $f_t$ is of type I,
then $f_t=(1/w(f_t))\sum_{s\in S_t} f_s$.

Recall that $A\con \mc T$ is called an \emph{antichain} if for
every $t\neq t'\in A$, neither $t \preceq t'$ nor $t'\preceq t$.
Given $t,t'\in \mc T$, we define $t\wedge t'=\max\conj{v\in \mc
T}{v \preceq t,t'}$. Notice that $\mc A\con \mc T$ is an antichain
iff $t\wedge t' \precneqq t,t'$ for every $t\neq t'\in \mc T$.

\defi
Fix a tree-analysis $(f_t)_{t\in\mc T }$ of $f$ as above.  Given a finitely supported vector $x$,
a set $\mc A\con \mc T$ is called a \emph{regular antichain} for $x$ and $(f_t)_{t\in \mc T}$ if

\noindent \emph{(a.1)} for every $t\in \mc A$, $f_t$ is not of  type II,

\noindent \emph{(a.2)}  $f_{t_1\wedge t_2}$ is of type II for every $t_1\neq t_2\in \mc
A$, and

\noindent \emph{(a.3)}  $\ran f_t\cap \ran x\neq \buit$,  for every $t\in \mc A$.

$\mc A$ is a\emph{ maximal antichain} for $x$  if in addition $\mc A$ satisfies

\noindent{\emph{(a.4)}}  for every $t\in \mc T$ if $\supp f_t\cap \ran x\neq \buit$, then there is
some $u\in \mc A$ comparable with $t$.

Let $(x_k)_{k=1}^n$ be a block sequence, and let $\boldsymbol{\mc A}=(\mc A_k)_{k=1}^n$  be such
that each $\mc A_k$ is a regular antichain for the vector $x_k$ and the tree-analysis $(f_t)_{t\in
\mc T}$. For a given $t\in \mc T$, we define
$$  D_t^{\boldsymbol{\mc A}}=  \bigcup_{u \succeq t}\conj{k\in[1,n]}{u\in \mc A_k}, \,  E_t^{\boldsymbol{\mc A}}= D_t^{\boldsymbol{\mc A}}\setminus
(\bigcup_{s\in S_t} D_s^{\boldsymbol{\mc A}}).$$
Whenever there is no possible confusion we simply
write $D_t$ and $E_t$ to denote $D_t^{\boldsymbol{\mc A}}$ and $E_t^{\boldsymbol{\mc A}}$
respectively.

$\boldsymbol{\mc A}=(\mc A_k)_{k=1}^n$ is called a \emph{(maximal) regular array} for
$(x_k)_{k=1}^n$ and $(f_t)_{t\in \mc T}$ if each $\mc A_k$ is a (maximal) regular antichain for
$x_k$ and $(f_t)_{t\in \mc T}$, and in addition

\noindent \emph{(a.5)} for every $t\in \bigcup_{k}\mc A_k$ such that $f_t$ is of type I, either
$t$ is a \emph{catcher}, i.e.,
  $D_s=\buit$ for
every $s\in S_t$, or for every $k\in E_t$, $t$ is a \emph{splitter} of $x_k$, i.e.,  for every
$k\in E_t$ there are at least $s_1\neq s_2\in \mc S_t$ such that $\ran f_{s_i}\cap \ran
x_k\neq \buit$.

We denote by $\mm S({\boldsymbol{\mc A}})$ and $\mm
C({\boldsymbol{\mc A}})$ the set of splitter nodes and catcher
nodes of $\boldsymbol{\mc A}$, respectively. Notice that if
$t_i\in \mc A_{k_i}$ ($i=1,2$) are catcher nodes, then they are
incomparable, and that $\mc A_k=\mm S({\boldsymbol{\mc A}})\cup
\mm C({\boldsymbol{\mc A}})$.

Note that if no $f_t$ ($t\in \mc T$) is of type II then $\#\mc
A_k\le 1$ for all $k$, and so the tree-analysis below becomes much
simpler.\fdefi
\defi\emph{(The functor $\mc A(x,\mm C)$.)}
Given  a block vector $x$ and $\mm C\con \mc T$ consisting of nodes of type I, let $\mc A(x, \mm
C)$ be the set of nodes $t\in \mc T$ such that

\noindent \emph{(A.1)} $f_t$ is not of type II,

\noindent \emph{(A.2)} $\ran f_t\cap \ran x\neq \buit$,

\noindent \emph{(A.3)} for every $s\preceq t$ if $s\in S_u$ and $f_u$ is of type I, then
for every $s'\in S_u\setminus \{s\} $,   $\ran f_{s'}\cap \ran x=\buit$,

\noindent  \emph{(A.4)} if $f_t$ is of type I  and $t\notin \mm C$,  then $t$ is a splitter of
$x$.

\noindent \emph{(A.5)} for every $u \precneqq t$, $u\notin \mm C$.
\fdefi

\prop\label{maxantex}
 $\mc A=\mc A(x, \mm C)$   is a maximal regular
antichain  such that $\{t\in \mc A\setminus C\, : \, f_t\text{ of type I}\}\con \mm S(\mc A)$.
 Moreover, if $(x_k)_{k=1}^n$ is a block sequence, then the corresponding
 $\boldsymbol{\mc A}=({\mc A}(x_k,\mm C))_{k=1}^n$ is a maximal  regular array  such that

\noindent  \emph{(a)} $\conj{t\in \bigcup_{k} \mc A_k\setminus C}{f_t\text{ of type I}}\con
\mm S(\boldsymbol{\mc A})$, and

\noindent\emph{(b)} $\mm C \con\mm C(\boldsymbol{\mc A})$ and for every $t\in \mm C$, $E_t$ is
an interval of integers.
\fprop
\prue
Fix $t\neq t'\in \mc A_k$.  $f_{t\wedge t'}$  being of type II follows from the facts that if $u
\precneqq t $, then $u\notin \mm C$, by \emph{(A.5)}, hence if $f_u$ is of type I, then
\emph{(A.3)} implies that $u$ is not splitter of $x$.  We show the maximality of $\mc A$: Fix
$t\in \mc T$ such that $\supp f_t\cap \ran x\neq \buit$. Let $t_0\succeq t$ be such that $f_{t_0}$
is of type 0 and $\supp f_{t_0}\con \ran w_k$, and set $b=[0,t_0]=\conj{v\in \mc T}{v \pe t_0}$
which is a $\pe$-well ordered set and $t\in b$. We distinguish two cases: Suppose first that
$b\cap \mm C=\buit$. Let $u_0=\min\conj{u\in b}{v \text{ satisfies } \emph{(A.1)}, \,
\emph{(A.4)}}$. Notice that $u_0$ exists since $t_0$ satisfies \emph{(A.1)} and \emph{(A.4)}. The
minimality of $u_0$ shows that $u_0$ satisfies \emph{(A.3)}, hence $u_0\in \mc A$. Suppose now
that $b\cap \mm C\neq \buit$, and set $v_0=\min b\cap \mm C$. It is not difficult to show that
$u_0=\max\conj{u \pe v_0 }{u \text{ satisfies } \emph{(A.1)}, \, \emph{(A.4)}}$ is in $\mc A $
(notice that $v_0$ satisfies \emph{(A.1)} and \emph{(A.4)}, hence $u_0$ is well defined.)

Repeating this procedure for each vector in a given a block sequence $(x_k)_{k=1}^n$, one gets
that the array $(\mc A(x_k,\mm C))_{k=1}^n$  is maximal and regular. Finally suppose that $t\in
\mm C$ and suppose that $k_1<k_2<k_3$ with $k_1,k_3\in E_t$. It is routine  to check that $t$
satisfies \emph{(A1)-(A.5)} for $x_{k_2}$, hence it follows that $k_2\in E_t$.
\fprue

\prop\label{estosonbes} Suppose that $\boldsymbol{\mc A}=(\mc A_k)_{k=1}^n$ is a regular array for a
block sequence  $(x_k)_{k=1}^n$ and  $(f_t)_{t\in \mc T}$. Then:

\noindent (b.0) If $t\in \mc A_k$ is a splitter or if $f_t$ is of type 0, then $\supp f_t\cap
\ran x_k\neq \buit$.

\noindent (b.1) If $f_t$ is  of type I, then $\{D_s\}_{s\in S_t}\cup \conj{\{k\}}{k\in E_t}$
is a block family, and if $t$ is a splitter, then $\#E_t\le \#S_t-1$.

Suppose that in addition $\boldsymbol{\mc A}=(\mc A_k)_{k=1}^n$ is maximal for $(x_k)_{k=1}^n$.

\noindent (b.2) Let $t\in \mc A_k$,  $u \precneqq s \preceq t$ with $f_u$ of type I, and
$s\in S_t$. Then for every $s'\in S_u\setminus\{s\}$ $\ran f_{s'}\cap \ran x_k=\buit$.

\noindent (b.3) Suppose that $f_t$ is of type II, $k\in D_t$ and $s\in S_t$. If  $\supp
f_t\cap \ran x_k\neq \buit$, then $k\in D_s$.
\fprop
\prue

\emph{(b.0)}:  If  $f_t$ is of type 0,   the conclusion is clear. If
  $t$ is a splitter,  let $s_1\neq s_2\in S_t$ be such that $f_{s_1}<f_{s_2}$ and $\ran f_{s_1}\cap \ran x_k, \ran
f_{s_2}\cap \ran x_k\neq \buit$. Then $\max \supp f_{s_1}\in \ran x_k$.

\emph{(b.1)}: For the first part, if  $t$ is a catcher,  there is nothing to prove, so we
assume
 $t$ is a splitter. First we show that $\{D_s\}_{s\in S_t}\cup \conj{\{k\}}{k\in E_t}$ is a
disjoint family. If $k\in E_t\cap D_s$ for some $s\in S_t$, then there is some $u \succeq s$
with $u\in \mc A_k$. But $t\in \mc A_k$ and $t \precneqq u$, a contradiction. Suppose that
$k\in D_s\cap D_{s'}$ with $s\neq s'\in S_t$. Then there are $u,u'\in \mc A_k$ such that $u
\succeq s$, $u'\succeq s'$. Hence $u\wedge u'=t$ but $f_t$ is of type I, contradicting
\emph{(a.2)}. For the second part, suppose that $k_1<k_2<k_3$ are such that $k_1,k_3 \in D_s$
for some $s\in S_t$. This implies that $\ran x_{k_1}\cap \ran f_s,\ran x_{k_3}\cap \ran
f_s\neq \buit$, and hence $\ran x_{k_2}\con \ran f_s$. This implies that $\ran x_{k_1}\cap
\ran f_{s'}=\buit$ for every $s'\in S_t\setminus \{s\}$. Since $t$ is a splitter, $k_2\notin
E_t$, and, by \emph{(a.3)}, $k_2\notin D_{s'} $ for every $s'\in S_t\setminus \{s\}$.

Let $S_t=\{s_1<\dots < s_d\}$ be ordered such that $f_{s_i}<f_{s_j}$ whenever $i<j$. For $k\in
E_t$, the set $H_k=\conj{i\in [1,d]}{\ran x_k \cap \ran f_{s_i}\neq \buit}$ has at least two
elements. We claim that the mapping $k\mapsto \max H_k\in \{2,\dots,d\}$ is one-to-one. To see
this note that  for $k<k'$ we obtain that $H_k \cap H_{k'}= \{\max H_k\}$ if $\max H_k=\max
H_{k'}$, and $H_k<H_{k'}$ otherwise.

\emph{(b.2)}: Fix $s'\in S_t\setminus \{s\}$, and suppose that $\ran f_{s'}\cap \ran x_k\neq
\buit$. Since $\ran f_{s}\cap \ran x_k\neq \buit$, we get that $\supp f_{s'}\cap \ran x_k\neq
\buit$. By the maximality of $\mc A_k$, there is $t'\in \mc A_k$ comparable with $s'$. Since
$\mc A_k$ is an antichain, we get that $t' \succeq s'$, and hence $t \wedge t'=u$. But $f_{u}$
is of type I, a contradiction.

\emph{(b.3)}: This follows using \emph{(a.4)}, and \emph{(a.1)}, \emph{(a.2)}.
\fprue

\subsubsection{Assignments, filtrations, and their relationships}
\defi
Given a block sequence $(x_{k})_{k=1}^n$, and  a regular array $\boldsymbol{\mc A}=(\mc
A_k)_{k=1}^n$ for $(x_{k})_{k=1}^n$, a sequence $(g_{k,t}^{\boldsymbol{\mc A}})_{t\in \mc
A_k,k}\con c_{00}(\N)$ is called a $\boldsymbol{\mc A}$\emph{-assignment} provided that $\supp
g_{k,t}\con \{k\}$ for every $k$ and $t\in \mc A_k$. The property \emph{(b.1)} ensures that every
$\boldsymbol{\mc A}$-assignment $(g_{k,t}^{\boldsymbol{\mc A}})_{t\in \mc A_k,k}$ naturally
\emph{filters down} to the whole tree $(G_{k,t}^{\boldsymbol{\mc A}})_{t\in \mc T}$ as follows: If
$k\notin D_t^{\boldsymbol{\mc A}}$, then $G_{k,t}^{\boldsymbol{\mc A}}=0$, and if $t\in \mc A_k$,
then $G_{k,t}^{\boldsymbol{\mc A}}=g_{k,t}^{\boldsymbol{\mc A}}$. Suppose that $k\in
D_t^{\boldsymbol{\mc A}}\setminus D_s^{\boldsymbol{\mc A}}$. If $f_t$ is of type I, then we define
recursively $G_{k,t}^{\boldsymbol{\mc A}}=(1/w(f_t))G_{k,s}^{\boldsymbol{\mc A}}$, where $s\in
S_t$ is the unique $s=s(k,t)\in S_t$ such that $k\in D_s^{\boldsymbol{\mc A}}$ (by \emph{(b.1)}).
If $f_t$ is of type II, $f_t=\sum_{s\in S_t}\la_sf_s$, then we simply set
$G_{k,t}^{\boldsymbol{\mc A}}=\sum_{s\in S_t}\la_s G_{k,s}^{\boldsymbol{\mc A}}$. For $t\in \mc
T$, let
$$
G_t^{\boldsymbol{\mc A}}=\sum_{k\in
D_t^{\boldsymbol{\mc A}}}G_{k,t}^{\boldsymbol{\mc A}}.$$ We call $(G_t^{\boldsymbol{\mc A}})_{t\in
\mc T}$ the \emph{filtration} of $(g_{k,t}^{\boldsymbol{\mc A}})_{t\in \mc A_k,k}$. Whenever there
is no possible confusion, we write $g_{k,t}$, $G_{k,t}$ and $G_t$ instead of the respective
$g_{k,t}^{\boldsymbol{\mc A}}$, $G_{k,t}^{\boldsymbol{\mc A}}$ and $G_t^{\boldsymbol{\mc A}}$.
\fdefi

\prop\label{estoesc}
Fix $t\in \mc T$. Then we have the following:

\noindent (c.1)  For every $k$, $\supp g_{k,t}\con \{k\}$. Hence  $\supp g_t\con D_t$.

\noindent (c.2) If $f_t$ is not of type II, then $G_t=\sum_{k\in
E_t}g_{k,t}+(1/w(f_t))\sum_{s\in S_t}G_s$.

\noindent (c.3) If $f_t$ is of type II,$f_t=\sum_{s\in S_t}\la_sf_s$, then $G_t=\sum_{s\in S_t}\la_sG_s$.
\fprop
\prue (c.1) is clear.
(c.2): If $f_t$ is of type 0, this is clear. Suppose that $f_t$ is of type I. Then by definition
\begin{align}
G_t= &   \sum_{k\in E_t}G_{k,t}+\sum_{k\in D_t\setminus E_t}G_{k,t} = \sum_{k\in
E_t}g_{k,t}+\sum_{s\in S_t}\sum_{k \in D_s}G_{k,t}=\nonumber \\= &   \sum_{k\in
E_t}g_{k,t}+\sum_{s\in S_t} \frac{1}{w(f_t)}\sum_{k \in D_s}G_{k,s} =  \sum_{k\in
E_t}g_{k,t}+\frac{1}{w(f_t)}\sum_{s\in S_t} G_{s}.
\end{align}

(c.3): Suppose that $f_t$ is of type II, i.e.,  $f_t=\sum_{s\in S_t}\la_sf_s$, and suppose that $k\in D_t$.
Then, by (c.1), $G_t(e_k)=G_{t,k}(e_k)=\sum_{s\in S_t}\la_s G_{k,s}(e_k)=(\sum_{s\in S_t} \la_s G_{s})(e_k)$.
If $k\notin D_t$, then $G_t(e_k)=0$, and $\sum_{s\in S_t}\la_sG_s(e_k)=0$.\fprue

\defi{\emph{(Canonical Assignment)}} \label{canonassig}
Suppose that $\boldsymbol{\mc A}=(\mc A_k)_{k}$ is a regular array for $(x_k)_{k=1}^n$ and
$(f_t)_{t\in \mc T}$. Let $f_{k,t}=f_t(x_k)e_k^*$ for $k\in[1,n]$  and  $t\in \mc A_k$. This is
the \emph{$\boldsymbol{\mc A}$-canonical assignment.}
\fdefi
\nota Note that if the array $\boldsymbol{\mc A}$ is maximal, then filtering down the canonical
assignment we get $f_t(w_k)=F_{k,t}(e_k)$,  for every $t\in \mc T$, and $k\in D_t$: If $k\in E_t$,
  this is just by definition. Suppose  $k\notin E_t$. If $f_t$ is of type I, then $F_{k,t}
(e_k)=(1/w(f_t))F_{k,s}(e_k)$, where $s\in S_t$ is unique such that $k\in D_s$. By the
maximality of $\mc A_k$, we get   that $\supp f_{s'}\cap \ran w_k=\buit$ for every $s'\in
S_t\setminus \{s\}$ (by \emph{(b.2)}), hence $f_t(x_k)=(1/w(f_t))f_s(x_k)=(1/w(f_t))F_{k,s}
(e_k)=F_{k,t} (e_k)$,  by the inductive hypothesis. If $f_t=\sum_{s\in S_t}\la_s f_s$ is of
type II, then by the maximality of $\mc A_k$, $f_t(x_k)=\sum_{s\in S_t, k\in D_s} \la_s
f_s(x_k)=\sum_{s\in S_t, k\in D_s} \la_s F_{k,t}(e_k)=F_{k,t} (e_k)$, the last equality
because $F_{k,u}=0$ if $k\notin D_u$.

We obtain that $f_t(\sum_{k\in D_t}b_kx_k)=F_t(\sum_{k\in D_t}b_k e_k)=F_t(\sum_{k=1}^n b_k e_k)$
for every sequence of scalars $(b_k)_{k=1}^n$. The last equality follows from $\supp G_t\con D_t$.
In particular, $f(\sum_{k=1}^n b_k x_k)=G_\buit(\sum_{k=1}^n b_k e_k)$, since
$D_\buit=\conj{k}{\supp f\cap \ran x_k\neq \buit}$, by maximality of $\boldsymbol{\mc A}$.

 \fnota

\prop\label{desigualdades} Suppose that $\boldsymbol{\mc A}=(\mc A_k)_{k=1}^n$ is a regular array (not necessarily  maximal)
for $(x_k)_{k=1}^n$ and $(f_t)_{t\in \mc T}$. Fix scalars $(b_k)_{k=1}^n$, $(c_k)_{k=1}^n$ and
suppose that $(g_{k,t})_{t\in \mc A_k ,k}$,   $(h_{k,t})_{t\in \mc A_k ,k}$  are $\boldsymbol{\mc
A}$-assignments.

\noindent{(1)}  If for every  $t\in \mc A_k$ $g_{k,t} (b_ke_k)\le h_{k,t} (c_ke_k)$,  then for
every $t\in \mc T$, $G_{k,t} (b_ke_k)\le H_{k,t} (c_ke_k)$.

\noindent{(2)} $\nrm{G_{k,u}(e_k)}_\infty\le \max\conj{\nrm{g_{k,t}}_{\infty}}{t\in \mc A_k}$, for
every $u\in \mc T$.

\noindent{(3)} If for every $t\in \bigcup_{k=1}^n \mc A_k$ $\sum_{k\in E_t}g_{k,t} (b_k e_k)\le
\sum_{k\in E_t}h_{k,t} (c_k e_k)$, then for every $t\in \mc T$, $G_t (\sum_{k\in D_t}b_ke_k)
    \le H_t(\sum_{k\in D_t}c_k e_k)$.

\noindent{(4)}  $\nrm{G_u}_{\infty}\le \nrm{\sum_{t\in A_k, k} g_{k,t}}_\infty$ for every $u\in
\mc T$.

\fprop
\prue
This follows  from Proposition \ref{estoesc}.
\fprue

\subsubsection{Two successive filtrations}\label{twofiltra}
In some applications of the theory one needs to do the process of assignment and filtration twice
starting with different arrays of antichains. To see this, suppose that $\boldsymbol{\mc C}$ and
$\boldsymbol{\mc D}$ are  regular arrays for $(x_k)_{k=1}^n$ and $(f_t)_{t\in \mc T}$. Then we can
naturally define a $\boldsymbol{\mc D}$-assignment $(g_{k,t}^{\boldsymbol{\mc D}})_{t\in \mc
D_k,k}$ by taking the filtration $g_{k,t}^{\boldsymbol{\mc D}}=G_{k,t}^{\boldsymbol{\mc C}}$. For
this to work, one needs the following special relationship between $\boldsymbol{\mc C}$ and
$\boldsymbol{\mc D}$.

\defi
We write $\boldsymbol{\mc C}\nprec \boldsymbol{\mc D}$ if for every $k$, every $c\in \mc C_k $ and
every $d\in \mc D_k$, we have that $c \nprec d$. A $\boldsymbol{\mc C}$-assignment
$(g_{k,t}^{\boldsymbol{\mc C}})_{k\in \mc C_k,k}$ is called \emph{coherent} provided that
$g_{k,t}^{\boldsymbol{\mc C}}=0$ whenever $f_t(w_k)=0$.
\fdefi

\prop \label{iurppp} Suppose that $\boldsymbol{\mc C} \nprec \boldsymbol{\mc D}$ are two regular
arrays for $(x_k)_{k=1}^n$ and $(f_t)_{t\in \mc T}$, and suppose that $\boldsymbol{\mc D}$ is in
addition maximal. Fix a coherent $\boldsymbol{\mc C}$-assignment $(g_{k,t}^{\boldsymbol{\mc
C}})_{t\in \mc C_k,k}$. Then,

\noindent{(a)} For every $k\in D_t^{\boldsymbol{\mc C}}\cap D_t^{\boldsymbol{\mc D}}$,
$G_{k,t}^{\boldsymbol{\mc C}}=G_{k,t}^{\boldsymbol{\mc D}}$.

\noindent (b) $g_\buit^{\boldsymbol{\mc C}}=g_\buit^{\boldsymbol{\mc D}}$. \fprop

\prue (a): If $k\in E_t^{\boldsymbol{\mc D}}$, this is just by definition. Suppose   $f_t$ is
of type I, and suppose that $k\in D_{s}^{\boldsymbol{\mc D}}$, for some $s\in S_t$. Then
$G_{k,t}^{\boldsymbol{\mc D}}=(1/w(f_t))G_{k,s}^{\boldsymbol{\mc D}} $. Since
${\boldsymbol{\mc D}}$ is a maximal regular array, by Proposition \ref{estosonbes}
\emph{(b.2)},  $\ran f_{s'}\cap \ran w_k=\buit $ for every $s'\in S_t\setminus \{s\}$.  If
$k\in D_s^{\boldsymbol{\mc C}}$, then we are done by the inductive hypothesis. So, suppose
$k\in E_t^{\boldsymbol{\mc C}}$, i.e., $t\in {{\mc C}}_k$. Hence, $t\precneqq u$ for some
$u\in \mc D_k$ (because $k\in D_s^{\boldsymbol{\mc D}}$),  contradicting our assumption that
$\boldsymbol{\mc C}\nprec \boldsymbol{\mc D}$. If $f_t=\sum_{s\in S_t}\la_s f_s$ is a
sub-convex combination, then
\begin{equation}
G_{k,t}^{\boldsymbol{\mc D}}=\sum_{s\in S_t,k\in D_s^{\boldsymbol{\mc
D}}}G_{k,s}^{\boldsymbol{\mc D}}=\sum_{s\in S_t,k\in D_s^{\boldsymbol{\mc D}}\cap
D_s^{\boldsymbol{\mc C}}}G_{k,s}^{\boldsymbol{\mc D}}=\sum_{s\in S_t,k\in D_s^{\boldsymbol{\mc
D}}\cap D_s^{\boldsymbol{\mc C}}}G_{k,s}^{\boldsymbol{\mc D}}=G_{k,t}^{\boldsymbol{\mc C}}.
\end{equation}
To see the last equality note that if $k\notin D_s^{{\boldsymbol{\mc D}}}$, then, by the
maximality of ${\boldsymbol{\mc D}}$, $\supp f_u\cap \ran w_k=\buit$ for every $u\succeq s$,
so, by the coherence of the assignment, $G_{k,t}^{\boldsymbol{\mc C}}=0$; if $k\notin
D_s^{{\boldsymbol{\mc C}}}$, then $k\notin D_u^{\boldsymbol{\mc C}}$ for all $u\succeq s$, and
so $g_{k,u}^{\boldsymbol{\mc C}}=0$ for all $u\succeq s$ $u\in \mc C_k$, giving us
$G_{k,s}^{\boldsymbol{\mc D}}=0$.

\noindent (b): Note now that
\begin{equation}
g_\buit^{\boldsymbol{\mc D}}=\sum_{k\in D_\buit^{\boldsymbol{\mc
C}}}g_{k,\buit}^{\boldsymbol{\mc D}}=\sum_{k\in D_{\buit}^{\boldsymbol{\mc C}}\cap
D_{\buit}^{\boldsymbol{\mc C}}}g_{k,\buit}^{\boldsymbol{\mc D}}=g_{\buit}^{\boldsymbol{\mc
C}}.
\end{equation}
For if $k\in D_\buit^{\boldsymbol{\mc C}}\setminus D_\buit^{\boldsymbol{\mc D}}$, then
by the maximality of ${\boldsymbol{\mc D}}$, for all $u\in \mc T$, $\supp f_u\cap \ran
w_k=\buit$, hence, by the coherence of the $\boldsymbol{\mc C}$-assignment
$g_{k,u}^{\boldsymbol{\mc C}}=0$ for all $u$, and hence $g_{k,\buit}^{\boldsymbol{\mc D}}=0$;
if $k\in D_\buit^{\boldsymbol{\mc D}}\setminus D_\buit^{\boldsymbol{\mc C}}$, then   every
$\mc C_k=\buit$, and so $g_{k,\buit}^{\boldsymbol{\mc C}}=0$. \fprue

Let us now give the two main applications of this general theory of tree-analysis.

\subsection{The proof of the basic inequality}\label{dwaknvnbvb}

Recall that $W$ is the minimal subset of $c_{00}(\N)$ containing $\{\pm e_k^*\}_k$, and
closed under $(m_j^{-1},n_{4j})$-operations.  Fix  a $(C,\vep)$-RIS $(x_k)_{k=1}^n$, and
fix $({j_k})_{k=1}^n$ witnessing that $(x_k)_{k=1}^n$ is a $(C,\vep )$-RIS, i.e.,

\noindent a) $\nrm{x_k}\le C$,

\noindent b) $|\supp x_k| \le m_{j_{k+1}}\vep$ and

\noindent c) For all type I functionals $\phi$ of $K$ with $w(\phi)<m_{j_k}$, $|\phi(x_k)|\le
C/w(\phi)$.  Fix a sequence $(b_k)_{k=1}^n$ of scalars, $\max_k |b_k|\le 1$, and  $f\in
K_{\ou}$. Let $(f_t)_{t\in \mc T}$ be a tree-analysis of $f$. Consider the maximal regular
array $\boldsymbol{\mc A}=(\mc A(x_k,\mm C))_{k=1}^n$, where $\mm C$ is the set of nodes $t$
such that $f_t$ is of type I and $w(f_t)=m_{j_0}$.

We introduce the following two $\boldsymbol{\mc A}$-assignments $(g_{k,t})_{t\in \mc A_k,k}$,
and $(r_{k,t})_{t\in \mc A_k,k}$. Fix $k$ and $t\in \mc A_k$. If $t_t$ is of type 0, then we
set $g_{k,t}=e_k^*$ and $r_{k,t}=0$. Suppose that $t$ is of type I, and $w(f_t)\neq m_{j_0}$.
Let
\begin{equation}
l_t=\min\conj{k\in E_t}{w(f_t)\le m_{j_{k}}}
\end{equation}
 if this exists, and $l_t=\infty$ otherwise. Then let

\begin{align*}
g_{k,t}=\left\{\begin{array}{ll}
\frac1{w(f_t)}e_k^* &\text{ if $k>l_t$}\\
0& \text{ if $k< l_t$}\\
e_k^* &\text{ if $k=l_t$}
\end{array}
\right. & \, r_{k,t}=\left\{\begin{array}{ll}
0&\text{ if $k>l_t$}\\
\vep e_k^*& \text{ if $k< l_t$}\\
0&\text{ if $k=l_t$}
\end{array}
\right.
\end{align*}
Suppose now that $w(f_t)=m_{j_0}$. Notice that $E_t$ is an interval. Set
\begin{equation}
k_t=\max\conj{l\in D_t}{|b_l|=\nrm{(b_i)_{i\in E_t}}_\infty}.
\end{equation}
Then let
\begin{align*}
g_{k,t}=\left\{\begin{array}{ll}
e_k^* &\text{ if $k=k_t$}\\
0& \text{ if $k\neq k_t$}\\
\end{array}
\right. & \, r_{k,t}= \vep e_k^*
\end{align*}
Let $(G_t)_{t\in \mc T}$ and $(R_t)_{t\in \mc T}$ be the corresponding filtrations.

\clam[D] Fix $t\in \mc T$. Then:

\noindent (d.1) $\nrm{R_t}_\infty\le \vep$.

\noindent (d.2) $|f_t(\sum_{k\in D_t}b_kx_k)|\le C(G_t+R_t)(\sum_{k\in D_t}|b_k|e_k)$.

\noindent (d.3) For every $t$ for which $f_t$ is of type I, either
 $G_t \in \conv \conj{h\in W}{w(h)=w(f_t)}$  or $G_t =e_{k}^*+h_t$ for some
 $k\notin \supp h_t$ and $h_t\in \conv \conj{h\in W}{w(h)=w(f_t)}$.
\fclam \prucl \emph{(d.1)} follows  from Proposition \ref{desigualdades}, and \emph{(d.2)} follows
also from the same proposition applied to the canonical $\boldsymbol{\mc A}$-assignment, the
assignment $(C(g_{k,t}+h_{k,t}))_{t\in \mc A_k,k}$, and the sequences of scalars $(b_k)_k$ and
$(|b_k|)_k$.

\emph{(d.3)}:  If $w(f_t)=m_{j_0}$, then $t$ is a catcher and $G_t=\sum_{k\in
E_t}g_{k,t}=e_{k_t}^*\in W$.  Suppose that $t$ is of type I, $w(f_t)\neq m_{j_0}$. By \emph{(c.2)}
and the particular $\boldsymbol{\mc A}$-assignment, we know that either
$G_{t}=(1/w(f_t))(\sum_{k\in E_t, k>l_t}e_k^*+\sum_{s\in S_t} G_{s})$ or $G_{t} =e_{l_t}^*+h_t$,
where $h_t=(1/w(f_t))(\sum_{k\in E_t, k>l_t}e_k^*+\sum_{s\in S_t}G_{s} )$. Assume this last case
holds.

\noindent \textit{Subcase 1a}. For every $s\in S_t$ the functional
 $f_{s}$ is not of type II. From the inductive hypothesis, we  have that for  every $s\in S_{t}$, $G_{s}=h_{s}$ or
 $G_{s}=e^{*}_{l_s}+h_{s}$, $h_{s}\in W$. For $s\in S_{t}$ such that   $G_{s} =e^{*}_{l_{s}}+h_{s}$, set
 $I_{s}^{1}=\{n\in\mathbb{N}: n<l_{s}\}$ and  $I_{s}^{2}=\{n\in\mathbb{N}: n>l_{s}\}$. We set
 $h_{s}^{1}=I_{s}^{1}h_{s}$,  $h_{s}^{2}=I_{s}^{2}h_{s}$. Then, for every $s\in  S_{t}$ the functionals $h_{s}^{1}$,
 $e^{*}_{l_{s}}$, and  $h_{s}^{2}$ are successive  and belong to $W$. By \emph{(b.1)}, for
 $s\not=s^{\prime}\in S_{t}$ the corresponding  functionals together with $\{e_{k}^*\}_{k\in E_t,k>l_t}$
  form a block family, and we obtain  that
  \begin{equation}   \label{sixes}
\#\{e_k^*\}_{k\in E_t,k>l_t}+\#\conj{e_{l_s}^*}{s\in S_t}+\#\conj{h_s^1}{s\in S_t}+
\#\conj{h_s^2}{s\in S_t}\le 4 \#S_t.
\end{equation}
Therefore, $(1/w(f_t))(\sum_{k\in E_t, k>l_t}e_k^*+\sum_{s\in S_t}G_{s})\in W$.

 \noindent\textit{Subcase 1b}. There are $s\in S_t$ for which  $f_s $
 is of type II. Let $B_1$ be the set of immediate successors $s$ of $t$ such that $f_s$ is of type II, and $B_2=S_t\setminus S_1$.
 Observe that every sub-convex combination $f_s =\sum_{u\in S_s }r_{u}f_{u}$
 satisfies that $f_u$ is of type I. We may assume, allowing repetitions if needed, that  for every
 $s\in S_t$ such that $f_s$ is of type II, $f_s=(1/k)\sum_{q=1}^k f_{s,q}$, where each $f_{s,q}\in \{f_{u}\, : \,u\in
 S_s\}$. For each  $q=1,2,\ldots ,k$ we set
 $h_t^q=   ( {1}/{m_j})(  \sum_{l\in E_t, l>l_t} e_{l}^*+\sum_{s\in B_1}G_s
 +\sum _{s\in B_2}G_{s,q})$, where $G_{s,q}=G_{u}$ for ${u}\in S_s $ such that
 $f_{s,q}=f_{u}$. A similar argument as in  the previous subcase shows that $h_t^q\in W$
 with $w(h_t^q)=m_j$  for $q=1,2,\ldots ,k$
 and  $h_t=({1}/{k})\sum_{q=1}^kh_t^q$, as required.\fprucl

The particular case  $t=\emptyset$, the root of $\mc T$,  gives us the conclusion of the Basic
Inequality.

\nota Note that a finer assignment using the same array of antichains will actually give us the
conclusion of the Basic Inequality for a bit smaller auxiliary space $T[(m_j^{-1},2n_j)_j]$.
\fnota

\subsection{The proof of the finite interval representability of $J_{T_0}$}
The general scheme of the proof is quite similar to the proof of Basic Inequality though the
input block sequence of vectors is slightly differently  chosen. Notice however that the
finite interval representability involves two inequalities needed  for showing that the
representing operator as well as its inverse are uniformly bounded. Thus, while in the proof
of the Basic Inequality we could afford to go the auxiliary space $T[(m_j^{-1},4n_{j})_j]$
this is no longer possible in this case. In other words, we need to improve on the counting
inequality (\ref{sixes}). It is exactly for this reason that we introduce below two arrays of
antichains and use two successive filtrations as explained above in Subsection
\ref{twofiltra}.

Fix a transfinite block sequence $(x_\al)_{\al<\ga}$, $n\in \N$,   a sequence $I_1\le I_2\le \dots
\le I_n$ of successive, not necessarily distinct, infinite intervals of $\ga$, and $\vep>0$. Let
$j_0$ be such that $m_{2j_0+1}>100 n /\vep$  and set $l=n_{2j_0+1}/m_{2j_0+1}$. Find a
$(1,j_0)$-dependent sequence $(z_1,\psi_1,\dots,z_{n_{2j_0+1}},\psi_{n_{2j_0+1}})$ such that (a)
$\ran \psi_i\con \ran z_i$ for every $i=1,\dots,n_{2j_0+1}$ and (b) $(z_k)_{k=(i-1)l+1}^{il}\con
\langle x_\al \rangle_{\al\in I_i}$ for every $i=1,\dots,n$. Let
\begin{align*}
\phi= \frac1{m_{2j_0+1}}\sum_{i=1}^{n_{2j_0+1}} \psi_i,
\end{align*}
and for each $k=1,\dots,n$ we set
\begin{align*}
w_k= &   \frac{m_{2j_0+1}}l \sum_{i=(k-1)l+1}^{kl} z_i  \text{ and }\phi_k=
\frac1{m_{2j_0+1}} \sum_{i=(k-1)l+1}^{kl} \psi_i \in K_{\ou}.
\end{align*}

\prop\label{estonw} Fix  $k=1,\dots,n$. Then

\noindent{(1)} $\ran \phi_k\con \ran w_k$, $\phi_k w_k=1$ and    $1\le \nrm{w_k}\le  24$.

\noindent{(2)} For every $f\in K_{\ou}$ of type I with $w(f)>m_{2j_0+1}$,
$|f(w_k)|\le {1}/{m_{2j_0+1}^2}$.

\noindent{(3)} Let $f\in K_{\ou}$ be of type I, $f=(1/w(f))\sum_{i=1}^d f_i$ with $w(f)=m_{2j+1}$
for $j<j_0$ and $d\le n_{2j+1}$. Let $ d_0=\max\conj{i\le d}{w(f_i)<m_{2j_0+1}}$, and set
$f_\mathrm{L}=1/{m_{2j+1}}\sum_{i=1}^{d_0-1}f_i$  and $ f_\mathrm{R}=1/{m_{2j+1}}\sum_{i=d_0+1}^{d}f_i$. Then
$ |f_\mathrm{L}(w_k)|\le  1/{m_{2j_0+1}^2} \text{ and }|f_\mathrm{R}(w_k)|\le {1}/{m_{2j_0+1}}$.

\noindent{(4)} Let $f= (1/w(f))\sum_{i=1}^d f_i$ with $w(f)=m_{2j_0+1}$ and $d\le n_{2j_0+1}$ be
such that $\#\{i\in [1,d]\, :\, w(f_i)=w(\psi_i)\text{ and  }\supp z_i\cap \supp f_i\neq \buit\}\le
2$. Then, $|f(w_k)|\le 1/m_{2j_0+1}^2$.
\fprop

\prue First of all, note that $(z_i)_{i=(k-1)l+1}^{kl}$ is a $(12, 1/n_{2j_0+1})$-RIS. Note
also that  \emph{(1)} and \emph{(2)} follow  from Proposition \ref{niceasris}.  \emph{(3)} By the
properties of special sequences,
\begin{equation}
\# \bigcup_{i=1}^{d_0-1}\supp f_i \le w(f_{d_0})<m_{2j_0+1}.
\end{equation}
So, $ |f_\mathrm{L}(w_k)|\le \nrm{f_0}_{\ell_1}\nrm{w_k}_\infty\le {m_{2j_0+1}^3}/{n_{2j_0+1}}\le
1 /{m_{2j_0+1}^2}$. Let us now estimate $|f_\mathrm{R}(w_k)|$.  To save on notation we only
estimate for $k=1$. Set
\begin{equation*}
 F_0=\{r\in [1,l] \, : \, \#(\conj{i\in [d_0+1,d]}{\range z_r \cap \supp f_i \neq \buit})\ge 2\},\,
F_1= [1,l] \setminus F_0.
\end{equation*}
Notice that $|F_0|\le n_{2j+1}-1$. For $i=0,1$ let $w^i= ({m_{j_1}}/{l}) \sum_{k\in F_i } z_k$.
Since $f_\mathrm{R}\in K_{\ou}$ and since $(z_k)_k$ is a $(12,1/n_{2j_0+1})$-RIS, we have that
\begin{equation}  \label{mjwere1}
|f_\mathrm{R}(w^0)|\le \nrm{w^0}\le \frac{m_{2j_0+1}}{l}\sum_{k\in F_0 }\nrm{ z_k}\le \frac{m_{2j_0+1}}{l}6
n_{2j+1}.
\end{equation}
To estimate $|f_\mathrm{R}(w^1)|$ we use the basic inequality. For each $i=d_0+1,\dots,d$, let
\begin{align*}
H_i=& \conj{k\in F_1}{\range z_k \cap \supp f_i \neq \buit}.
\end{align*}
Note that $\{H_i\}_i$ is a partition of $F_1$ and is a block family. For  $i=d_0+1,\dots,d$,
 we set $w^{1,i}={m_{j_1}}/{l} \sum_{k\in H_i}z_k$.
Clearly $w^1=w^{1,d_0+1}+\dots +w^{1,d}$ and hence
\begin{equation}
|f_\mathrm{R}(w^{1})|\le \sum_{i=d_{0}+1}^d |f_\mathrm{R}(w^{1,i})|= \frac1{m_{2j+1}} \sum_{i=d_0+1}^d
|f_i(w^{1,i})|.
\end{equation}
Let us estimate now $| f_i(w^{1,i})|$, for $i=d_0+1,\dots,d$. For a fixed such  $i$,
applying again the basic inequality, we obtain that $| f_i(w^{1,i})|\le
12(g_1^i+g_2^i)({m_{2j_0+1}}/{l}\sum_{k\in H_i} e_k)$, where in the worst case,
$g_1^i=h^i+e_{k_i}^*$, with $h^i\in W$, and $h^i\in\conv_\Q\conj{h\in W}{w(h)=w(f_i)}$.
Since the auxiliary space is 1-unconditional, by Proposition \ref{estbasis}, $|h^i(
({m_{2j_0+1}}/{l})\sum_{k\in H_i} e_k)|\le {m_{2j_0+1}}/{w(f_i)}$. Note that
$\nrm{g_2^i}_\infty \le 1/n_{2j_0+1}$. Putting all these inequalities together we get
\begin{align} \label{mjwere2}
|f_\mathrm{R}(w^{1})|\le &    \frac{12}{m_{2j+1}}\left( \sum_{i=d_0+1}^d
\frac{m_{2j_0+1}}{w(f_i)}+\frac{m_{2j_0+1}n_{2j+1}}{l}+\frac {m_{2j_0+1}}{n_{2j_0+1}}
\right)\le \nonumber
\\
\le &   \frac{12}{m_{2j+1}}\left( \sum_{i=d_0+1}^d
\frac{m_{2j_0+1}}{w(f_i)}+\frac{m_{2j_0+1}^2n_{2j+1}}{n_{2j_0+1}}+\frac{
m_{2j_0+1}}{n_{2j_0+1}} \right).
\end{align}
Using (\ref{mjwere1}) and (\ref{mjwere2}) we obtain
\begin{equation}  \label{uuuno2}
|f_\mathrm{R}(w_1)|\le  \frac{12 m_{2j_0+1}}{m_{2j+1}}\left(
\frac{2n_{2j+1}m_{2j_0+1}}{n_{2j_0+1}}+\frac1{n_{2j_0+1}}+\sum_{i=d_0+1}^d \frac1{w(f_i)}\right)\le
\frac{1}{m_{2j_0+1}}.
\end{equation}

 \emph{(4)} Let $E= \conj{i\in [1,d]}{w(f_i)=w(\psi_i)\text{ and  }\supp z_i\cap \supp f_i\neq \buit}$. By
our assumptions, $\# E\le 2 $.  For $i\in [(k-1)l,kl]\setminus E$ the properties of the dependent
sequences yield that $|f(z_i)|\le {1}/{n_{2j_0+1}}$. Hence, $|f(w_k)|\le 2\cdot 24
 {m_{2j_0+1}}/{l}+ {m_{2j_0+1}}/{n_{2j_0+1}}\le 1/{m_{2j_0+1}^2}$.
\fprue

\lema\label{theotherin}
For the above defined sequence  $( w_k)_k$  we have that
\begin{equation}
\nrm{\sum_{k=1}^n b_k w_k}\le 121 \nrm{\sum_{k=1}^n b_k v_k}_{J_{T_0}}
\end{equation}
  for every choice of
scalars $(b_k)_{k=1}^n$.
\flema
\prue Fix a sequence $(b_k)_{k=1}^n$ of scalars with $\max_k |b_k|\le 1$, an  $f\in K_{\ou}$,
and its  tree  $(f_t)_{t\in \mc T}$.

\noindent{\emph{Antichains.}} A node $t\in \mc T$ is called \emph{relevant} if \emph{(1)}
$w(f_t)\le m_{2j_0+1}$, and \emph{(2)} if $u\precneqq t$ is its immediate predecessor, if $f_u$ is
of type I, and if $w(f_u)=m_{2j+1}<m_{2j_0+1}$, then $t= s(u)=\max\conj{s\in
S_u}{w(f_s)<m_{2j_0+1}}$, where the maximum is taken according to the block
  ordering $S_u=\{s_1<\dots<s_d\}$. Let $\mm C$ be the set of nodes $t$ which are either
non-relevant, or such that $f_t$ is of type I and $w(f_t)=m_{2j_0+1}$. Let ${\boldsymbol{\mc
B}}=(\mc B_k)_{k=1}^n$ where $\mc B_k=\mc A(w_k,\mm C)$ for $k=1,\dots,n$ (see \emph{(A.1)-(A.5)}
in Proposition \ref{maxantex} above). For each $k$, let $\mc B_k^\mm{unc}= \mm S(\mc B_k)\setminus
\mm C$ be the set of splitters that are not in $\mm C$, $\mc B_k^\mm{cnd}=\mc B_k \cap \mm C$, and
$\mc B_k^\mm{at}=\mc B_k\setminus (\mc B_k^\mm{unc}\cup \mc B_k^\mm{cnd})$.

Fix $u\in \mc B_k^\mm{unc}$, and observe that $u$ is a splitter of $x_k$ for every $k\in E_u$.
List all $s\in S_u$ such that $\ran f_s\cap \ran w_k\neq \buit$ ,$\{s_{k,1},\dots,s_{k,d}\}$
ordered
 according to  the block ordering $f_{s_{k,1}}<\dots <f_{s_{k,d}}$. Set
\begin{align*}
w_{k,u}^{\mathrm{in}}= &w_k |[\min \supp w_k,\max \supp f_{s_{k,1}}]  \\
w_{k,u}^\mathrm{fin}= & w_k -w_{k,t}^\mathrm{in}.
\end{align*}
For $\star  \in   \{\mathrm{in},\mathrm{fin}\}$, let  $\mc B_{k,u}^\star=\mc A(w_{k,u}^\star, \mm
C^{nr})$, where $\mm{C}^{nr}$ is the set of non-relevant nodes of $\mc T$. Set $\mc B_k^\star
=\bigcup_{u\in \mc B_k^\mm{unc}}\mc B_{k,u}^\star$. Observe that $\boldsymbol{\mc B^\star}=(\mc
B_k^\star)_{k=1}^n$ is a regular (not necessarily maximal) arrow for $(w_k)_{k=1}^n$ and
$(f_t)_{t\in \mc T}$, whenever $\star\in \{\mm{in},\mm{fin},\mm{cnd},\mm{at}\}$.

\noindent{\emph{Assignments and filtrations}}. Consider the following $\boldsymbol{\mc
B^\star}$-assignments $(g_{k,t}^{\star})_{k\in \mc B_k^\star,k}$ where $\star\in
\{\mm{in},\mm{fin},\mm{cnd}\}$, and $(r_{k,t}^{\star})_{k\in \mc B_k^\star,k}$  where
$\star\in\{\mm{in},\mm{fin},\mm{cnd},\mm{at}\}$: Fix $k$, and  $t\in \bigcup_{\star\in
\{\mm{in},\mm{fin},\mm{cnd},\mm{at}\} }\mc B_k^\star$.

\noindent{(a)} Suppose that $f_t$ is of type 0. Then we set
$r_{k,t}^{\mm{at}}=(1/m_{2j_0+1})e_{k}^*$ if $t\in \mc B_k^\mm{at}$, and  we set
$g_{k,t}^{\star}=0$ and $r_{k,t}^{\star}=(1/m_{2j_0+1})e_{k}^*$, if $t\in \mc B_k^\star$ for some
$\star\in\{\mm{in},\mm{fin}\}$.

\noindent{(b)} Suppose that $t$ is non-relevant. Then clearly $t\notin  \mc B_{k}^\mm{at}$. Fix
$\star\in\{\mm{in},\mm{fin},\mm{cnd}\}$. We set $g_{k,t}^{\star}=0$ in all cases. Suppose that
$w(f_t)>m_{2j_0+1}$. Then we set $r_{k,t}^{\star}=(1/m_{2j_0+1})e_{k}^*$ for
$\star\in\{\mm{in},\mm{fin}\}$, and $r_{k,t}^\mm{cnd}=(\sgn(b_k)/m_{2j_0+1})e_{k}^*$. Finally, if
$t\neq s(u)$, where $u$ is the immediate predecessor of $t$ (see the definition of relevant node),
then we set $r_{k,t}^{\star}=\nrm{f_t(w_k)}e_{k}^*$ for  $\star\in\{\mm{in},\mm{fin}\}$ and
$r_{k,t}^{\mm{cnd}}=\sgn(b_k)\nrm{f_t(w_k)}e_{k}^*$.

\noindent{(c)} Suppose now that $t$ is relevant. There are two subcases.

\noindent{(c.1)} $w(f_t)=m_{2j_0+1}$. If $t\in \mc B_k^\star$, for $\star\in
\{\mm{in},\mm{fin}\}$, then we set $g_{k,t}^{\star}=(1/w(f_t))e_k^*$ and $r_{k,t}^{\star}=0$.
Suppose that  $t\in \mc B_k^\mm{cnd}$. Suppose that $f_t=\pm I
(1/m_{2j_0+1})\sum_{i=1}^{n_{2j_0+1}} g_i$, where $I\con \ou$ is an interval, and
$\Phi=(g_1,\dots,g_{n_{2j_0+1}})$ is a $2j_0+1$-special sequence. Set
$\Psi=(\psi_1,\dots,\psi_{n_{2j_0+1}})$. Consider $I_t=\conj{i\in
[1,\ka_{\Phi,\Psi}-1]}{Ig_i\neq 0} =[k({t,1}),k({t,2})]$, and let
$\vep_t=\sgn(\sum_{k=k(t,1)+1}^{k(t,2)-1}b_k)$. If $k=k(t,i)$ for $i=1,2$, then we set
$g_{k,t}^{\mm{cnd}}=\sgn{(b_{k(t,i)})} e_{k(t,i)}^*$ and $r_{k,t}^{\mm{cnd}}=0$. We set
$g_{k,t}^{\mm{cnd},\mc B_k^\mm{cnd}}=\vep_t e_k^*$ and $r_{k,t}^{\mm{cnd},\mc B_k^\mm{cnd}}=0$
if $k\in (k(t,1),k(t,2))$. We set   $g_{k,t}^{\mm{cnd}}=0$, and
$r_{k,t}^{\mm{cnd}}=\sgn(b_k)(1/m_{2j_0+1})e_k^*$ otherwise.

\noindent{(c.2)} Suppose that $w(f_t)\neq m_{2j_0+1}$. Then $t\in \mc B_k^\star$, for some
$\star\in \{\mm{in},\mm{fin}\}$. Set $g_{k,t}^{\star}=(1/w(f_t))e_{k}^*$ and $r_{k,t}^{\star}=0$
for all cases, except for $w(f_t)=m_{2j+1}<m_{2j_0+1}$. In this case, we observe that since $t$ is
splitting there are at least two immediate successor $s_1\neq s_2\in S_t$ such that $\ran
w_{k,u}^\star\cap \ran f_{s_i}\neq \buit$ ($i=1,2$) for some $u\in \mc B_k^\star$. This implies
that there is at most one $k\in E_t^{\mc B^\star}$ such that $\ran f_{s(t)}\cap \ran
w_{k,v}^\star\neq \buit$ for $v\in \mc B_k^\mm{unc}$, and $t\in \mc B_{k,v}^\star$. Then we set
$g_{k,t}^{\star}=(1/m_{2j_0+1})e_k^*$ and $r_{k,t}^{\star }=0 $ if $k$ is this one, and
$g_{k,t}^{\star }=0$ and $r_{k,t}^{\star }= (1/m_{2j_0+1})e_k^*$ otherwise.

Let $(G_t^\star)_{t\in \mc T}$ , $(R_t^\star)_{t\in \mc T}$ be the corresponding filtrations.
Recall that given a regular array $\boldsymbol{\mc A}=(\mc A_k)_k$  for $(x_k)_k$ and $(f_t)_{t\in
\mc T}$ the canonical $\boldsymbol{\mc A}$-assignment $(f_{k,t}^{\boldsymbol{\mc A}})_{t\in \mc
A_k,k}$ is defined by $f_{k,t}^{\boldsymbol{\mc A}}=f(x_k)e_k^*$. It was shown in Remark
\ref{canonassig} that if in addition $\boldsymbol{\mc A}$ is maximal, then for every
$(a_k)_{k=1}^n$ and every $t\in \mc T$, $F^{\boldsymbol{\mc A}}(\sum_{k\in D_t^{\boldsymbol{\mc
A}}}a_k e_k)=f_t(\sum_{k\in D_t^{\boldsymbol{\mc A}}}a_k w_k)$.

\clam  Fix $t\in \mc T$, and for $\star\in\{\mm{in},\mm{fin},\mm{cnd},\mm{at}\}$ let
$D_t^\star=D_t^{{\boldsymbol{\mc B^\star}}}$. Then:

\noindent{(e.1)} $|F_t^{\star}(\sum_{k\in D_t^{\star}}b_ke_k) |\le 24
(G_t^{\star}+R_t^{\star})(\sum_{k\in D_t^{\star}}|b_k|e_k)$ for $\star\in\{\mm{in},\mm{fin}\}$.

\noindent{(e.2)} $|F_t^{\mm{cnd}}(\sum_{k\in D_t^{\mm{cnd}}}b_ke_k)
|\le 24 (G_t^{\mm{cnd}}+R_t^{\mm{cnd}})(\sum_{k\in D_t^{\mm{cnd}}}b_ke_k)$.

\noindent{(e.3)} $|F_t^{\mm{at}}(\sum_{k\in D_t^\mm{at} })b_ke_k|\le 24 |R_t^{\mm{at}}(\sum_{k\in
D_t^\mm{at}})b_ke_k|$.

\noindent{(e.4)}  $G_t^{\star}\in W(T_0)$ for $\star\in\{\mm{in},\mm{fin}\}$, and
$G_t^{\mm{cnd}}\in 3W(J_{T_0})$.

\noindent{(e.5)} $\nrm{R_{t}^{\mm{at}}}_{\infty}\le 1/m_{2j_0+1}$. For  $\star\in \{\mm{in},\mm{fin},\mm{cnd}\}$,
  either $t$ is non-relevant, $w(f_t)<m_{2j_0+1}$
and  $G_k^{\star}=\sum_{k\in E_t^{\mc B^\star}}\nrm{f_t(w_k)}e_k^*$ or $\nrm{R_{t}^{\star}}_{\infty}\le
1/m_{2j_0+1}$.

\fclam

\prucl \emph{(e.1)-(e.3)} are immediate applications of Proposition \ref{desigualdades}.
\emph{(e.4)}: Most of the cases follow immediately by definition of the corresponding
assignments. We sketch the non-trivial ones: Suppose that $t$ is relevant. If
$w(f_t)=m_{2j_0+1}$, then $D_t^{ \mm{cnd}}=E_t^{ \mm{cnd}}$, and the corresponding assignment
gives that $G_t^{\mm{cnd}}=\la_1 e_{k(t,1)}^* + \la_2 e_{k(t,2)}^*+\vep_t\sum_{k\in E_t^{
\mm{cnd}}\cap ( k(t,1),k(t,2))} e_{k}^*\in 3W(J_{T_0})$, where
 $\la_i=\sgn(b_{k(t,i)})\chi_{E_t^{ {\mm{cnd}}}}(k(t,i))$, for $i=1,2$, and where $\chi_E$ denotes the characteristic
 function of $E$. Fix $\star\in\{\mm{in},\mm{fin}\}$. We claim that $\# D_t^{ \star}\le 1$: Suppose not, and say
that $k_1<k_2\in D_t^{ \star}$. Then since $w(f_t)=m_{2j_0+1}$ there are $u_i\in \mc
B_{k_i}^\mm{unc}$ and  $s_i\in \mc B_{k_i,u_i}^\star$ ($i=1,2$) such that $u_1,u_2 \precneqq t
\preceq s_1,s_2$. If $\star=\mm{in}$, then since $\ran f_t\con \ran f_{s(k_1,u_1)}$, it follows
that $\ran f_t<\ran w_{k_2}$, and since $\ran f_{s_2}\con \ran f_{t}$, we obtain that $\ran
f_{s_2}\cap \ran w_{k_2}=\buit$, contradicting  the fact that $s_2\in \mc B_{k_2,u_2}^\mm{in}$. If
$\star=\mm{fin}$, in a similar manner we obtain that $\ran w_{k_1}\cap \ran f_{s_1}=\buit$, a
contradiction. Hence, either $G_{t}^{\star}=0$, or $ G_{t}^{\star}=(1/m_{2j_0+1})e_k^*$, certainly
in $W(T_0)$ considered as sub-convex combinations.

Suppose now that $w(f_t)\neq m_{2j_0+1}$. There are three subcases to consider: If
$w(f_t)>m_{2j_0+1}$, then $t$ is non-relevant, hence a catcher node, and   $G_t^{\star}=\sum_{k\in
E_t^\star}g_{k,t}^\star=0$. If    $w(f_t)=m_{2j} $ with $j\le j_0$, then the inductive hypothesis
gives that $G_t^{\mm{cnd}}\in 3W(J_{T_0})$ (since $E_t^{ \mm{cnd}}=\buit$). Fix
$\star\in\{\mm{in},\mm{fin}\}$. Observe that for every $k\in E_t^{ \star}$ there is $s\in S_t$
such that $\ran f_s\con \ran w_k$, in which case $D_s^{ \star}=\buit$, and so $\#E_t^{
\star}+\#\conj{s\in S_t}{G_s^{\mm{cnd}}\neq 0}\le \# S_t$, and then,
$G_t^{\star}=(1/w(f_t))(\sum_{k\in E_t^{ \star}}e_k^*+\sum_{s\in S_t} G_s^\star )\in W(T_0)$.

If $w(f_t)=m_{2j+1}<m_{2j_0+1}$, then using that there is at most one immediate successor $s(t)$
of $t$ which is relevant we obtain that either $G_{t}^{\mm{cnd} }=0$, or $G_{t}^{\mm{cnd}
}=(1/m_{2j})G_{s(t)}^{\mm{cnd}}$, and for  $\star\in\{{\mm{in},\mm{fin}}\}$, either
$G_{t}^{\star}=(1/m_{2j+1})e_k^*$, or $G_{t}^{ \star}=(1/m_{2j+1})G_{s(t)}^{\star}$.

\emph{(e.5)}:    $\nrm{R_t^{\mm{at}}}_{\infty}\le 1/m_{2j_0+1}$ follows from Proposition
\ref{desigualdades}, since  this is so for the corresponding assignment of which $R_t^\mm{at}$
is a filtration. Suppose that $\star\in \{\mm{in},\mm{fin},\mm{cnd}\}$. The proof is by
backwards induction over $t$.  Again we concentrate in non-trivial cases. Suppose that $f_t$
is of type I and $t$ is  relevant. Then if $w(f_t)=m_{2j}$ with $j\le j_0$, then the desired
result follows from the definition of the corresponding assignments, and inductive hypothesis.
Suppose that $w(f_t)= m_{2j+1}$ with $j<j_0$. Then $R_{t}^{\star}=\sum_{k\in E_t^{
\star}}r_{k,t}^{\star}+(1/w(f_t))\sum_{s\in S_t}R_s^{\star}$. By the definition of the
assignments, $\nrm{r_{k,t}^{\star}}_{\infty}\le 1/m_{2j_0+1}$ for every $k\in E_t^{\star}$.
Observe that all $s\in S_t$, except possibly one, $s(t)$, are  non-relevant and that
$r_{k,s}^{\star}=\nrm{f_s(w_k)}e_{k}^*$ for every $k\in E_s^{\star}=D_s^{\star}$. Hence, for
every $s\in S_t\setminus \{s(t)\}$,
$\nrm{(1/w(f_t))R_{s}^{\star}}_\infty=\max\conj{(1/w(f_t))\nrm{f_s(w_k)} }{k\in
E_s^{\star}}\le 1/m_{2j_+1}$; the last inequality follows from Proposition \ref{estonw}. By
the inductive hypothesis $\nrm{R_{s(t)}^{\star}}_\infty\le 1/m_{2j_0+1}$, so we are done.

Suppose that $t$ is non-relevant. The case $w(f_t)>m_{2j_0+1}$ is immediate. Suppose that
$w(f_t)=m_{2j}$ and $t\neq s(u)$, where $u$ is the immediate predecessor of $t$ (see the
definition of relevant node). Notice that $t$ is a catcher, so $E_t^{ \star}=D_t^{ \star}$, and
$R_{t}^{\star}=\sum_{k\in E_t^{ \star}} \nrm{f_t(w_k)}e_k^*$, as desired. \fprucl

We are now ready to finish the proof using the part \ref{twofiltra} of the general theory above.
Notice that for each $\star\in \{\mm{in},\mm{fin},\mm{cnd},\mm{at}\}$, ${\boldsymbol{\mc
B^\star}}\nprec {\boldsymbol{\mc B}}$, and that the canonical assignments of $\boldsymbol{\mc
B^\star}$ are coherent. Let $(h^\star_{k,t})_{t\in \mc B_k,k}$ be the assignments induced by the
canonical ${\boldsymbol{\mc B^\star}}$-assignments $(f^\star_{k,t})_{t\in \mc B^\star_k,k}$, for
$\star\in\{\mm{in},\mm{fin},\mm{cnd},\mm{at}\}$.

\clam For every $t\in \mc T$,
$H_t^{\mm{in}}+H_t^{\mm{fin}}+F_t^{\mm{cnd}}+F_t^{\mm{at}}=F_{t}^{{\boldsymbol{\mc B}}}$, the
canonical assignment of ${{\boldsymbol{\mc B}}}$.\fclam \prucl We show that for every $t\in B_k$,
$h_{k,t}^{\mm{in}}+h_{k,t}^{\mm{fin}}+h_{k,t}^{\mm{cnd}}+h_{k,t}^{\mm{at}}=f_t(w_k)e_{k}^*$. The
only non trivial case is if $t\in \mc B_k^\mm{unc}$. Notice that since $B_{k,t}^\star$ is a
maximal antichain for $w_{k,t}^\star$ and $(f_s)_{s \succeq t}$, we obtain that
$h_{k,t}^{\star}=F_{k,t}^{\star}=f_t(w_{k,t}^\star)e_k^*$. Hence,
$f_{k,t}^{\mm{in}}+f_{k,t}^{\mm{fin}}=(f_{t}(w_{k,t}^\mm{in})+f_{t}(w_{k,t}^\mm{fin}))e_k^*=f_t(w_k)e_k^*$,
and $h_{k,t}^\star=0$ for $\star\in\{\mm{cnd},\mm{at}\}$.  \fprucl

Finally, by Proposition \ref{iurppp}, $H_{\emptyset}^{\star}=F_{\emptyset}^{\star}$, for $\star\in
\{\mm{in},\mm{fin},\mm{cnd},\mm{at}\}$. Hence,
\begin{align}  \label{otherinequ}
|f(\sum_{k=1}^n b_k w_k )|=&  |F_{\buit}^{{{\boldsymbol{\mc B}}}}(\sum_{k=1}^n b_k e_k
)|\le \sum_{\star\in
\{\mm{in},\mm{fin},\mm{cnd},\mm{at}\}}|F_{\buit}^{\star }(\sum_{k\in D_\buit^{\star}} b_k e_k) |\le \nonumber\\
\le &    24(5\nrm{ \sum_{k=1}^nb_k e_k}_{J_{T_0}}+4\nrm{\sum_{k=1}^n b_k
e_k}_{\infty})\le 120\nrm{ \sum_{k=1}^n e_k}_{J_{T_0}}+\vep.
\end{align}
\fprue

 \cor\label{mayassme} The natural isomorphism $F:\langle
w_1,\dots,w_n \rangle\to \langle v_1,\dots,v_n \rangle$ defined by $F(w_i)=v_i$ satisfies that
$\nrm{F}\le 1$ and $\nrm{F^{-1}}\le 120+\vep$. Consequently, $J_{T_0}$ is finite interval
representable on the basis $(e_\al)_{\al<\ou}$ of $\eqs_{\ou}$ with a constant $\boldsymbol{C}<
121$. \fcor
\prue
Proposition \ref{lkio1} shows that $\nrm{F}\le 1$;  the other inequality follows from Lemma
\ref{theotherin}.
\fprue

\section{The unconditional counterpart} We produce a space $\eqs_{\omega_1}^u$ which is the counterpart of
$\eqs_{\omega_1}$ in the frame
 of the spaces with an unconditional basis, as in \cite{G2}. This space is defined  as was
   $\eqs_{\omega_1}$ by a
 norming family of functionals $K_{\omega_1}^u$  satisfying \emph{(1)}-\emph{(4)} from Subsection
 \ref{definitionnorming}, and in addition the following condition
 \begin{enumerate}
 \item[\emph{(5)}]  It is closed under the restriction of all functionals with odd weight  to every subset of $\omega_1$
 \end{enumerate}
Although $\eqs_{\omega_1}^u $ belongs to the class of  spaces with an unconditional basis its
study uses the same tools used in the study of $\eqs_{\omega_1}$. For example,  given a
bounded operator $T:\eqs_{\omega_1}^u\to \eqs_{\omega_1}^u$ the transfinite sequence
$(d(Te_{\gamma},\R e_{\gamma}))_{\gamma<\ou}$ belongs to $c_0(\ou)$, and the operator $T$ is
strictly singular if and only if the sequence $(\|Te_{\gamma}\|)_{\gamma<\ou}$ belongs to
$c_0(\ou)$.
\nota\label{somethatremains}
\noindent 1. The basic inequality (Lemma \ref{bin}) still remains true provided that
(\ref{negligible}) holds for an arbitrary subset $E\con [1,n]$, not only for intervals.

\noindent 2. For every block sequence $(y_n)_n$ of $\eqs_{\ou}^u$ and every $j$ there is a
$(6,j)$-exact pair $(y,\phi)$ with $y\in \langle y_n\rangle_n$ (indeed, what one locates first
are $2-\ell_1^n$ averages.)
\fnota

The next result is the corresponding analogue from \cite{go-ma2}.
\prop\label{gowerstrick}
Let $T:\eqs_{\ou}^u\to \eqs_{\ou}^u$ be bounded, and let $(x_n)_n$ be a RIS of $\eqs_{\ou}^u$. For
each $n$, let $B_n\cup C_n=\supp x_n$ be a partition. Then $\lim_{n\to \infty} C_n T B_n x_n=0$.
\fprop
\prue(Sketch)  Assume not. Notice that since $(x_n)_n$ is a block sequence, so is $(C_n TB_n x_n)_n$. Going
to a subsequence if needed we assume that $\inf_n \nrm{C_n TB_nx_n}\ge \vep >0$. Since for every $\phi\in
K_{u,\ou}$, the restriction $A\phi\in K_{u,\ou}$ for every subset $E\con \ou$, we have that the sequence
$(B_n x_n)_n$ is also RIS. Now for each $n$, choose $f_n\in K_{u,\ou}$ such that $\supp f_n\con C_n$ and $f_n
(C_nTB_nx_n)\ge \vep$. Let  $j$ be such that $\nrm{T}<m_{2j+1}\vep$, and find appropriate
$(2j_i)_{i=1}^{n_{2j+1}}$ such that
\begin{equation}   (\frac{m_{2j_1}}{n_{2j_1}}\sum_{k\in F_1}B_k x_k,
\frac1{m_{2j_1}}\sum_{k\in F_1}f_k,\dots., \frac{m_{2j_{n_{2j+1}}}}{n_{2j_{n_{2j+1}}}}\sum_{k\in
F_{n_{2j+1}}}B_k x_k, \frac1{m_{2j_{n_{2j+1}}}}\sum_{k\in F_{n_{2j+1}}}f_k)
\end{equation}  is a $(0,j)$-dependent sequence, for $F_1<\dots <F_{n_{2j+1}}$, each $\#F_i=n_{2j_i}$.
Then, $\nrm{T x}\ge {\vep}/{m_{2j+1}}>\nrm{T}\nrm{x}$ where $x=1/n_{2j+1}\sum_{i=1}^{n_{2j+1}}
(m_{2j_i}/n_{2j_i}\sum_{k\in F_i}B_k x_k)$, a contradiction.
\fprue

\prop\label{unc1}
Let $T:\eqs_{\ou}^u\to \eqs_{\ou}^u$ be bounded such that for all $\al<\ou$,
$e_{\al}^*Te_{\al}=0$. Then  $\lim_{n\to \infty }T x_n=0$ for every RIS $(x_n)_n$.
\fprop
\prue
For each $n$, let $A_n=\supp x_n$.
\clam
$\lim_{n\to \infty}A_nT x_n=0$.
\fclam
\prucl
Notice that
\begin{equation}
A_n Tx_n=\left\{\begin{array}{ll} \frac{2 L_n(2L_n-1)}{L_n^2}  \frac1{\#
P_n}\sum_{(B,C)\in P_n} B T C x_n &
\text{if } \# A_n \text{ even}\\
\frac{2L_n (2L_n+1) ((L_n+1)^2+1)}{(L_n+1)^2 (L_n^2+1)} \frac1{\# P_n}\sum_{(B,C)\in P_n}
B T C x_n & \text{if } \# A_n \text{ odd}
\end{array}\right.
\end{equation}
where $ L_n$ is the entire part of $\#A_n /2$, and
\begin{equation}
P_n=\left\{\begin{array}{ll} \conj{(B,C)}{B\cup C=A_n,\, B\cap C=\buit, \, \#B =\supp x_n /2 } &
\text{if } \#A_n \text{ even}\\
\conj{(B,C)}{B\cup C=A_n,\, B\cap C=\buit, \, |\#B -\#C|=1  } & \text{if } \#A_n \text{
odd}
\end{array}\right.
\end{equation}
Hence, $A_n Tx_n =(\la_n/{\# P_n})\sum_{(B,C)\in P_n} B T C x_n$ with $1\le \la_n \le 4$. By
Proposition \ref{gowerstrick}, $A_nTx_n\to_n 0$, as desired. \fprucl Now suppose that
$\lim_{n\to \infty}Tx_n\neq 0$. W.l.o.g. we may assume that $(Tx_n)_n$ is a block sequence and
with support disjoint from $(x_n)_n$ (let $\ga_0$ be the minimal $\ga<\ou$ such that there is
some infinite $A$ such that $\inf_{n\in A}\nrm{P_\ga T x_n}>0$; now, replacing $T$ by
$P_{\ga_0} T$, and going to a subsequence $(x_n)_{n\in A}$ we may assume that $(Tx_n)_n$ is a
block sequence. By the previous Claim we obtain that $A_nTx_n\to_n 0$, so we may assume that
$(Tx_n)_n$ and $(x_n)_n$ are disjointly supported). Now it is easy to produce, for large
enough $j$, a $(0,j)$-dependent sequence $(y_1,\phi_1,\dots,y_{n_{2j+1}},\phi_{n_{2j+1}})$
such that $ \nrm{T ((1/n_{2j+1})\sum_i y_i)}>\nrm{T}\nrm{(1/n_{2j+1})\sum_i y_i}$, a
contradiction. \fprue In the same way one can show the following useful result.
\prop\label{preblockev} For every $X\hookrightarrow \eqs_{\ou}^u$ generated by a block
sequence $(x_n)_n\con \eqs_{\ga}^u$, every bounded $T:X\to \eqs_{[\ga,\ou)}^u$ is strictly
singular. Indeed, $\lim_{n\to \infty} Ty_n=0 $ for every RIS $(y_n)_n$ in $X$. \qed \fprop

\cor\label{blocksverunc} For every $X\hookrightarrow \eqs_{\ou}^u$ and every $\vep>0$,
there is some block sequence $(z_n)_n$ of $\eqs_{\ou}^u$ and some Schauder basis
$(x_n)_n\con X$ such that $\nrm{z_n-x_n}\le \vep$.
\fcor
\prue Fix $X\hookrightarrow \eqs_{\ou}^u$. By standard facts of transfinite block sequences
(see Proposition \ref{gliding}) we can find some $\la<\ou$  a block sequence $(w_n)_n$ of
$\eqs_{\la}^u$, and a sequence $(y_n)_n\con X$ such that $\sum_n\nrm{P_\la y_n -z_n}\le
\vep/2$. W.l.o.g.  (going to a block subsequence if needed) we may assume that $(z_n)_n$ is a
RIS. Consider $U:\overline{\langle w_n\rangle}_n\to \eqs_{\ou}^u$ defined by
$U(w_n)=P_{[\la,\ou)}x_n$. Since $P_\la| \overline{\langle y_n\rangle_n}$ is an isomorphism,
$U$ is bounded. By Proposition \ref{preblockev}, $U$ is strictly singular. Hence we can find a
block subsequence $(z_n)_n$ of $(w_n)_n$ and the corresponding block subsequence $(x_n)_n$ of
$(y_n)_n$ such that $\nrm{z_n-x_n}\le \vep$. \fprue
\cor\label{0ss} If $T:\eqs^u_{\ou}\to \eqs^u_{\ou}$ is bounded and for all $\al$ we have
that $e_\al^* T e_{\al}=0$, then $T$ is strictly singular. \fcor \prue Let $X
\hookrightarrow \eqs_{\ou}^u$, and fix $\vep>0$. Choose  some RIS $(z_n)_n$ and some
sequence $(x_n)_n\con X$ such that $\sum_n \nrm{z_n-x_n}\le \vep/\nrm{T}$. By Proposition
\ref{unc1}, $\lim_{n\to \infty} Tx_n=0$. Hence we can find $x\in \langle x_n\rangle_n$
such that $\nrm{Tx}\le \vep$. \fprue

\cor
For any $T:\eqs^u_{\ou}\to \eqs^u_{\ou}$, there is some diagonal operator $D_T$ such that $S=T-D_T$ is
strictly singular, $S e_\al=0$ for all $\al<\ou$ and $S$ has separable range. \fcor
\prue
Let $D_T:\eqs_{\ou}^u\to \eqs_{\ou}^u$ be defined  for $\al<\ou$ by $D_T(e_\al)=e^*_{\al}(T e_\al)
e_\al$. $D_T$ is bounded and by Corollary \ref{0ss},  $T-D_T$ is strictly singular.
\fprue

\cor
For any infinite $A\con \ou$, the space $\eqs_{A}^u$ is reflexive with an
unconditional basis
and
\begin{equation}
\mc{L}(\eqs^u_A)\cong \mc{D}(\eqs^u_A)\oplus \mc{S}(\eqs^u_A).
\end{equation}
\qed \fcor Here $\mc{D}(\eqs^u_A)$ denotes the space of the diagonal operators and
$\mc{S}(\eqs^u_A)$ is the space of strictly singular operators $S$ with separable range
such that $e_{\al}^*(Se_{\al})=0$ for every $\al\in A$.

\cor
For any infinite $A\con \ou$, $\eqs^u_A$ is not isomorphic to a proper subspace of itself.
\fcor
 \prue Let $X \hookrightarrow \eqs^u_A$, $T: \eqs^u_A \to X$ be an isomorphism and let
$U= i_{X,\eqs^u_A}\circ T$ be a semi-Fredholm operator with $\al(U)=0$. Then $U=D_U+S$, $D_U$
diagonal such that $D_T(e_\al)=e_\al^* (T e_\al)e_\al$, and $S$ strictly singular. Since $D_U$
is a strictly singular perturbation of the semi-Fredholm operator $U$ with $\al(U)=0$, $D_U$
is semi-Fredholm, and $\al(D_U)<\infty$. But $\ker D_U=\overline{\langle
\conj{e_\al}{Te_\al=0} \rangle}$. So, $D_U \eqs^u_A=\overline{\langle \conj{e_\al}{Te_\al\neq
0} \rangle}$ which has co-dimension equal to $\al(D_U)$, hence  $D_U$  and  $U$ are Fredholm
with index 0. Since $U$ is 1-1, this implies that $X=\eqs^u_A$, as desired. \fprue

\cor
Let $A,B$ two infinite sets of countable ordinals such that $A\cap B$ is finite. Then
every bounded operator $T:\eqs^u_A\to \eqs^u_B$ is strictly singular.\qed \fcor

\cor There is a nonseparable reflexive space $X$ with an unconditional basis such that

\noindent(a) $X$ is not isomorphic to any of its proper subspaces.

\noindent (b) Every bounded linear operator $T:X\to X$ is of the form $D+S$ with $D$ a
diagonal operator and $S$ a strictly singular operator with separable range.

\noindent(c) For every $I_1$, $I_2$ infinite disjoint subsets of $\omega_1$ the spaces
$\eqs_{I_1}$, $\eqs_{I_2}$ are totally incomparable. \fcor

Suppose now that in addition $\ro$ is universal.
\cor
For every interval $I$ of ordinals, $(e_{\al})_{\al\in I}$ is nearly subsymmetric.
Moreover, for any two minimal intervals $I=[\al,\al+\om)$, $J=[\be,\be+\om)$,
$\eqs^u_{I}$ is an asymptotic version of $\eqs ^u_J$.  \fcor

So, if we consider the version of $\eqs^u_{\ou}$ obtained  by a
universal $\ro$-function then the unconditional basis
$(e_{\al})_{\al<\ou}$ is nearly subsymmetric and for any pair of
disjoint minimal infinite intervals $I_1$, $I_2$ $\eqs^u_{I_1}$ is
an asymptotic version of $\eqs^u_{I_2}$, while they are totally
incomparable.

\prop
The unconditional counterpart $\eqs_{\ou}^u$ is arbitrarily distortable.
\fprop
\prue
The norms $(\nrm{\cdot}_{u,j})_j$ arbitrarily distort the space $\eqs_{\ou}^u$, since $(6,j)$-exact pairs exist in
every block sequence and by Corollary \ref{blocksverunc} every subspace $X\hookrightarrow \eqs_{\ou}^u$
``almost" contains a block sequence.
\fprue

\end{document}